%% file: ETD-TDSR.tex
\documentclass[onefignum,onetabnum]{siamart190516}
\input{ex_shared}

\ifpdf
\hypersetup{
  pdftitle={Energy-dissipative spectral renormalization exponential integrator method for gradient flow problems},
  pdfauthor={Dianming Hou, Lili Ju, and Zhonghua Qiao}
}
\fi

\begin{document}
\graphicspath{{figures/},}
\maketitle

\begin{abstract}
 In this paper, we present a novel spectral renormalization exponential integrator method for solving gradient flow problems. Our method is specifically designed to simultaneously satisfy  discrete analogues of the energy dissipation laws and achieve high-order accuracy in time. To accomplish this, our method first incorporates the energy dissipation law into the target gradient flow equation by introducing a time-dependent spectral renormalization (TDSR) factor. Then, the coupled equations are discretized using the spectral approximation in space and the exponential time differencing (ETD) in time. Finally, the resulting fully discrete nonlinear system is decoupled and solved using the Picard iteration at each time step. Furthermore, we introduce an extra enforcing term into the system for updating the TDSR factor, which greatly relaxes the time step size restriction of the proposed method and enhances its computational efficiency. Extensive numerical tests  with various gradient flows  are also presented  to demonstrate the accuracy and effectiveness of our method as well as its high efficiency when combined with an adaptive time-stepping strategy for long-term simulations.
\end{abstract}



\section{Introduction}
\setcounter{equation}{0}

Partial differential equations (PDEs) with gradient flow structure are commonly used to model physical phenomena in many scientific and engineering fields, including materials science \cite{All79,Cahn59,Eld02}  and fluid dynamics \cite{And97,Gur96}. These PDE systems are derived from the energy variational principle of total free energy in different Sobolev spaces. As a result, gradient flow models typically take the following general form:
\be\label{prob0}
\dps\frac{\partial \phi}{\partial t}=-\mbox{grad}_{H} E[\phi],\quad \x\in\Omega, t> 0,
\ee
where $\phi(\x,t)$ denotes the scalar-valued phase function defined over a domain $\Omega\subseteq \mathbb{R}^{d}$\ ($d =1,2,3$) at time $t$, $E[\phi]$  is the free energy functional associated with the corresponding physical problem, and
$\mbox{grad}_{H} E[\phi]$ is the functional derivative of $E$ with respect to $\phi$ in the Sobolev space $H$.
This gradient flow model satisfies the energy dissipation law:
 \bq\label{EDlaw}
 \frac{d}{dt}E[\phi]
 =\Big(\mbox{grad}_{H} E[\phi], \frac{\partial \phi}{\partial t}\Big)_{H}
 =-\|\mbox{grad}_{H} E[\phi]\|_{H}^{2},
 \eq
 where $(\cdot,\cdot)_{H}$ and $\|\cdot\|_H$ are the inner product and associated norm of the Sobolev space $H$, respectively.
This implies that the phase solution $\phi$ evolves in a way that decreases the energy functional $E[\phi]$ over time. 
For evolution gradient-flows-structure equations arising in physical applications, it is important to design numerical schemes that can enforce the corresponding physics into simulations, particularly for the energy dissipation law \eqref{EDlaw}.
Many efforts have been devoted to developing energy stable numerical methods in the long-standing and active research field of gradient flows. These methods include, but are not limited to, the convex splitting method \cite{BB97,BLWW13,Ell93,Eyr98}, the linear stabilization method \cite{JLQZ18,LQT16,SY10,Xu06}, the Lagrangian multiplier method \cite{Bad11,GT13}, the Invariant Energy Quadratization (IEQ) method \cite{Yan16,Zha17}, and the Scalar Auxiliary Variable (SAV) method \cite{SX18,Shen17_1}.

The convex splitting method was first introduced by Elliott and Stuart in \cite{Ell93} to numerically study a class of semilinear parabolic  equations. The main idea of the convex splitting method for gradient flow models is to split the free energy functionals into two parts: the convex part and the concave one. The convex and concave parts are then treated implicitly and explicitly, respectively, to derive unconditionally energy stable or energy dissipative numerical schemes, as seen in \cite{BLWW13,Eyr98,GLWW14} and the references cited therein. Although first- and second-order unconditionally stable convex splitting schemes have been obtained for some specific energy functionals of gradient flows, there is no general convex-concave splitting form of the energy functional that allows for high-order time-stepping schemes that are unconditionally energy stable.  Another widely-used approach is the linear stabilization method, which adds one or two linear stabilization terms into the corresponding time-stepping schemes based on backward differentiation formulation (BDF), Crank-Nicolson, or exponential time differencing (ETD) methods to greatly relax the time step size constraints. In this type of method, the nonlinear term is usually treated explicitly, leading to a linear system solved at each time step. Although unconditional or conditional stability of the linear stabilization method for some gradient flows  have been successfully established (see e.g., \cite{DJLQ18,LQT16,LQW21,LQW22}), there does not exist a general framework for stability analysis of these existing schemes, particularly for high-order ones.
Recently, the SAV approach \cite{Shen17_1} and its variants \cite{ALL19,HAX19,JZZ22} have been developed to design unconditionally energy stable linear schemes for gradient flow problems. In particular, SAV schemes with high-order accuracy and energy stability have been of great interest, due to the long-time nature of the gradient flow coarsening process. In \cite{ALL19}, an arbitrarily high-order extrapolated and linearized Runge--Kutta SAV (RK--SAV) method was constructed for the Allen--Cahn and Cahn--Hilliard  equations. The proposed RK-SAV method is unconditionally stable with respect to a modified discrete energy, and the corresponding error estimate was also rigorously derived. Huang et al. \cite{HSY21} introduced a novel SAV approach to construct an implicit-explicit linear and unconditionally energy stable BDF$k$ $(1 \leq k \leq 5)$ method for gradient flows. However, its unconditional stability only indicates dissipation in time of the auxiliary variable without information on the phase variables.

The steady-state spectral renormalization method was introduced by Ablowitz and Musslimani for the first time in \cite{AM05} to compute self-localized states of nonlinear waveguides. Since then, it has been utilized in a variety of contexts, such as nonlinear optics \cite{AABI10,AABI12}, Bose-Einstein condensation \cite{HI12}, and water waves \cite{AFM06}. Built on this idea, a type of time-dependent spectral renormalization (TDSR) approach was developed in \cite{CFM22,CM17} to incorporate intrinsic physics in the form of conservation laws or dissipation rate equations in the development of numerical schemes. This approach has achieved the impressive capacity for accurately and efficiently capturing certain relevant physical properties regardless of the accuracy of the solution in \cite{CM17}, particularly preserving certain conservation laws of the model. However, it also encounters convergence issues when solving the nonlinear systems of the  physical evolution equation and the  ordinary differential equation involving the renormalization factors with some iterative methods. This leads to the efficiency problem of time step sizes not being large in order to ensure the convergence of the iteration. Such convergence issue often becomes particularly severe for cases where the dissipative properties (such as the energy) of the model are enforced in the TDSR approach. Furthermore, for the dissipative model, the $p$-th order TDSR scheme is numerically stable only in the sense that the error of the computed energy is of order $p$ in time, which doesn't indicate the dissipation property of the model at the discrete level.

A common feature of the gradient flow problem \eqref{prob0} is that its evolution process often takes a considerably long time to reach the steady state. Moreover,  it usually undergoes both fast and slow energy-changing stages throughout the evolution process, which imposes the need for adaptive time stepping during the simulation. Therefore, there is a high demand for the development of high-order, structure-preserving, efficient numerical schemes with variable time steps for the gradient flow \eqref{prob0}.  To address these issues, we propose in
this paper a high-order accurate and energy dissipative method, called TDSR-ETD, for solving general gradient flow problems, which combines the spectral renormalization method for handling the energy dissipation law and the exponential integrator methods for accurate and stable time integration. As motivated by the idea of \cite{CFM22,CM17}, we further apply the proposed  TDSR-ETD  method to a broader class of phase-field models with energy dissipation properties.

The rest of the paper is organized as follows. Section 2 provides a detailed illustration of the proposed TDSR-ETD method for the $L^2$ (Allen-Cahn type) and $H^{-1}$  (Cahn-Hilliard type) gradient flows with respect to a classic free energy functional, under the periodic or homogeneous Neumann boundary conditions. This section includes discussions of the corresponding numerical schemes, solution algorithms, and physical properties, such as energy dissipation and mass conservation. In Section 3, we further explore the application of our TDSR-ETD method to other two types of gradient flow problems, the molecular beam epitaxial model and the phase-field crystal model. Section 4 presents extensive numerical experiments and comparison tests to demonstrate the accuracy and efficiency of the proposed method. Finally, some concluding remarks are given in Section 5.

\section {The spectral renormalization exponential integrator method}\label{sect2}
\setcounter{equation}{0}
To illustrate the proposed spectral renormalization exponential integrator method, we  take the gradient flow model \eqref{prob0} with respect to the following classic free energy functional:
\be\label{EF}
E[\phi]=\int_{\Omega}\Big[\frac{\varepsilon^{2}}{2}|\nabla\phi|^{2}+F(\phi)\Big]d\x,
\ee
where  the parameter $\varepsilon>0$ is related to the interfacial width and  $F: \mathbb{R}\rightarrow  \mathbb{R}$ denotes a nonlinear potential function. Two types of boundary conditions usually will be considered for the above gradient flows, the periodic boundary condition and the homogenous Neumann boundary condition, respectively.
Taking the Sobolev space $H$ to be
$L^2(\Omega)$ or $H^{-1}(\Omega)$ in \eqref{prob0}  gives us the following time-dependent PDEs:
\brr\label{prob}
\begin{cases}
\dps\frac{\partial \phi}{\partial t}
=\mathcal{G}_{H} \mu,\qquad\qquad\quad\;\;  \x\in\Omega, t> 0,\\
\dps\mu=-\varepsilon^{2}\Delta\phi+f(\phi),\qquad \x\in\Omega, t> 0,
\end{cases}
\err
with the initial value  $\phi(\x,0)=\phi_{0}(\x)$ for any $\x\in\overline\Omega$, where $f=F'$ and
\bq\label{ef}
\mathcal{G}_H:=
\begin{cases}
\begin{array}{ll}
-I,&\quad \mbox{if $H:=L^{2}(\Omega)$},\\
\Delta,&\quad \mbox{if $H:=H^{-1}(\Omega)$.}
\end{array}
\end{cases}
 \eq
We call the above equations \eqref{prob}  the $L^2(\Omega)$ or  $H^{-1}(\Omega)$  gradient flows with respect to the energy functional \eqref{EF}, respectively (equivalently, the Allen-Cahn type  equation \cite{All79}  or the Cahn-Hilliard type equation \cite{Cahn59}, respectively).

 Inspired by the work of \cite{CFM22,CM17},
the key idea of the proposed spectral renormalization exponential integrator method is to introduce an extra scalar variable, that is the so-called TDSR factor $R(t)$, to incorporate the energy dissipation law \eqref{EDlaw} into  the gradient flow problem \eqref{prob}. In particular, we rewrite the solution of the gradient flow problem \eqref{prob}
as
\begin{equation}\label{ttt}
\phi(\x,t)=R(t)\psi(\x,t),
\end{equation}
and thus $\mu\big(R\psi\big)=-\varepsilon^{2}\Delta(R\psi)+f(R\psi)$.
Without loss of generality, we focus our following discussion  on
the two-dimensional case $(d=2,\;\x=(x,y))$ and assume $\Omega:=(0,L)^{2}$, but all results derived below can be straightforwardly  extended to the three-dimensional case. Let $T>0$ be a given terminal time and $\{\Dt_{n}=t_{n}-t_{n-1}>0\}_{n=1}^K$ be a general partition of the time interval $[0,T]$ such that $t_0=0$ and $\sum_{n=1}^{K} \Dt_n=T$.

\subsection{The $L^2$ gradient flow with the periodic boundary condition}\label{scase1}
\label{case_I}
Based on the $L^2$ gradient flow  equation \eqref{prob} (i.e., $\mathcal{G}_H = -I$) and the energy dissipation law \eqref{EDlaw}, we have the following coupled system for
$(R(t),\psi(\x,t))$:
\begin{subequations}\label{prob1}
 \begin{align}
&\dps\frac{\partial (R\psi)}{\partial t}=\varepsilon^{2}\Delta(R\psi)-f(R\psi),\qquad (\x,t)\in \Omega\times(0,T],\label{111}\\
&\dps\frac{dE\big[R\psi\big]}{dt}=-\|\mu\big(R\psi\big)\|^{2},\qquad\quad \quad t\in(0,T],\label{prob1-b}
 \end{align}
\end{subequations}
with the initial value  $\dps R(0)=R_{0}$ and $\psi(\x,0)=\phi_{0}(\x)/R_{0}$ for any $\x\in\overline\Omega$.
It is easy to see that if $R_0 = 1$, then the  coupled system \eqref{prob1} has a unique solution of $R(t) \equiv 1$ and $\psi(\x,t) = \phi(\x,t)$, where $\phi(\x,t)$ is the solution of the $L^2$ gradient flow equation \eqref{prob}. Thus, we always set the initial condition $R_{0}=1$ in what flows. For numerical stabilization, a linear splitting  \cite{DJLQ_rev} is often applied to the gradient flow equation \eqref{111} so that we deal with a transformed equation as
\begin{equation}\label{prob1a}
\dps\frac{\partial (R\psi)}{\partial t}=(\varepsilon^{2}\Delta-s)(R\psi)-(f(R\psi)-sR\psi),\
\end{equation}
where $s\geq0$ is a constant stabilizing parameter.

\paragraph{\bf\em Fourier spectral discretization in space}
 Let us apply the Fourier spectral method for the spatial discretization of the system \eqref{prob1}.
 Note that other space discretization methods can also be used, such as finite difference, finite element or finite volume methods.
The trial function space of the Fourier spectral method is defined as:
\beq
X_{N}(\Omega):={\rm span} \{ e^{i2(kx+ly)\pi/L},\;-N\leq k,l\leq N\},
\eeq
where $N$ is a positive  integer and $i=\sqrt{-1}$. Let $\Pi_{FS}$ denote the projection operator from $L^{2}(\Omega)$ to $X_{N}(\Omega)$ by
$$\Pi_{FS}v:=\dps\sum_{k,l=-N}^{N}(\widehat{\Pi}_{FS}v)_{k,l}e^{i2(kx+ly)\pi/L},\quad \forall\, v\in L^{2}(\Omega)$$
with $$(\widehat{\Pi}_{FS}v)_{k,l}=\frac{(v,e^{i2(kx+ly)\pi/L})}{\|e^{i2(kx+ly)\pi/L}\|^{2}},\quad -N\leq k,l\leq N.$$
Then, we obtain the semi-discretized (in space) system of  \eqref{prob1a} as: given $\dps R(0)=1$ and $\psi_{N}(0)=\Pi_{FS}\phi_{0},$
find $R(t)$ and $\psi_{N}(t)\in X_{N}(\Omega)$ for any  $t\in(0,T]$ such that
\begin{subequations}\label{sem-d1}
 \begin{align}
&\dps\frac{\partial(R\psi_{N})}{d t}=(\varepsilon^{2}\Delta-s)(R\psi_{N})-\Pi_{FS}\mathcal{N}(R\psi_N),\\
&\dps\frac{dE\big[R\psi_{N}\big]}{dt}=-\|\mu\big(R\psi_{N}\big)\|^{2},
  \end{align}
\end{subequations}
where $\mathcal{N}(R\psi_N)=f(R\psi_N)-sR\psi_N.$
Correspondingly,  $\phi_{N}(\x,t):=R(t)\psi_{N}(\x,t)$ is a semi-discrete approximation to the solution of the $L^{2}$ gradient flow equation \eqref{prob}.
 For any $\psi_{N}(t)\in X_{N}(\Omega),$ we can express it as
\bq\label{eq11}
\psi_{N}(\x,t)=\dps\sum_{k,l=-N}^{N}\widehat{\Psi}_{k,l}(t)e^{i2(kx+ly)\pi/L},
\eq
where $\widehat{\Psi}(t) = (\widehat{\Psi}_{k,l}(t))$ is time-dependent coefficient matrix of dimension $(2N+1)\times(2N+1)$.
Substituting the above expression into the semi-discrete system \eqref{sem-d1} and using the Duhamel's principle for  $t\in(t_{n},t_{n+1})$ ($n\geq 0$ and $t_0=0$), we obtain
\bry
\begin{cases}
\dps R(t_{n+1})\widehat{\Psi}_{k,l}(t_{n+1})=e^{\Dt_{n+1}\mathcal{L}_{k,l}}R(t_{n})\widehat{\Psi}_{k,l}(t_{n})-\Big(I^{\widehat{\Pi}_{FS}\mathcal{N}}_{1}\Big)_{k,l},\quad -N\leq k,l\leq N,\\
\dps E\big[R(t_{n+1})\psi_{N}(t_{n+1})\big]-E\big[R(t_{n})\psi_{N}(t_{n})\big]=-I^{\mu}_{2},
\end{cases}
\ery
 where
  \bry
 \mathcal{L}_{k,l}&\;=\dps-\varepsilon^{2}[(2k\pi/L)^{2}+(2l\pi/L)^{2}]-s,\\ \\
\Big(I^{\widehat{\Pi}_{FS}\mathcal{N}}_{1}\Big)_{k,l}&\dps\,=\int_{t_{n}}^{t_{n+1}}e^{(t_{n+1}-\tau)\mathcal{L}_{k,l}}\Big(\widehat{\Pi}_{FS}\mathcal{N}(R(\tau)\psi_{N}(\tau))\Big)_{k,l}d\tau
\ery
for $-N\leq k,l\leq N$, and
$I^{\mu}_{2}\dps\,=\int_{t_{n}}^{t_{n+1}}\|\mu\big(R(\tau)\psi_{N}(\tau)\big)\|^{2}d\tau.$

\paragraph{\bf\em ETD multistep approximation in time}
For any function $u(t)$ defined on $[0,T]$,  denote $P_{r,n}u(t)$ as its Lagrange interpolation polynomial of degree $r$ using the values of $u(t)$ at $t_{n+1},t_{n},\cdots,t_{n+1-r},$
then we have
\beq
 P_{r,n}u(t)=\sum_{j=-1}^{r-1}\omega_{r,j}(t)u(t_{n-j}),\quad t\in [t_{n},t_{n+1}]
\eeq
with $\omega_{r,j}(t)= \prod_{l=-1,l\neq j}^{r-1}\frac{t-t_{n-l}}{t_{n-j}-t_{n-l}}.$ As below, we list some of the polynomial $P_{r,n}u(t)$ up to $r=2$:
\bry
&P_{0,n}u(t)=u(t_{n+1}),\\
&P_{1,n}u(t)=\eta u(t_{n+1})+(1-\eta) u(t_{n}),\\
& P_{2,n}u(t)=\frac{(\gamma_{n+1}\eta+1)\eta}{1+\gamma_{n+1}} u(t_{n+1})+(1-\eta)(1+\gamma_{n+1}\eta) u(t_{n}))+\frac{\gamma^{2}_{n+1}(\eta-1)\eta}{1+\gamma_{n+1}} u(t_{n-1}),
\ery
where $\eta=(t-t_{n})/\Dt_{n+1}$ and $\gamma_{n+1}=\Dt_{n+1}/\Dt_{n}$. Then, the ETD Multistep approximations (in the spirit of Adam-Moulton) of  the integrations $(I^{u}_{1})_{k,l}$ and $I^{u}_{2}$ are given respectively by
\bry
(I^{u}_{1})_{k,l}&\dps\approx\sum_{j=-1}^{r-1}u_{k,l}(t_{n-j})\int_{t_{n}}^{t_{n+1}}e^{(t_{n+1}-\tau)\mathcal{L}_{k,l}}\omega_{r,j}(\tau)d\tau\\
&\dps:=\sum_{j=-1}^{r-1}\alpha_{k,l}^{(r,j)}u_{k,l}(t_{n-j}),\quad -N\leq k,l\leq N,\\
I^{u}_{2}&\dps\approx\sum_{j=-1}^{r-1}u(t_{n-j})\int_{t_{n}}^{t_{n+1}}\omega_{r,j}(\tau)d\tau:=\sum_{j=-1}^{r-1}\beta^{(r,j)}u(t_{n-j}).
\ery
Below  we give some values of $\{\alpha^{(r,j)}_{k,l}\}_{k,l=-N}^N$ and $\{\beta^{(r,j)}_{k,l}\}_{k,l=-N}^N$ for $r=0,1,2$:
\bryl
\dps\alpha^{0,-1}_{k,l}=\alpha_{0};\\
\dps\alpha^{1,-1}_{k,l}=\alpha_{1},~\alpha^{1,0}_{k,l}=\alpha_{0}-\alpha_{1};\\
\alpha^{2,-1}_{k,l}=\frac{\alpha_{1}+\gamma_{n+1}\alpha_{2}}{1+\gamma_{n+1}}, ~\alpha^{2,0}_{k,l}=\alpha_{0}-\alpha_{1}+\gamma_{n+1}(\alpha_{1}-\alpha_{2}), ~\alpha^{2,1}_{k,l}=-\frac{\gamma_{n+1}^{2}}{1+\gamma_{n+1}}(\alpha_{1}-\alpha_{2})\\
\eryl
with
\bryl
\begin{cases}
\alpha_{0}=-\frac{1}{\mathcal{L}_{k,l}}(1-e^{\mathcal{L}_{k,l}\Dt_{n+1}}), ~~\alpha_{1}=-\frac{\Dt_{n+1}-\alpha_{0}}{\mathcal{L}_{k,l}\Dt_{n+1}},~~\alpha_{2}=-\frac{\Dt_{n+1}-2\alpha_{1}}{\mathcal{L}_{k,l}\Dt_{n+1}}, &\mbox{if } \mathcal{L}_{k,l}\neq0,\\
\alpha_{0}=\Dt_{n+1},~~ \alpha_{1}=\frac{\Dt_{n+1}}{2},~~\alpha_{2}=\frac{\Dt_{n+1}}{3}, & \mbox{ otherwise,}\\
\end{cases}
\eryl
and
\bryl
\dps\beta^{0,-1}=\Dt_{n+1};\\
\dps \beta^{1,-1}=\beta^{1,0}=\frac{\Dt_{n+1}}{2};\\
\dps\beta^{2,-1}=\frac{\Dt_{n+1}(2\gamma_{n+1}+3)}{6(1+\gamma_{n+1})},~~\beta^{2,0}=\frac{3+\gamma_{n+1}}{6}\Dt_{n+1},
~~\beta^{2,1}=-\frac{\gamma^{2}_{n+1}\Dt_{n+1}}{6(1+\gamma_{n+1})}.
\eryl

Then, the fully discrete TDSR-ETD scheme with $(r+1)$-th order in time for the continuous  system \eqref{prob1} (and the semi-discrete problem \eqref{sem-d1})  reads: for $n=r,r+1,\cdots, K-1$, find $R^{n+1}$ and $\widehat{\Psi}^{n+1}$ such that
\begin{subequations}\label{full-ds}
 \begin{align}
		&\dps R^{n+1}\widehat{\Psi}^{n+1}_{k,l}=e^{\Dt_{n+1}\mathcal{L}_{k,l}}R^{n}\widehat{\Psi}^{n}_{k,l}-\Big(\mathcal{I}^{\widehat{\Pi}_{FS}\mathcal{N},r}_{1,n+1}\Big)_{k,l},\quad -N\leq k,l\leq N,\label{full-ds1}\\
		&\dps E\big[R^{n+1}\psi^{n+1}_{N}\big]-E\big[R^{n}\psi^{n}_{N}\big]=-\max\{\mathcal{I}^{\mu,r}_{2,n+1},0\},\label{full-ds2}
  \end{align}
\end{subequations}
where
\begin{subequations}\label{int_num}
 \begin{align}
&\Big(\mathcal{I}^{\widehat{\Pi}_{FS}\mathcal{N},r}_{1,n+1}\Big)_{k,l}\dps \;=\sum_{j=-1}^{r-1}\alpha_{k,l}^{(r,j)}\Big(\widehat{\Pi}_{FS}\mathcal{N}\big(R^{n-j}\psi^{n-j}_{N}\big)\Big)_{k,l},\label{int_num1}\\
&\mathcal{I}^{\mu,r}_{2,n+1}\dps\;=\sum_{j=-1}^{r-1}\beta^{(r,j)}\|\mu\big(R^{n-j}\psi_{N}^{n-j}\big)\|^{2}.\label{int_num2}
\end{align}
 \end{subequations}
 Note that $\mathcal{I}^{\mu,r}_{2,n+1}$ approximates $I^{\mu}_{2}$ with $(r+1)$-th order in time and $I^{\mu}_{2}\ge 0$, therefore
$$|\max\{\mathcal{I}^{\mu,r}_{2,n+1},0\}-I^{\mu}_{2}|\leq  |\mathcal{I}^{\mu,r}_{2,n+1}-I^{\mu}_{2}|,$$
and consequently $\max\{\mathcal{I}^{\mu,r}_{2,n+1},0\}$ is also an  $(r+1)$-th order approximation of $I^{\mu}_{2}$. Furthermore, when $r=1$ and $2$,
it is easy to verify that $\beta^{r,j}>0$ for $j=-1,\cdots,r-1$ and thus $\max\{\mathcal{I}^{\mu,r}_{2,n+1},0\}=\mathcal{I}^{\mu,r}_{2,n+1}$.

\begin{remark}
To start  the $(r+1)$-th order TDSR-ETD scheme \eqref{full-ds}, 
the needed values of $\{(R^{k},\widehat{\Psi}^{k})\}_{k=1}^r$  can be computed
using  the TDSR-ETD schemes with $r$-th order or some other methods.
\end{remark}

The proposed TDSR-ETD scheme \eqref{full-ds} is a nonlinear scheme and has to be solved with some iterative algorithms at each time step. When a physically valid solution of  \eqref{full-ds} exists, due the nonnegativity of $\max\{\mathcal{I}^{\mu,r}_{2,n+1},0\}$, it follows directly from \eqref{full-ds2} and $\phi^{n}_{N}:=R^{n}\psi^{n}_{N}$ that
\bq
E\big[\phi^{n+1}_{N}\big]\leq E\big[\phi^{n}_{N}\big],\quad n\geq 0
\eq
i.e., the TDSR-ETD scheme \eqref{full-ds} is  dissipative with respect to the original energy \eqref{EF} in the discrete sense. However,
it usually requires the time step size not be too large to ensure  the existence of the solution and the convergence of the iterative solution process for the scheme \eqref{full-ds}.

\paragraph{\bf \em Extra forcing term for the TDSR factor and Picard iteration}
 We know $R^{n}$ is an approximation to $1$ since $R(t)\equiv1$ and $\psi(\x,t)=\phi(\x,t)$ is the solution of the continuous  problem \eqref{prob1}.
In order to enforce that the numerical approximation for the TDSR factor $R(t)$ remains close to 1 and relax the solvability of the coupled system, we additionally incorporate  the following forcing term for $R(t)$
 \bq\label{ef1}
 \frac{dR^2}{dt}=0, \quad t> 0,
 \eq
 into  \eqref{prob1-b} to obtain a slightly modified system as follows:
\begin{subequations}\label{prob12}
 \begin{align}
&\dps\frac{\partial (R\psi)}{\partial t}-\varepsilon^{2}\Delta(R\psi)+f(R\psi)=0,\qquad (\x,t)\in \Omega\times(0,T],\\
&\dps\frac{d\big(E\big[R\psi\big]+\theta R^2\big)}{dt}=-\|\mu\big(R\psi\big)\|^{2},\qquad\; t\in(0,T],\label{prob12-2}
\end{align}
\end{subequations}
where $\theta\geq0$ is a constant enforcing parameter.
Note that the above system \label{eq1} is again equivalent to the original problem \eqref{prob} and the coupled system \eqref{prob1}.
The motivation is to simultaneously enforce $R(t)$ to stay close to $1$ and (slightly) relax the energy dissipation law so that the time step size restriction for updating  $R(t)$ can be improved.
By following the similar derivations for  \eqref{prob1}, the fully discrete $(r+1)$-th order in time TDSR-ETD scheme for the revised system \eqref{prob12} reads:  for $n=r,r+1,\cdots, K-1$, find $R^{n+1}$ and $\widehat{\Psi}^{n+1}$ such that
\begin{subequations}\label{full-ds-enf}
 \begin{align}
		&\dps R^{n+1}\widehat{\Psi}^{n+1}_{k,l}=e^{\Dt_{n+1}\mathcal{L}_{k,l}}R^{n}\widehat{\Psi}^{n}_{k,l}-\Big(\mathcal{I}^{\widehat{\Pi}_{FS}\mathcal{N},r}_{1,n+1}\Big)_{k,l},\label{full-ds-enf1}\\
		&\dps E\big[R^{n+1}\psi^{n+1}_{N}\big]-E\big[R^{n}\psi^{n}_{N}\big]+\theta((R^{n+1})^2-(R^{n})^2)=-\max\{\mathcal{I}^{\mu,r}_{2,n+1},0\}.\label{full-ds-enf2}
  \end{align}
\end{subequations}
Then, we have the following modified discrete energy dissipation law  for the TDSR-ETD scheme \eqref{full-ds-enf}
\bq\label{equ1}
E\big[\phi^{n+1}_{N}\big]+\theta \big[(R^{n+1})^2-1\big]\leq E\big[\phi^{n}_{N}\big]+\theta \big[(R^{n})^2-1], \quad n\geq 0.
\eq
{We note that the modified discrete energy $E\big[\phi^{n}_{N}\big]+\theta \big[(R^{n})^2-1\big]$ is of $(r+1)$th-order in time approximation to the original discrete energy $E\big[\phi^{n}_{N}\big]$ when $R^{n}$ approximates 1 with $(r+1)$th-order accuracy.} Moreover,  $\{E\big[\phi^{n}_{N}\big]\}_{n\geq 0}$ is uniformly bounded from above  if and only if  the modified discrete energy is also bounded from above.

Finally, from $t_n$ to $t_{n+1}$, the solution $(R^{n+1},\widehat{\Psi}^{n+1})$ of \eqref{full-ds-enf}  can be efficiently solved with the Picard iteration and the decoupling of the system is as follows:

\noindent{\bf Algorithm 1.}
Set $m=0$, $R^{n+1,(0)}=R^{n}$ and $\widehat{\Psi}^{n+1,(0)} = \widehat{\Psi}^{n}$. Then
\begin{itemize}
\item[(i)] compute $R^{n+1,(m+1)}$ by using Newton's method with initial guess $R^{n+1,(m)}$ to solve the scalar nonlinear system (derived from \eqref{full-ds-enf2})
\begin{align*}
&E\big[R^{n+1,(m+1)}\psi^{n+1,(m)}_{N}\big]+\theta (R^{n+1,(m+1)})^2=E\big[R^{n}\psi^{n}_{N}\big]+\theta (R^{n})^2\nonumber\\
&\quad-\max\Big\{\beta^{(r,-1)}\|\mu\big(R^{n+1,(m)}\psi_{N}^{n+1,(m)}\big)\|^{2}+\sum_{j=0}^{r-1}\beta^{(r,j)}\|\mu\big(R^{n-j}\psi_{N}^{n-j}\big)\|^{2},0\Big\};
\end{align*}
\item[(ii)] compute   $\widehat{\Psi}^{n+1,(m+1)}$ by solving the linear system (derived from \eqref{full-ds-enf1})
\begin{align*}
 &R^{n+1,(m+1)}\widehat{\Psi}^{n+1,(m+1)}_{k,l}+\alpha_{k,l}^{(r,-1)}\Big(\widehat{\Pi}_{FS}\mathcal{N}\big(R^{n+1,(m+1)}\psi^{n+1,(m)}_{N}\big)\Big)_{k,l}\nonumber\\
 &\quad=e^{\Dt_{n+1}\mathcal{L}_{k,l}}R^{n}\widehat{\Psi}^{n}_{k,l}
 -\sum_{j=0}^{r-1}\alpha_{k,l}^{(r,j)}\Big(\widehat{\Pi}_{FS}\mathcal{N}\big(R^{n-j}\psi^{n-j}_{N}\big)\Big)_{k,l};
\end{align*}
\item[(iii)] set $m = m+1$. If not convergent, repeat Steps (i) and (ii); otherwise  return $R^{n+1}=R^{n+1,(m+1)}$ and $\widehat{\Psi}^{n+1} = \widehat{\Psi}^{n+1,(m+1)}$.
\end{itemize}
The convergence criterion is set to be $\|\psi^{n+1,(m)}_{N} - \psi^{n+1,(m+1)}_{N}\|_{\infty}\leq \hat\epsilon$ for some pre-given tolerance $\hat\epsilon>0$ in our numerical experiments.

\subsection{The $H^{-1}$ gradient flow with homogeneous Neumann boundary condition}\label{case_II}
Next we turn to present the TDSR-ETD method for the $H^{-1}$ gradient flow problem \eqref{prob} (i.e., $\mathcal{G}_H = \Delta$).
We similarly obtain the following coupled system for
$(R(t),\psi(\x,t))$:
\begin{subequations}\label{prob22}
 \begin{align}
&\dps\frac{\partial (R\psi)}{\partial t}=\Delta(\varepsilon^{2}\Delta(R\psi)-f(R\psi)),\qquad (\x,t)\in \Omega\times(0,T],\\
&\dps\frac{d\big(E\big[R\psi\big]+\theta R^2\big)}{dt}=-\|\nabla\mu\big(R\psi\big)\|^{2},\quad t\in(0,T],\label{prob12-b}
\end{align}
\end{subequations}
with the initial value  $\dps R(0)=1$ and $\psi(\x,0)=\phi_{0}(\x)$ for any $\x\in\overline\Omega$.
Instead of the periodic boundary condition discussed in the previous subsection, we now consider the homogeneous Neumann boundary condition
\bq\label{N_B}
\dps\frac{\partial\phi}{\partial \n}\big|_{\partial\Omega}=\frac{\partial\mu}{\partial \n}\big|_{\partial\Omega}=0,
\eq
where $\n$ denotes the outward unit normal vector on the boundary $\partial\Omega$.

We will use Legendre Galerkin (spectral) discretization in space for the coupled system \eqref{prob22} due to the homogeneous Neumann boundary condition.
Denote by $Q_{N}(\Omega)$ the space of polynomials of degree  less than or equal to $N$ with respect to each variable, and define
$Q^0_{N}(\Omega):=\{ v\in Q_{N}(\Omega): \frac{\partial v}{\partial \n}|_{\partial\Omega}=0\}$.
Let $\Pi_{LP}$ denote the usual interpolation operator from $L^{2}(\Omega)$ into $Q_{N}(\Omega)$.
Then we obtain the semi-discretized (in space) problem of  \eqref{prob22}:
 given $\dps R(0)=1$ and $\psi_{N}(0)=\Pi_{LP}\phi_{0}$,
find $R(t)$, $\psi_{N}(t)\in Q^0_{N}(\Omega)$, and $\mu_{N}(t) \in Q^0_{N}(\Omega)$ for any $q_{N}\in Q^0_{N}(\Omega)$ such that
\begin{subequations}\label{sem-lp}
 \begin{align}
&\dps\Big(\frac{\partial(R\psi_{N})}{\partial t},q_{N}\Big)=-(\nabla\mu_{N},\nabla q_{N}),\label{sem-lp-a} \\
&\dps(\mu_{N},q_{N})=\varepsilon^{2}(\nabla(R\psi_{N}),\nabla q_{N})+s(R\psi_{N},q_{N})+\big(\Pi_{LP}\mathcal{N}(R\psi_{N}),q_{N}\big),\label{sem-lp-b}\\
&\dps\frac{d\big(E\big(R\psi_N\big)+\theta R^2\big)}{dt}=-\|\nabla\mu\big(R\psi_{N}\big)\|^{2}.\label{sem-lp-c}
  \end{align}
\end{subequations}
Let us express $\psi_{N}(\x,t)$ and $\mu_{N}(\x,t)$ as
 \bq \label{eq2}
\psi_{N}(\x,t)\dps=\sum_{k,l=0}^{N-2}\widehat{\Psi}_{k,l}(t)h_{k}(x)h_{l}(y),\qquad \mu_{N}(\x,t)\dps=\sum_{k,l=0}^{N-2}\widehat{\mu}_{k,l}(t)h_{k}(x)h_{l}(y),
 \eq
where $$h_{k}(x)=L_{k}(\textstyle\frac{2x}{L}-1)-\frac{k(k+1)}{(k+2)(k+3)}L_{k+2}(\textstyle\frac{2x}{L}-1),$$
and $L_{k}(\cdot)$ denotes  the Legendre polynomial of degree $k$.
Then we deduce from \eqref{sem-lp-a} and \eqref{sem-lp-b} the following ordinary differential equation system:
\brl\label{sem-lp1}
\begin{cases}
\dps M\frac{d(R\widehat{\Psi})}{dt}M=-S\widehat{\mu}M-M\widehat{\mu}S, \\
\dps M\widehat{\mu}M=\varepsilon^{2}R\big(S\widehat{\Psi}M+M\widehat{\Psi}S\big)+sRM\widehat{\Psi}M
    +\widehat{\mathcal{N}},\\
\end{cases}
\erl
 where $M=(M_{k,l})$, $S=(S_{k,l})$ and $\widehat{\mathcal{N}}$ are  $(N-1)\times(N-1)$ matrices with the elements given by
 \beq
\dps M_{k,l}=(h_{l},h_{k}),\ \ S_{k,l}=(h'_{l},h'_{k}),\ \
\big(\widehat{\mathcal{N}}\big)_{k,l}=(\Pi_{LP}\mathcal{N}(R\psi_{N}),h_{l}(x)h_{k}(y))
 \eeq
for $0\leq k,l\leq N-2.$
 It is readily seen  \cite{STW10} that the mass matrix $M$ is a symmetric positive definite matrix
and the stiffness matrix $S$ is diagonal with positive elements.
Therefore, the ODE system \eqref{sem-lp1} can be rewritten as
\brl\label{sem-lp2}
\begin{cases}
\dps \frac{d(R\widehat{\Psi})}{dt}=-M^{-1}S\widehat{\mu}-\widehat{\mu}(t)(M^{-1}S)^{T}, \\
\dps\widehat{\mu}=\varepsilon^{2}R\big(M^{-1}S\widehat{\Psi}+\widehat{\Psi}(M^{-1}S)^{T}\big)+sR(t)\widehat{\Psi}
    + M^{-1}\widehat{\mathcal{N}} M^{-1}.\\
\end{cases}
\erl
Let $\Lambda$ be the diagonal matrix whose diagonal elements are  $\{\lambda_{k}>0\}_{k=0}^{N-2}$, the eigenvalues of $M^{-1}S$, and $P$ be the corresponding eigenvector matrix which is orthonormal, i.e.,
 \bq\label{EG}
  P^TM^{-1}SP=\Lambda.
  \eq

Define $\overline{\Psi}(t)=P^T\widehat{\Psi}(t)P$ and $\overline{\mu}(t)=P^T\widehat{\mu}(t)P$, then we derive from \eqref{sem-lp2} and  \eqref{EG} that
\bq\label{222111}
\begin{cases}
\dps \frac{d\big(R(t)\overline{\Psi}\big)}{dt}=-\Lambda\overline{\mu}-\overline{\mu}\Lambda, \\
\dps \overline{\mu}=\varepsilon^{2}R\big(\Lambda\overline{\Psi}+\overline{\Psi}\Lambda\big)+sR(t)\overline{\Psi}
    + H,\\
\end{cases}
\eq
where
$H(t):=(MP)^{-1}\widehat{\mathcal{N}}(t) (MP)^{-T}$.
Moreover, \eqref{222111} can be equivalently rewritten as
 \brl\label{eq3}
\begin{cases}
\dps \frac{d\big(R\overline{\Psi}_{k,l}\big)}{dt}=-(\lambda_{k}+\lambda_{l})\overline{\mu}_{k,l},\quad 0\leq k,l\leq N-2. \\
\dps\overline{\mu}_{k,l}=(\varepsilon^{2}\big(\lambda_{k}+\lambda_{l})+s)R\overline{\Psi}_{k,l}
    + (H)_{k,l},\quad 0\leq k,l\leq N-2.\\
\end{cases}
\erl
Therefore, we obtain an equivalent system for \eqref{sem-lp} as follows:
\begin{subequations}\label{eq3-1}
 \begin{align}
&\dps \frac{d\big(R\overline{\Psi}_{k,l}\big)}{dt}=-(\lambda_{k}+\lambda_{l})(\varepsilon^{2}\big(\lambda_{k}+\lambda_{l})+s)R\overline{\Psi}_{k,l}-\mathcal{H}_{k,l}, \quad 0\leq k,l\leq N-2,\\
&\dps\frac{d\big(E\big[R\psi_N\big]+\theta R^2\big)}{dt}=-\|\nabla\mu\big(R(t)\psi_{N}\big)\|^{2},
\end{align}
\end{subequations}
with $\mathcal{H}_{k,l}(t):=(\lambda_{k}+\lambda_{l})(H(t))_{k,l}$ and $\dps \psi_{N}(t)\dps=\sum_{k,l=0}^{N-2}(P\overline{\Psi}(t)P^{T})_{k,l}h_{k}(x)h_{l}(y).$

Similar to the process of deriving the fully discrete scheme \eqref{full-ds-enf}, by applying  Duhamel's principle to \eqref{eq3-1}, we can obtain the
fully discrete TDSR-ETD scheme with $r+1$ order in time for the $H^{-1}$ gradient flow \eqref{prob22} as follows:
for $n = r, r + 1,\cdots, K-1$, find $R^{n+1}$ and $\overline{\Psi}^{n+1}$ such that
\begin{subequations}\label{full-Ne}
 \begin{align}
		&\dps R^{n+1}\overline{\Psi}^{n+1}_{k,l}=e^{\Dt_{n+1}\mathcal{L}_{k,l}}R^{n}\overline{\Psi}^{n}_{k,l}-\Big(\mathcal{I}^{\mathcal{H},r}_{1,n+1}\Big)_{k,l},\label{full-Ne1}\\
		&\dps E\big[R^{n+1}\psi^{n+1}_{N}\big]-E\big[R^{n}\psi^{n}_{N}\big]+\theta((R^{n+1})^2-(R^{n})^2)=-\max\{\mathcal{I}^{\nabla\mu,r}_{2,n+1},0\}, \label{full-Ne3}
  \end{align}
\end{subequations}
with $\psi_{N}^{n+1}=\sum_{k,l=0}^{N-2}(P\overline{\Psi}^{n+1}P^{T})_{k,l}h_{k}(x)h_{l}(y)$, where
\begin{equation*}
\mathcal{L}_{k,l}=-(\lambda_{k}+\lambda_{l})[\varepsilon^{2}\big(\lambda_{k}+\lambda_{l})+s],\quad
\Big(\mathcal{I}^{\mathcal{H},r}_{1,n+1}\Big)_{k,l}\dps=\sum_{j=-1}^{r-1}\alpha_{k,l}^{(r,j)}\big(\mathcal{H}^{n-j}\big)_{k,l},\\
\end{equation*}
for $0\leq k,l\leq N-2$, and
$\mathcal{I}^{\nabla\mu,r}_{2,n+1}\dps=\sum_{j=-1}^{r-1}\beta^{(r,j)}\|\nabla\mu\big(R^{n-j}\psi_{N}^{n-j}\big)\|^{2}.$
From $t_n$ to $t_{n+1}$, the solution $(R^{n+1},\overline{\Psi}^{n+1})$ of \eqref{full-Ne}  can  be efficiently solved using Algorithm 1 in the same way.

The following modified discrete energy dissipation law is again deduced  for the TDSR-ETD scheme \eqref{full-Ne}
\begin{equation}
E\big[\phi^{n+1}_{N}\big]+\theta \big[(R^{n+1})^2-1\big]\leq E\big[\phi^{n}_{N}\big]+\theta\big[ (R^{n})^2-1\big], \quad n\geq 0,
\end{equation}
with $\phi^{n}_{N}=R^{n}\psi^{n}_{N}$.
In addition, the $H^{-1}$ gradient flow \eqref{prob22} is mass-conserved along the time (i.e., $\frac{d}{dt}\int_{\Omega}\phi d\x=0$),
it is easy to check the  TDSR-ETD scheme \eqref{full-Ne} conserves the mass in the discrete sense, that is
\begin{equation}
\int_{\Omega}\phi^{n}_{N}(\x)d\x = \int_{\Omega}\phi^{0}_{N}(\x)d\x, \quad\forall\,n\geq 0.
\end{equation}

\section{Application to other types of gradient flows}\label{sect3}
\setcounter{equation}{0}
 In addition to the above two classic gradient flows (Allen-Cahn and Cahn-Hilliard equations) with respect to \eqref{EF}, the proposed TDSR-ETD method can be naturally applied to many other types of gradient flow problems. In this section, we  will discuss its application to  two challenging phase field models, the molecular beam epitaxial (MBE) model \cite{GL97,li03} and the phase-field crystal (PFC) model \cite{EG04,Eld02}.

\subsection{The MBE model}
In this subsection, we aim to construct efficient energy dissipative numerical schemes  for the MBE model equipped with periodic boundary conditions via the TDSR-ETD method.
The MBE model can be viewed as the $L^{2}$-gradient flow with respect to the free energy functional defined by
\be\label{EF2}
E_M\big[\phi\big]
:=\int_{\Omega}\Big[\frac{\varepsilon^{2}}{2}(\Delta\phi)^{2}+F_M(\nabla\phi)\Big]d\x,
\ee
where the two commonly used nonlinear potential functionals $F_M: \mathbb{R}^{d}\rightarrow \mathbb{R}$ are given by
\bq\label{MBE_p}
F_M(\nabla\phi)=\begin{cases}
\begin{array}{ll}
\frac{1}{4}(|\nabla\phi|^{2}-1)^{2},&\quad \mbox{the case with slope selection,}\\
-\frac{1}{2}\ln(1+|\nabla\phi|^{2}),&\quad \mbox{the case without slope selection.}
\end{array}
\end{cases}
\eq
Then, the MBE model reads
\bq\label{MBE}
\frac{\partial\phi}{\partial t}=-(\varepsilon^{2}\Delta^{2}\phi+f_M(\nabla\phi)), \quad (\x,t)\in \Omega\times(0,T],
\eq
with initial condition $\phi(\x,0)=\phi_{0}(x)$ for any $\x\in\overline{\Omega}$, where $f_M(\nabla\phi)={\delta F_M(\nabla\phi)}/{\delta\phi}$. 
It satisfies the following energy dissipation law:
\beq
\frac{d}{dt}E_M\big[\phi\big]=-\|\varepsilon^{2}\Delta^{2}\phi+f(\nabla\phi)\|^{2}\leq0.
\eeq
Setting $\phi(\x,t)=R(t)\psi(\x,t)$ as before, and then we can rewrite the MBE model \eqref{MBE} into the following
equivalent coupled system for $(R(t),\psi(\x,t))$:
\begin{subequations}\label{prob12-1}
 \begin{align}
&\dps\frac{\partial (R\psi)}{\partial t}=-\varepsilon^{2}\Delta^{2}(R\psi)+s\Delta(R\psi)-\mathcal{N}_{M}(R\nabla\psi), \quad(\x,t)\in \Omega\times(0,T],\\
&\dps\frac{d\big(E_M\big[R\psi\big]+\theta R^2\big)}{dt}=-\|\mu_{M}\big(R\psi\big)\|^{2},\qquad\quad\qquad\quad t\in(0,T],
\dps
\end{align}
\end{subequations}
with $R(0)=1$ and $\psi(\x,0)=\phi_{0}(\x)$ for any $\x\in\overline{\Omega}$, where $\mathcal{N}_{M}(\cdot)$ and $\mu_{M}(\cdot)$ are given respectively by
\begin{equation*}
\mathcal{N}_{M}(\nabla(R\psi))=f_M(\nabla(R\psi))+s\nabla\cdot\nabla(R\psi),\quad
\mu_{M}\big(R\psi\big)=\varepsilon^{2}\Delta^{2}(R\psi)+f_M(\nabla(R\psi)).
 \end{equation*}

Let us apply the Fourier spectral method for the spatial discretization of the system \eqref{prob12-1} as done in subsection \eqref{scase1}
and set $\psi_{N}(0)=\Pi_{FS}\phi_{0}$. For  any $\psi_N\in X_N$ there is a representation \eqref{eq11}.
By following a similar process discussed in the previous section, we obtain the fully-discrete $(r+1)$-th order in time TDSR-ETD scheme for the MBE model \eqref{MBE} as follows:
for $n = r, r + 1, \cdots, K-1$, find $R^{n+1}$ and ${\widehat\Psi}^{n+1}$ such that
\begin{subequations}\label{full-ds-mbe}
 \begin{align}
		&\dps R^{n+1}\widehat{\Psi}^{n+1}_{k,l}=e^{\Dt_{n+1}(\mathcal{L}_{M})_{k,l}}R^{n}\widehat{\Psi}^{n}_{k,l}-\Big(\mathcal{I}^{\widehat\Pi_{FS}\mathcal{N}_{M},r}_{1,n+1}\Big)_{k,l},\label{full-ds-mbe1}\\
		&\dps E_M\big[R^{n+1}\psi^{n+1}_{N}\big]-E_M\big[\phi^{n}_{N}\big]+\theta((R^{n+1})^2-(R^{n})^2)=-\max\{\mathcal{I}^{\mu_{M},r}_{2,n+1},0\},\label{full-ds-mbe2}
  \end{align}
\end{subequations}
where
\bry
 \big(\mathcal{L}_{M}\big)_{k,l}&\;=-\varepsilon^{2}[(2k\pi/L)^{2}+(2l\pi/L)^{2}]^{2}-s[(2k\pi/L)^{2}+(2l\pi/L)^{2}],\\
\Big(\mathcal{I}^{\widehat{\Pi}_{FS}\mathcal{N}_M,r}_{1,n+1}\Big)_{k,l}\dps&\;\dps=\sum_{j=-1}^{r-1}\alpha_{k,l}^{(r,j)}\Big(\widehat{\Pi}_{FS}\mathcal{N}_M\big(\nabla(R^{n-j}\psi^{n-j}_{N})\big)\Big)_{k,l},
\ery
for $-N\leq k,l\leq N$, and
$
\mathcal{I}^{\mu_M,r}_{2,n+1}\dps\;=\sum_{j=-1}^{r-1}\beta^{(r,j)}\|\mu_M\big(R^{n-j}\psi_{N}^{n-j}\big)\|^{2}.
$
Obviously, the TDSR-ETD scheme \eqref{full-ds-mbe} is energy dissipative in the following sense:
 \begin{equation}
E_M\big[\phi^{n+1}_{N}\big]+\theta \big[(R^{n+1})^2-1\big]\leq E_M\big[\phi^{n}_{N}\big]+\theta \big[(R^{n})^2-1\big], \quad n\geq 0,
\end{equation}
where $\phi^{n}_{N} = R^n\psi^n_N$.

\subsection{The PFC model}
The  PFC model was first proposed in \cite{EG04,Eld02} and has been frequently used for modelling crystal growth at the atomic scale in both space and diffusive time scales.
It allows for the nucleation of crystallites at arbitrary locations and orientations while containing elastic and plastic deformations.
In this subsection, we consider the following phase field crystal equation with the periodic boundary condition:
\brl\label{PFC}
\begin{cases}
\dps\frac{\partial \phi}{\partial t}
=\Delta\mu_{p}(\phi), \qquad\qquad  \qquad\quad\;\;\, (\x,t)\in \Omega\times(0,T],\\
\dps\mu_{p}(\phi)=(\Delta+\sigma)^{2}\phi+\phi^{3}-\delta\phi,\quad (\x,t)\in \Omega\times(0,T],
\end{cases}
\erl
with initial condition $\phi(\x,0)=\phi_{0}(\x),$ where $\sigma$ and $\delta$ are two positive constants such that $0<\delta<\sigma^{2}$.
The above PFC model is mass conservative. Moreover, it satisfies the following energy dissipation law:
\beq
\frac{d}{dt} E_p[\phi]=-\|\nabla\mu_{p}(\phi)\|^{2}\leq0,
\eeq
where
\bq
E_p[\phi]=\int_{\Omega}\Big(\frac{1}{2}\phi(\Delta+\sigma)^{2}\phi+\frac{1}{4}\phi^{4}-\frac{\delta}{2}\phi^{2}\Big)d\x.
\eq
Setting $\phi(\x,t)=R(t)\psi(\x,t)$ as before, and then we can rewrite the PFC model \eqref{PFC} into the following
equivalent coupled system for $(R(t),\psi(\x,t))$:
\begin{subequations}\label{PFC1}
 \begin{align}
&\dps\frac{\partial (R\psi)}{\partial t}
=\Delta(\Delta+\sigma)^{2}(R\psi)+s\Delta(R\psi)-\mathcal{N}_{p}(R\psi),\qquad (\x,t)\in \Omega\times(0,T],\\\
&\dps
\frac{d}{dt} \big(E_p[R\psi]+\theta R^2\big)=-\|\nabla\mu_{p}(R\psi)\|^{2},\qquad\qquad\quad\qquad\quad t\in(0,T],
\dps
\end{align}
\end{subequations}
with $R(0)=1, \psi(\x,0)=\phi_{0}(\x)$ for any $\x\in\overline{\Omega}$,
where $$\mathcal{N}_{p}(R\psi)=-\Delta\big((R\psi)^{3}-(s+\delta)R\psi\big).$$

Let us again apply the Fourier spectral method for the spatial discretization of \eqref{PFC1}.
Then we obtain the $(r+1)$-th order in time TDSR-ETD scheme for the PFC model \eqref{PFC} as follows: for $n = r, r + 1, \cdots, K-1$, find $R^{n+1}$ and ${\widehat\Psi}^{n+1}$ such that
\begin{subequations}\label{full-ds-pfc}
 \begin{align}
		&\dps R^{n+1}\widehat{\Psi}^{n+1}_{k,l}=e^{\Dt_{n+1}\big(\mathcal{L}_{p}\big)_{k,l}}R^{n}\widehat{\Psi}^{n}_{k,l}-\Big(\mathcal{I}^{\widehat\Pi_{FS}\mathcal{N}_{p},r}_{1,n+1}\Big)_{k,l},\label{full-ds-pfc1}\\
		&\dps E_p\big[R^{n+1}\psi^{n+1}_{N}\big]-E_p\big[\phi^{n}_{N}\big]+\theta((R^{n+1})^2-(R^{n})^2)=-\max\{\mathcal{I}^{\mu_{p},r}_{2,n+1},0\},\label{full-ds-pfc2}
  \end{align}
\end{subequations}
where
\bry
 \big(\mathcal{L}_{p}\big)_{k,l}&\;=-\big((2k\pi/L)^{2}+(2l\pi/L)^{2}\big)\big[\big(-(2k\pi/L)^{2}-(2l\pi/L)^{2}+\sigma\big)^{2}+s\big],\\
\Big(\mathcal{I}^{\widehat{\Pi}_{FS}\mathcal{N}_p,r}_{1,n+1}\Big)_{k,l}\dps&\;\dps=\sum_{j=-1}^{r-1}\alpha_{k,l}^{(r,j)}\Big(\widehat{\Pi}_{FS}\mathcal{N}_p\big(R^{n-j}\psi^{n-j}_{N}\big)\Big)_{k,l},
\ery
for $-N\leq k,l\leq N$, and
$
\mathcal{I}^{\mu_M,r}_{2,n+1}\dps:=\sum_{j=-1}^{r-1}\beta^{(r,j)}\|\mu_p\big(R^{n-j}\psi_{N}^{n-j}\big)\|^{2}.
$
It is easy to check that the TDSR-ETD scheme \eqref{full-ds-pfc} is mass conservative and energy dissipative in the following sense:
 \beq
E_p\big[\phi^{n+1}_{N}\big]+\theta \big[(R^{n+1})^2-1\big]\leq E_p\big[\phi^{n}_{N}\big]+\theta \big[(R^{n})^2-1\big],\quad n\geq 0,
\eeq
where $\phi^{n}_{N} = R^n\psi^n_N$.

\section{Numerical experiments}
\setcounter{equation}{0}

In this section,  various numerical experiments and comparisons are presented to validate the theoretical results derived for the proposed TDSR-ETD schemes. In the following experiments, the tolerance error for the Picard iteration in Algorithm 1 at each time step is always set to be $\hat\epsilon = 10^{-7}$ unless specified otherwise.
 All the errors are computed under the $L^{\infty}$-norm measure throughout the numerical tests. Specifically, we test the TDSR-ETD schemes \eqref{full-ds-enf}  with $r=1$ and $r=2$, which are second-order and third-order accurate in time and thus referred to as TDSR-ETD2 and  TDSR-ETD3 respectively.

\subsection{Convergence tests}

\begin{example}\label{expl1}
Consider the gradient flow equation \eqref{prob} with the double-well potential $F(\phi) = \frac14(\phi^2-1)^2$:
\bq\label{ex2}
\begin{cases}
\begin{array}{ll}
\dps\frac{\partial \phi}{\partial t}+\mathcal{G}_{H}\big(\varepsilon^{2}\Delta \phi-\phi(\phi^{2}-1)\big)=0,
\quad &\x\in \Omega, t>0,\\
\phi(\x,0)=\phi_{0}(\x), \qquad\qquad\qquad\qquad\;\; &\x\in \Omega,
\end{array}
\end{cases}
\eq
where $\mathcal{G}_{H}=-I$ (the $L^2$ gradient flow, i.e., Allen-Cahn equation) or $\mathcal{G}_{H}=\Delta$ (the $H^{-1}$ gradient flow,  i.e., Cahn-Hilliard equation) as defined in \eqref{ef}, and
\beq
\phi_{0}(x,y)=
\begin{cases}
\begin{array}{ll}
\sin 2x\cos 3y,&\quad \mbox{if $\mathcal{G}_{H}=-I$,}\\
0.1(\cos 3x\cos 2y+\cos 5x\cos 5y),&\quad \mbox{if $\mathcal{G}_{H}=\Delta$.}
\end{array}
\end{cases}
 \eeq
 The periodic  and the homogeneous Neumann boundary conditions are imposed for the above $L^2$ and $H^{-1}$ gradient flows, respectively. We set the interfacial width parameter $\varepsilon^2=0.01$ with the computational domain $\Omega=(0,2\pi)^2$ for the $L^2$ gradient flow, and $\varepsilon^2=2.5\times 10^{-3}$ with $\Omega=(-1,1)^2$ for the $H^{-1}$ gradient flow.
\end{example}

\begin{figure*}[!t]
\begin{minipage}[t]{0.49\linewidth}
\centerline{\includegraphics[scale=0.3]{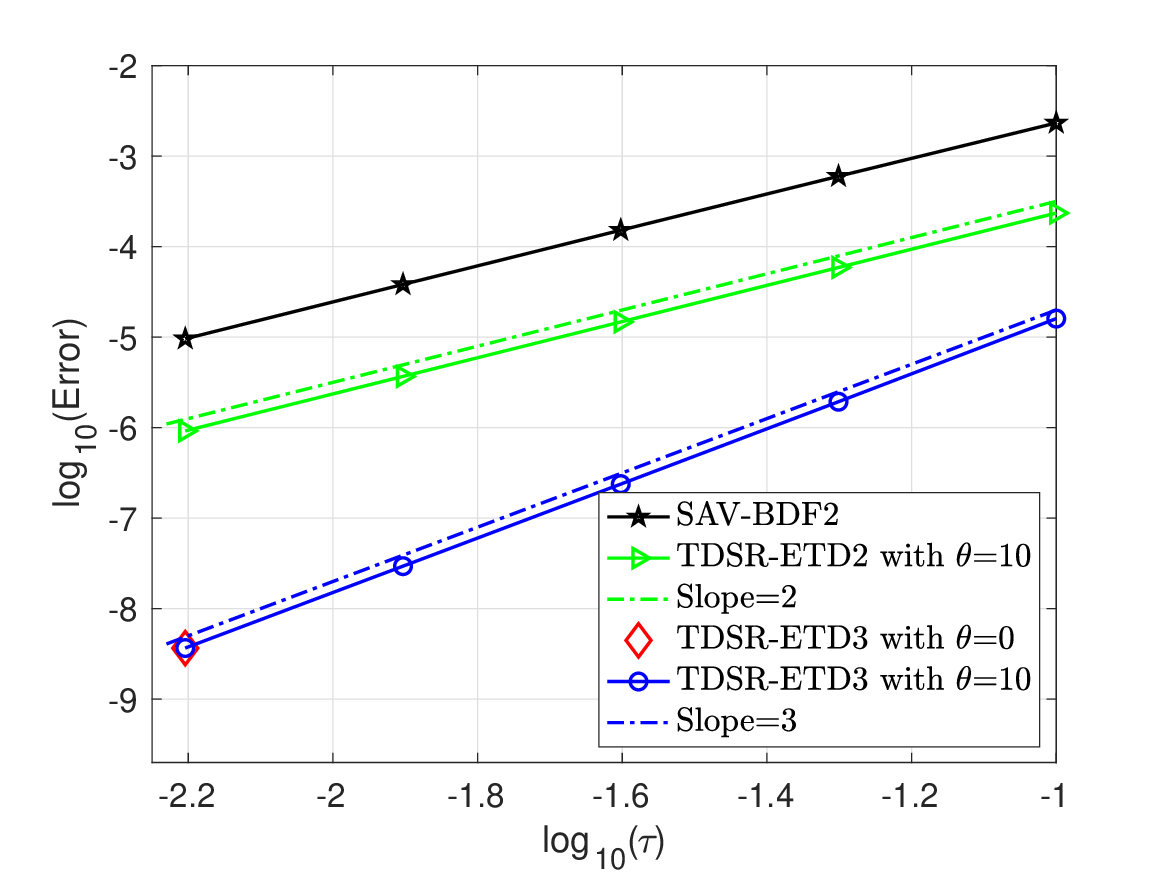}}
\centerline{\small (a) $\phi$}
\end{minipage}
\begin{minipage}[t]{0.49\linewidth}
\centerline{\includegraphics[scale=0.3]{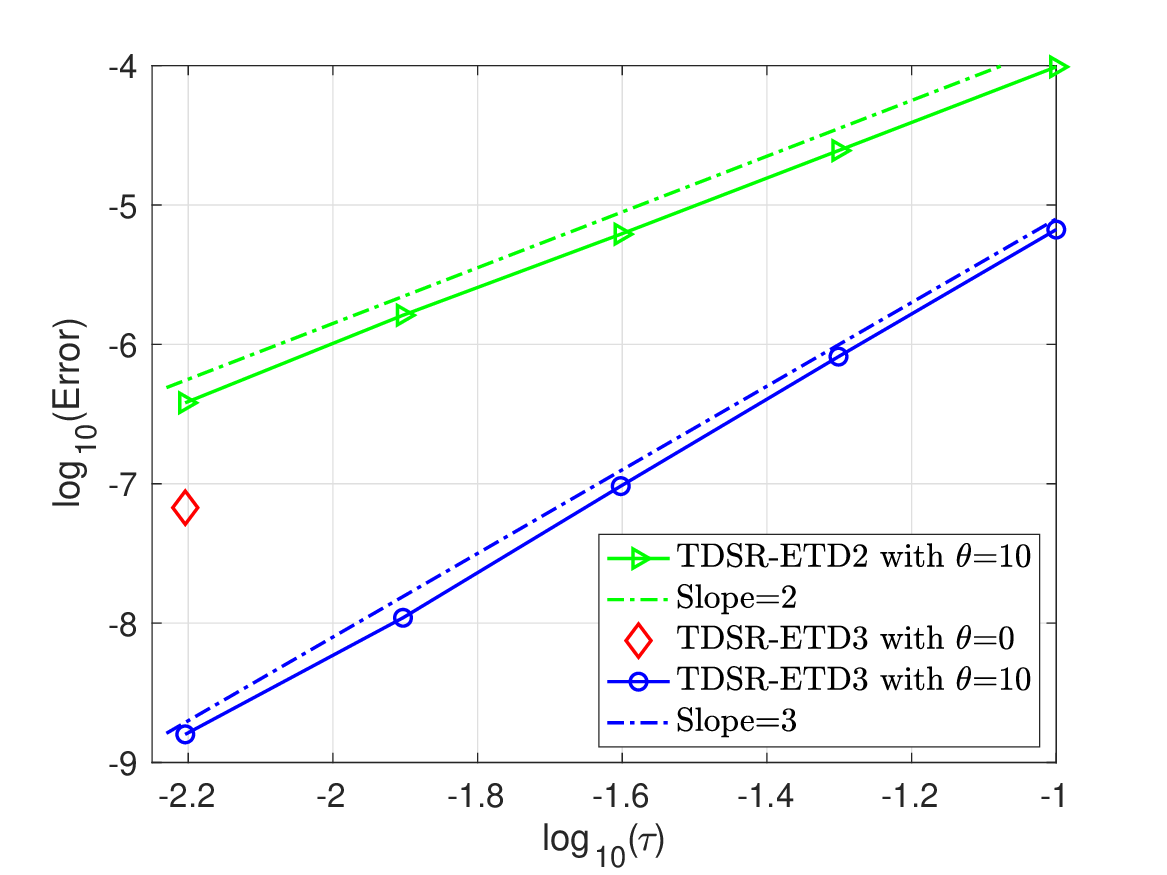}}
\centerline{\small (b)  $R$ }
\end{minipage}
\caption{(Example \ref{expl1}) Plots of the numerical errors of $\phi$  and $R$  with respect  to  the time step size $\tau$ produced by the TDSR-ETD2 and TDSR-ETD3 schemes \eqref{full-ds-enf}  for the $L^2$ gradient flow. For the purpose of comparison, the corresponding result of $\phi$ produced by the SAV-BDF2 scheme is added in the left panel.
}\label{fig1}
\end{figure*}

 We employ the proposed TDSR-ETD schemes \eqref{full-ds-enf} and \eqref{full-Ne} respectively, to compute the numerical solutions of \eqref{ex2} for these two types of gradient flows. The Fourier spectral method with $128\times128$ basis modes is applied in the spatial discretization of the $L^2$ gradient flow, and the spatial discretization used  for the $H^{-1}$ gradient flow is the Legendre spectral method with $256\times256$ basis modes. The uniform time step is used and we set the stabilizing parameter $s=2$ as done in \cite{DJLQ_rev} for this test.
It holds $R(t)\equiv1$ for the continuous problem but there are no exact solutions of $\phi$ or $\psi$ available for the tested two gradient flows. Thus we take the approximate solutions produced by the TDSR-ETD3 scheme with the time step size $\Dt=10^{-5}$ as the reference solution for computing  numerical solution errors on $\phi$ and $R$ at $T=1$.

 For the $L^{2}$ gradient flow ($\mathcal{G}_{H}=-I$, Allen-Cahn equation), we choose a set of uniformly refined time step sizes $\tau=0.1\times 2^{-k}, k=0,1,\cdots,4$.  We set $\hat\epsilon = 10^{-10}$ for the Picard iteration in this test to catch $R$ very accurately to produce the expected order of convergence for $R$, especially for the high-order scheme TDSR-ETD3. Figure \ref{fig1} presents the numerical errors of $\phi$  and $R$ computed by the TDSR-ETD2 and  TDSR-ETD3 schemes as functions of the time step size in log-log scale, in which two choices of  the forcing parameter are tested, $\theta=0$ and $\theta=10$ respectively. Moreover, for comparison, we also display in Figure \ref{fig1}-(a) the error behavior of $\phi$  by the SAV-BDF2 scheme.
 When  $\theta=0$, we find that  TDSR-ETD2 always fails due to non-convergence of Algorithm 1 and TDSR-ETD3 works only for the time step size $\tau= 1/160$. When $\theta=10$, both TDSR-ETD2 and TDSR-ETD3 work well for all tested time step sizes.
  This implies that the introduction of the enforcing parameter $\theta$ does greatly improve the time step size restriction for the TDSR-ETD schemes.
It is observed as expected  that the desired second-order accuracy in time of TDSR-ETD2  and third-order accuracy  of  TDSR-ETD3 are achieved for  numerical solutions of both $\phi$ and $R$, and the numerical errors of TDSR-ETD3 are among 10 to 100 times smaller than those of TDSR-ETD2, which demonstrate the advantage of higher-order accuracy. In addition, we see  from Figure \ref{fig1} that the numerical errors of TDSR-ETD2 always remain about 10 times smaller than those of SAV-BDF2 under the same time step size, which indicates TDSR-ETD2 performs much better than SAV-BDF2 even though they are both second-order accurate in time.

\begin{figure*}[!t]
\centerline{ \includegraphics[scale=0.3]{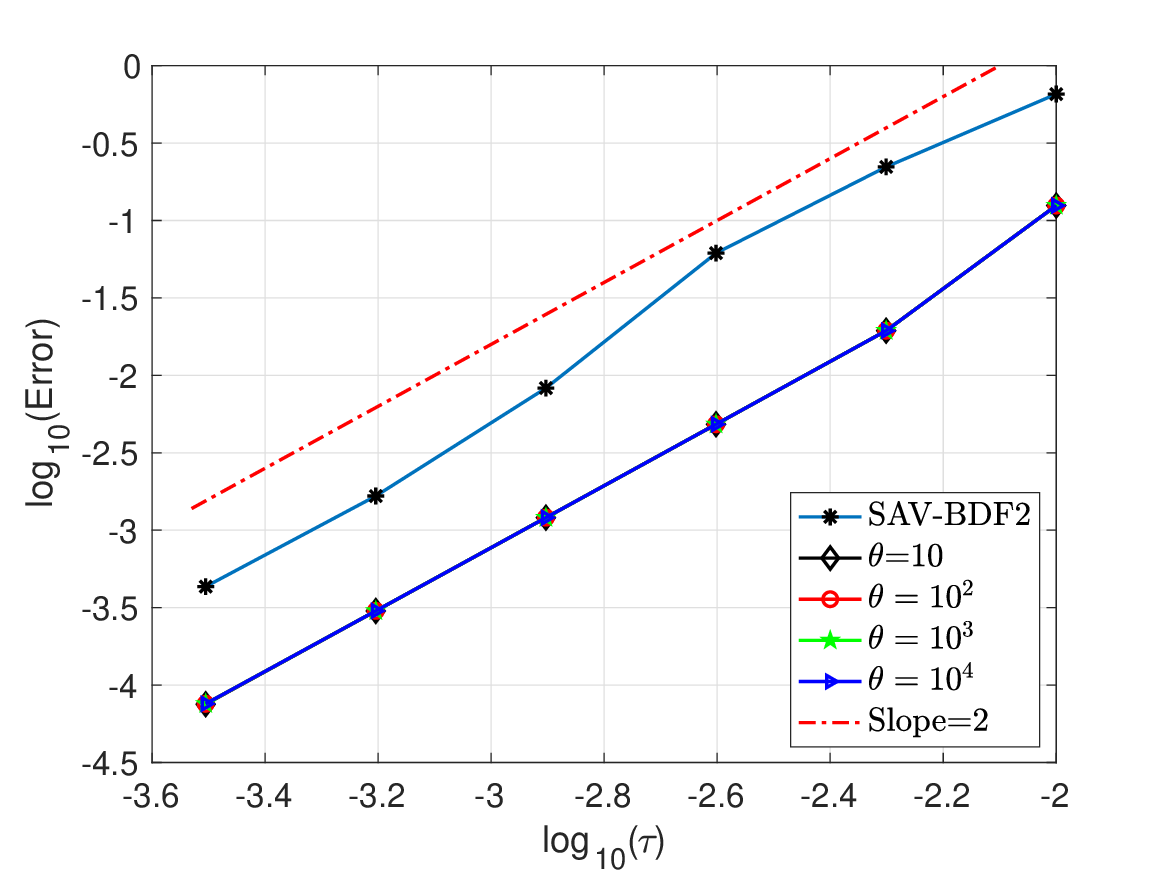}   \includegraphics[scale=0.3]{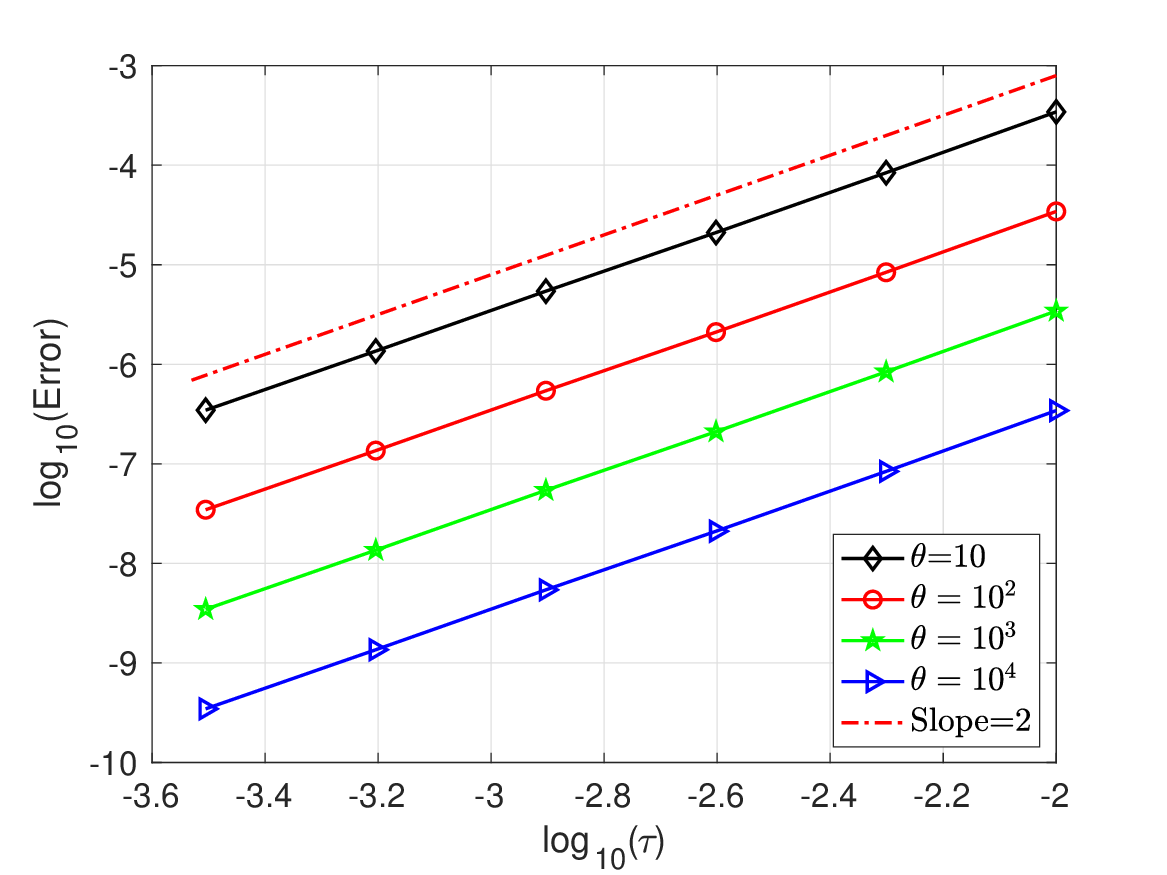} }
\centerline{\small (a) TDSR-ETD2 with $\theta=10, 10^{2},10^{3}$ and $10^{4}$ and SAV-BDF2. Left: $\phi$, right: R.}
\centerline{ \includegraphics[scale=0.3]{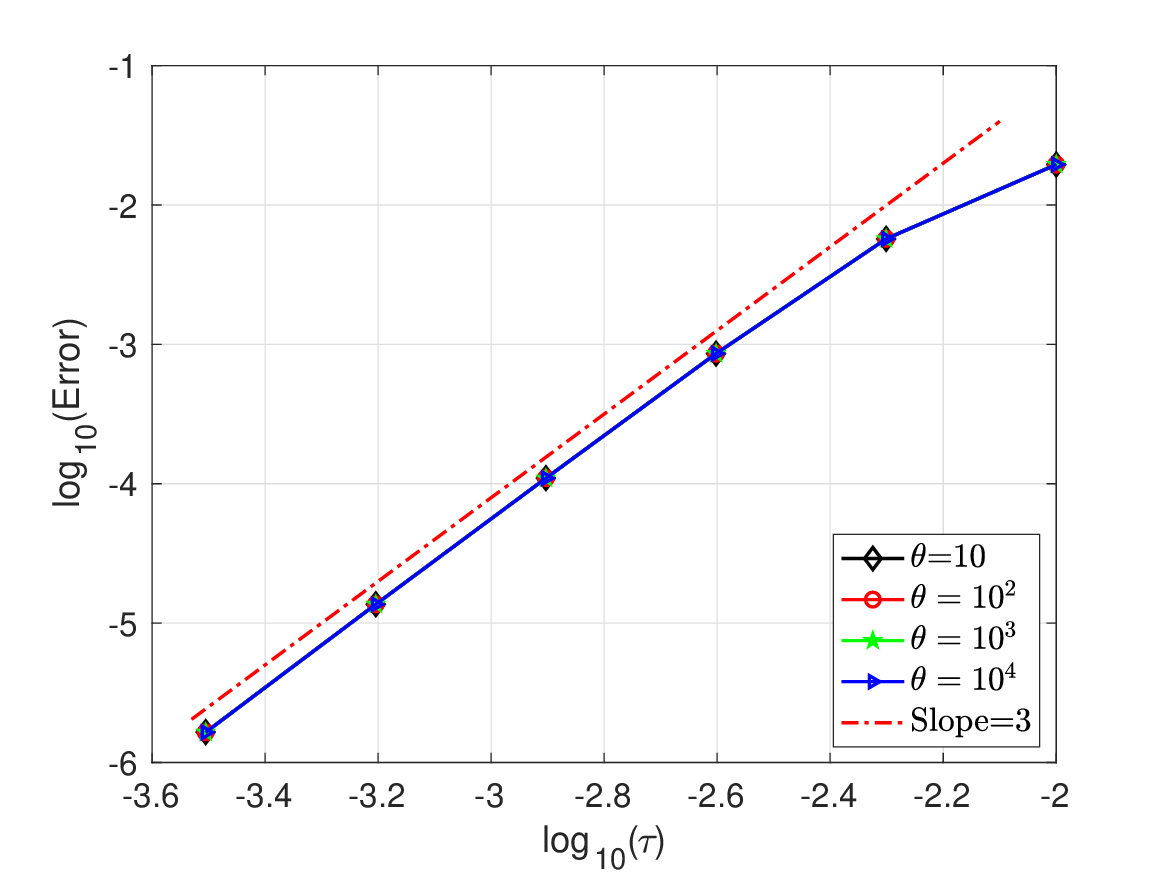}  \includegraphics[scale=0.3]{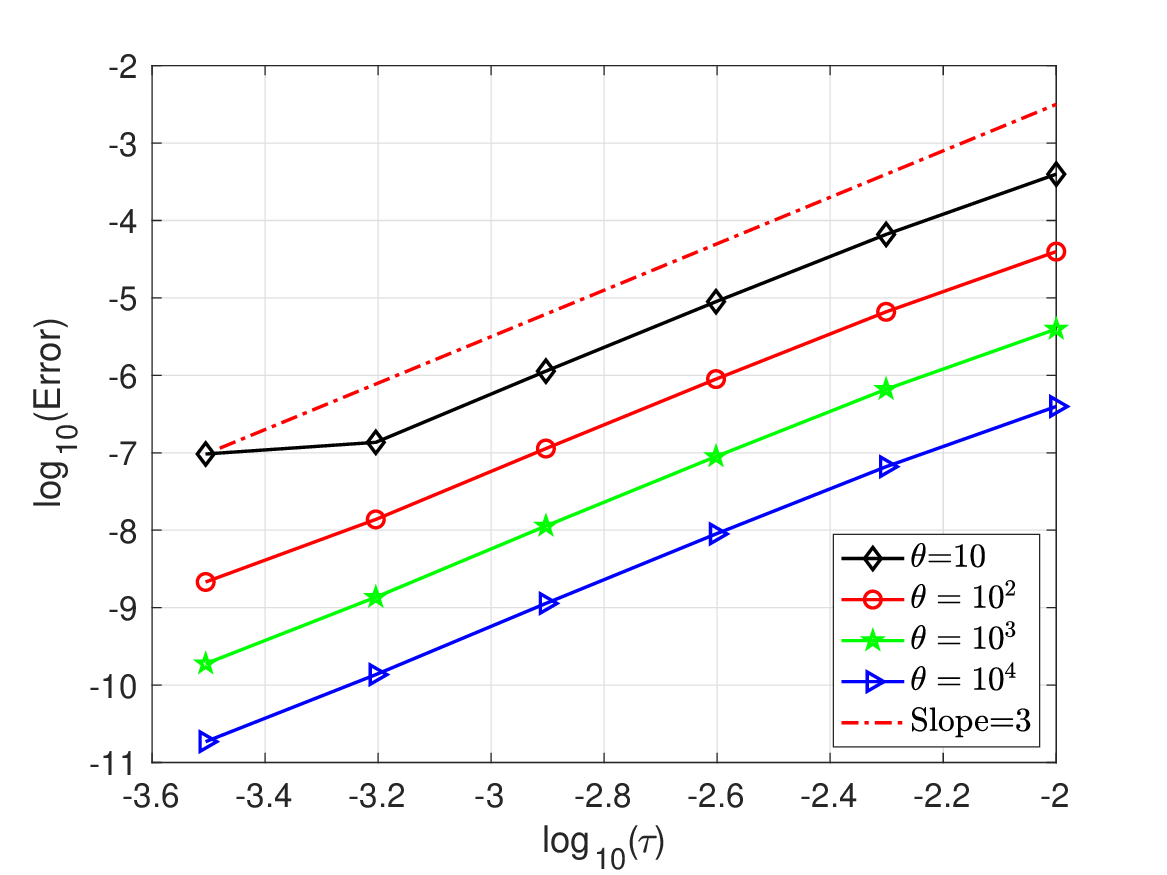} }
\centerline{\small (b) TDSR-ETD3 with $\theta=10, 10^{2},10^{3}$ and $10^{4}$. Left: $\phi$, right: R.}
\caption{(Example \ref{expl1}) Plots of the numerical errors of  $\phi$  and $R$ with respect to the time step size $\tau$ produced by the TDSR-ETD2 and TDSR-ETD3 schemes \eqref{full-Ne}   for the $H^{-1}$ gradient flow. For the purpose of comparison, the corresponding result of $\phi$ produced by the SAV-BDF2 scheme is added in the left panel of (a).
}\label{fig2}
\end{figure*}

\begin{figure*}[!ht]
\begin{minipage}[t]{0.49\linewidth}
\centerline{\includegraphics[scale=0.3]{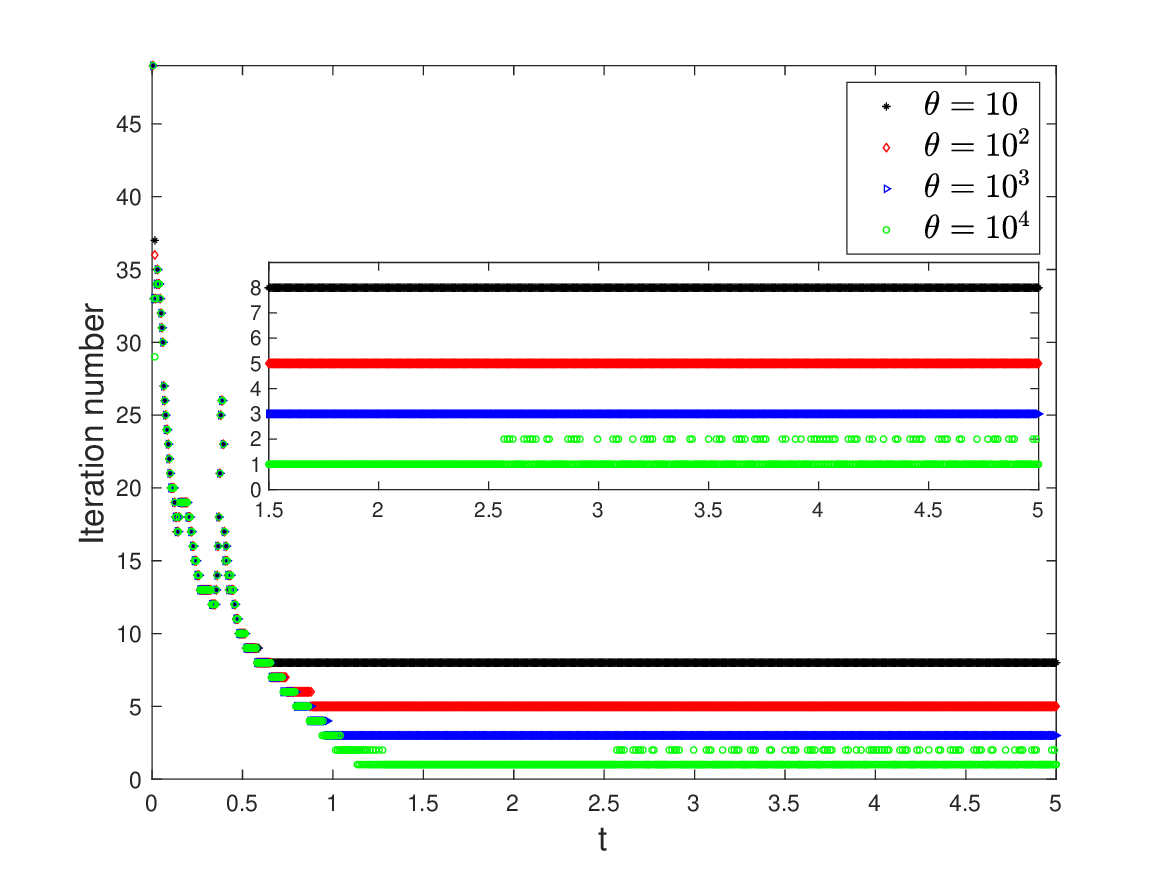}}
\centerline{\small (a) $\tau=5\times 10^{-3}$}
\end{minipage}
\begin{minipage}[t]{0.49\linewidth}
\centerline{\includegraphics[scale=0.3]{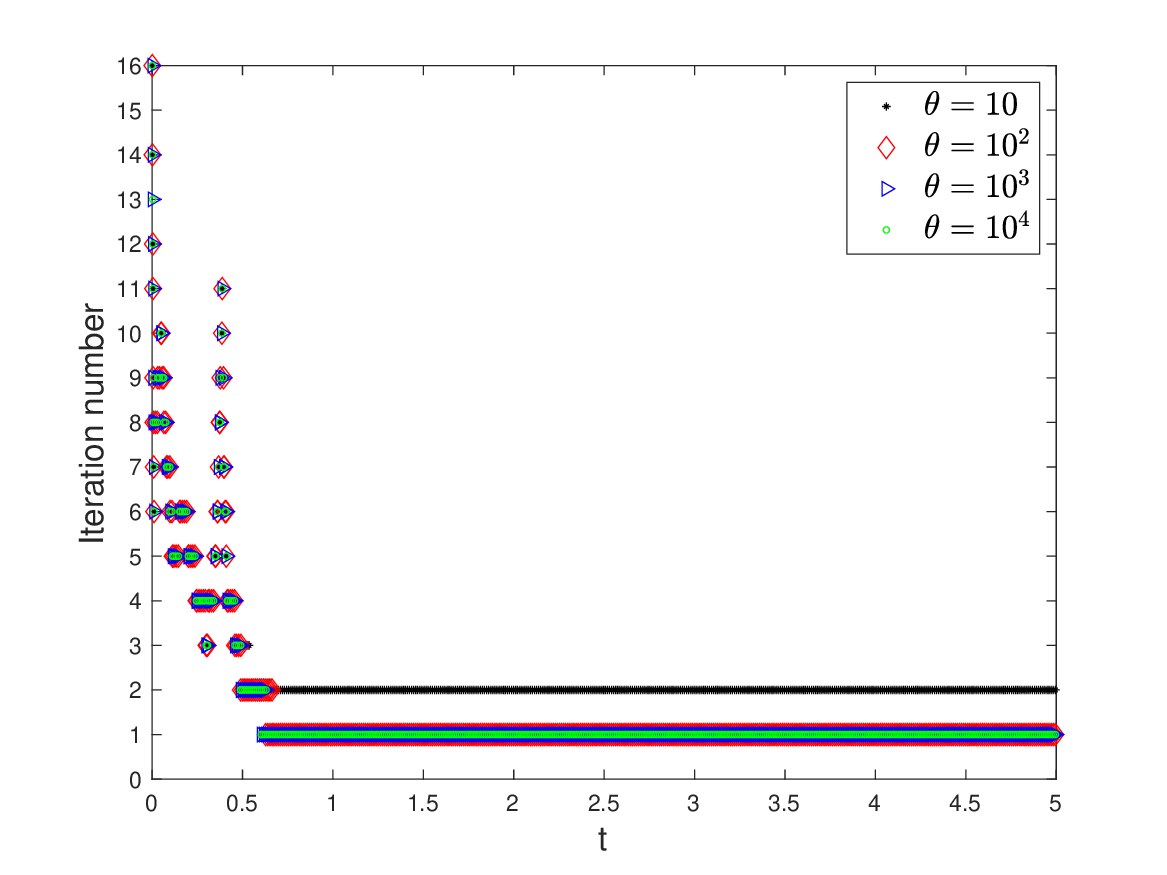}}
\centerline{\small (b) $\tau=10^{-3}$ }
\end{minipage}
\caption{(Example \ref{expl1}) Evolution  of the number of Picard iterations of Algorithm 1 at each time step (up to the time $T=5$) for the TDSR-ETD3 scheme \eqref{full-Ne}  with uniform time stepping  for the $H^{-1}$ gradient flow.
}\label{f1}
\end{figure*}

Next, we investigate the effect of the enforcing parameter $\theta$ on the convergence behaviors of $\phi$ and $R$, and on the needed numbers of  Picard iterations for convergence  in Algorithm 1 for the TDSR-ETD2 and TDSR-ETD3 schemes by testing the $H^{-1}$ gradient flow ($\mathcal{G}_{H}=\Delta$, Cahn-Hilliard equation). We choose several different enforcing parameters as $\theta=10, 10^{2}, 10^{3},$ and $10^{4}$ and a set of uniformly refined time step sizes $\tau=0.01\times 2^{-k}, k=0,1,\cdots,6$. Figure \ref{fig2} presents the error behaviors of the two TDSR-ETD schemes with respect to the time step sizes. Moreover, we also plot in Figure \ref{fig2}-(a)  the numerical errors  by the SAV-BDF2 scheme to make a comparison. It is shown as expected that TDSR-ETD2 and TDSR-ETD3 are of second-order and third-order accuracy in time for both $\phi$ and $R$, respectively, for all tested cases. Moreover, TDSR-ETD2 again achieves about 10 times better accuracy than SAV-BDF2 under the same time step size. In addition, we observe  from Figure \ref{fig2} that the accuracy of the numerical solutions for $\phi$ isn't  affected by the choice of these four values for $\theta$ but larger $\theta$ forces the numerical solution of $R$ to be closer to the exact solution $1$.
The evolution in time of the number of Picard iterations of Algorithm 1 at each time step (up to a longer  time $T=5$ for better illustration) is plotted in Figure \ref{f1} for the TDSR-ETD3 scheme with two time step sizes $\tau=5\times 10^{-3}$ and $\tau=10^{-3}$. It shows that a larger forcing parameter $\theta$ can  reduce the number of Picard iterations and is thus more efficient, particularly for the case of larger time step size such as $\tau=5\times 10^{-3}$ used in the simulation.
Therefore, we  will by default set $\theta=10^{4}$ in all remaining numerical experiments.

\begin{example}\label{expl2}
Consider the following MBE model with the periodic boundary condition:
\bq\label{ex3}
\begin{cases}
\begin{array}{ll}
\dps\frac{\partial \phi}{\partial t}+\varepsilon^{2}\Delta^{2} \phi+f_M(\nabla\phi)=0, \quad\qquad\qquad\quad\quad &\x\in (0,2\pi)^2, t>0,\\
\phi(\x,0)=0.1(\sin3x\sin2y+\sin5x\sin5y),\quad\;\,&\x\in (0,2\pi)^{2},
\end{array}
\end{cases}
\eq
where $\varepsilon^{2}=0.01$ and $f_M(\nabla\phi)=\delta F_M(\nabla\phi)/\delta\phi$ with  $F_M$ defined in \eqref{MBE_p}.
\end{example}

\begin{figure*}[!t]
\begin{minipage}[t]{0.49\linewidth}
\centerline{\includegraphics[scale=0.3]{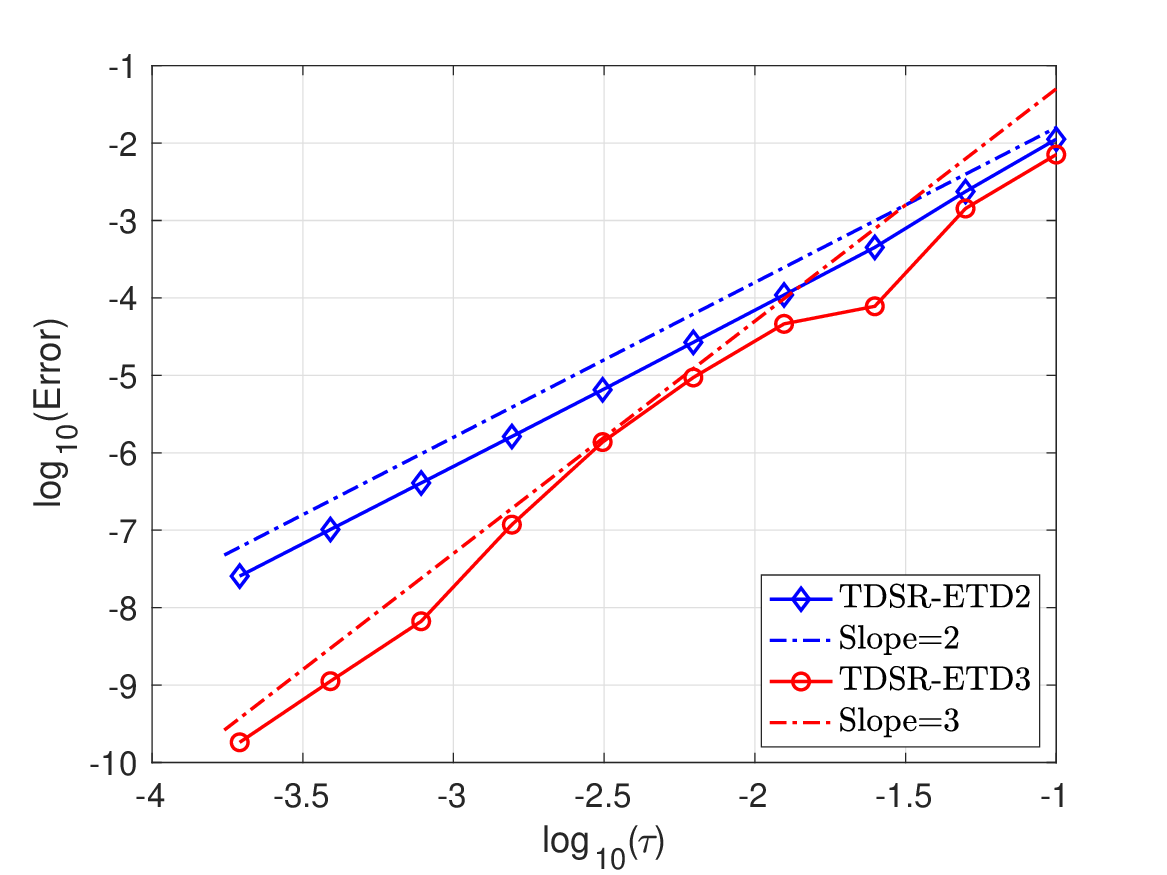}}
\centerline{\small (a) With slope selection}
\end{minipage}
\begin{minipage}[t]{0.49\linewidth}
\centerline{\includegraphics[scale=0.3]{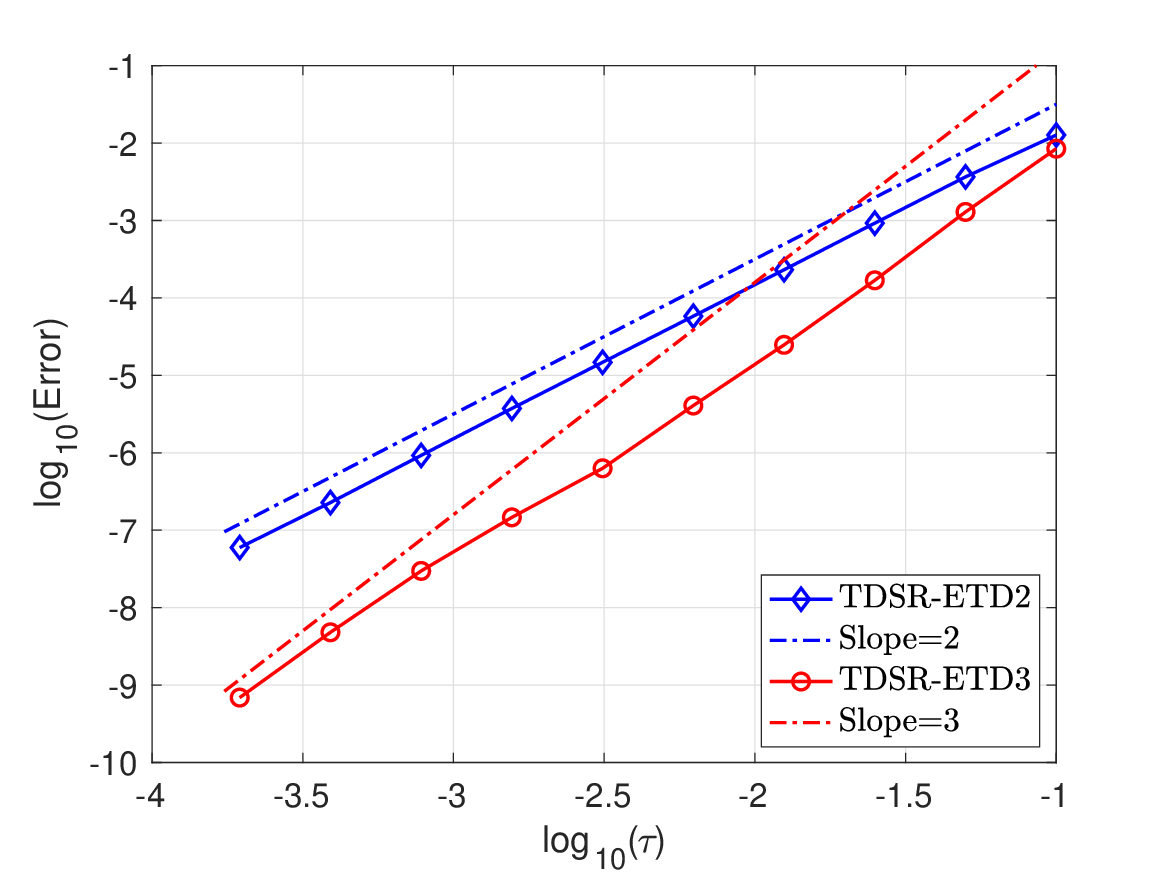}}
\centerline{\small (b)  Without slope selection }
\end{minipage}
\caption{(Example \ref{expl2}) Plots of the numerical errors of $\phi$  with respect  to  the time step size $\tau$ produced by the TDSR-ETD2 and TDSR-ETD3 schemes  \eqref{full-ds-mbe} for the MBE model.
}\label{fig3}
\end{figure*}

We adopt the TDSR-ETD scheme \eqref{full-ds-mbe} with $128\times128$ Fourier modes and uniform time steps, and set the stabilizing parameter $s=2$ and $s=1/8$ for the model with slope selection and without slope selection \cite{JLQZ18},  respectively.
Moreover, the numerical solution computed by the TDSR-ETD3 scheme with the  small  time step $\tau=10^{-5}$ is used as a reference solution at $T=1$ for calculating the numerical solution errors.
Figure \ref{fig3} plots the numerical errors of $\phi$ with respect to the time step size $\tau$ in the log-log scale for the MBE model with and without slope selection, respectively. The results are again in good agreement with the expected convergence rate in time, i.e., second order  for  TDSR-ETD2 and third order for TDSR-ETD3.

\begin{figure*}[!t]
\centerline{\includegraphics[scale=0.3]{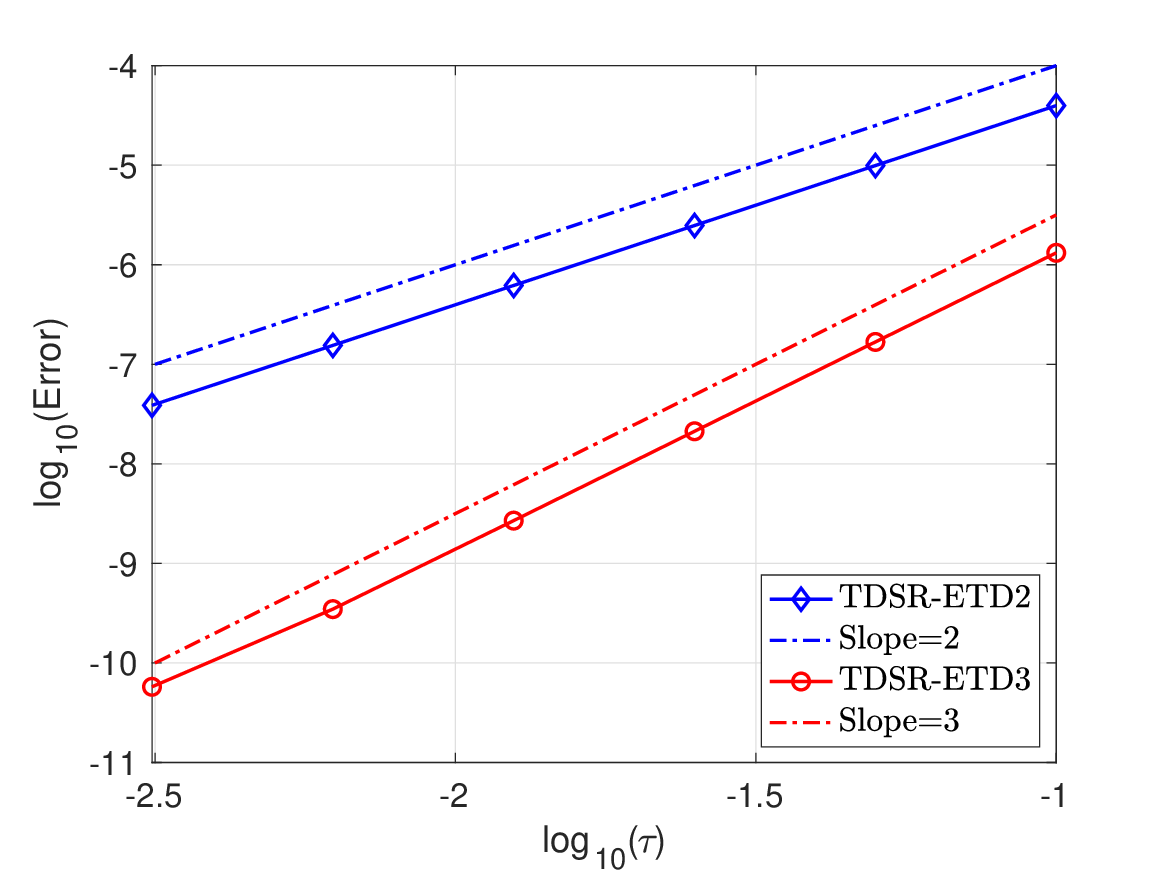}}
\vspace{-0.2cm}
\caption{(Example \ref{expl3}) Plots of the numerical errors of  $\phi$  and $R$ with respect to the time step size $\tau$ produced by the TDSR-ETD2 and TDSR-ETD3 schemes \eqref{full-ds-pfc} for the PFC model.
}\label{fig4}
\end{figure*}

\begin{example}\label{expl3}
Consider following PFC model \eqref{PFC} with the periodic boundary condition and $\sigma=1$ and $\delta=0.025$:
\bq\label{ex33}
\begin{cases}
\begin{array}{ll}
\dps\frac{\partial \phi}{\partial t}
=\Delta\big((\Delta+\sigma)^{2}\phi+\phi^{3}-\delta\phi\big),\quad\;\; &\x\in (0,32)^2, t>0,\\
\phi(\x,0)=\sin\Big(\frac{\pi x}{16}\Big)\cos\Big(\frac{\pi y}{16}\Big),\quad\qquad&\x\in (0,32)^{2}.
\end{array}
\end{cases}
\eq

\end{example}
We adopt  the TDSR-ETD scheme \eqref{full-ds-pfc} with  $256\times256$ Fourier modes and uniform time stepping. We set the stabilizing parameter  $s=\delta$.
To validate the temporal accuracy, a reference solution of the PFC model \eqref{PFC} at $T=1$ is computed by the TDSR-ETD3 scheme with a small enough time step size $\Dt=10^{-5}$.
In Figure \ref{fig4}, the numerical error as a function of the time step size is plotted in the log-log scale.
It is observed that TDSR-ETD2 and TDSR-ETD3  again achieve the expected second order and third order accuracy in time, respectively.

\subsection{Combined with a time-adaptive strategy}
One of the main advantages of the energy stable or energy dissipative numerical schemes is that they can be easily and naturally applied in a time-adaptive strategy.
In this subsection, we investigate the efficiency of the proposed TDSR-ETD schemes combined with a time-adaptive strategy.
There already exist several efficient time-adaptive strategies to the energy stable schemes,
including the ones \cite{GH11,STY16} according to local relative errors, and those \cite{HJQ22,HQ21,QZT11} based on the changing rates of free energy.
It is known that the coarsening dynamics usually experience large changing rates of free energy in some time intervals and small energy changing rates during some others.
Therefore, we here adopt the following robust time-adaptive strategy based on the energy variation \cite{QZT11}:
\bq\label{adp}
\dps\Dt_{n+1}=\max\big(\Dt_{min},\frac{\Dt_{max}}{\sqrt{1+\gamma|E^{'}(t_n)|^{2}}}\big),
\eq
where $\Dt_{min},\Dt_{max}$ are predetermined minimum and maximum time step sizes,
$\gamma$ is a positive constant to be determined.
Obviously, the adaptive time-stepping scheme will automatically select small time step sizes when the energy changing rate is large, and large time step sizes when small.  We will make use of the TDSR-ETD3 schemes with $\theta=10^{4}$ in combination with the above time-adaptive strategy \eqref{adp}.

\begin{figure*}[!t]
\centerline{ \includegraphics[scale=0.16]{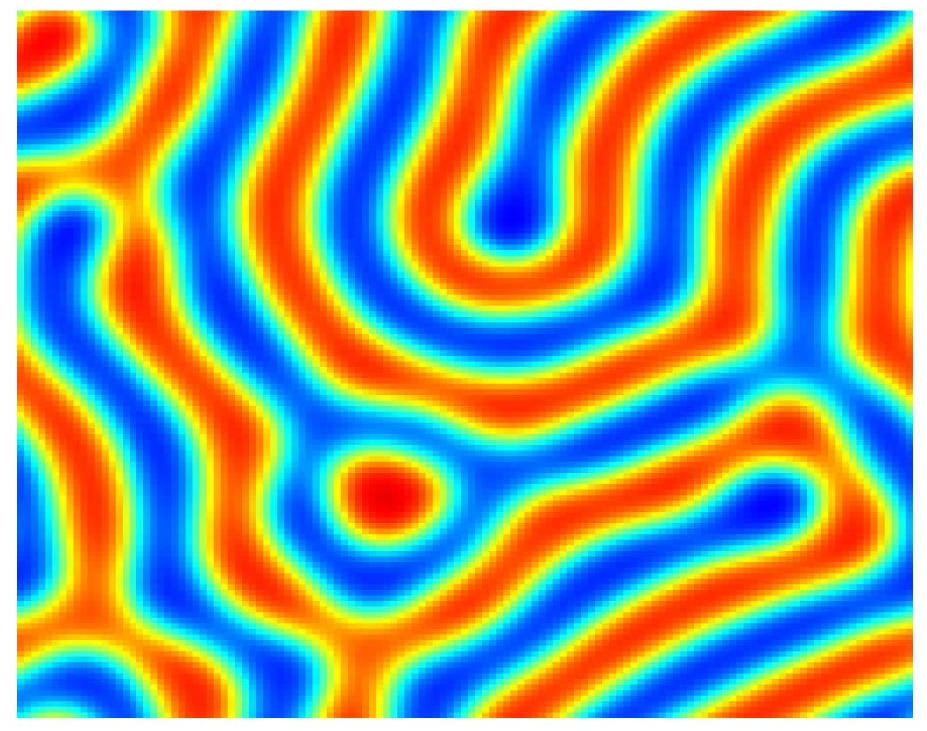}\includegraphics[scale=0.16]{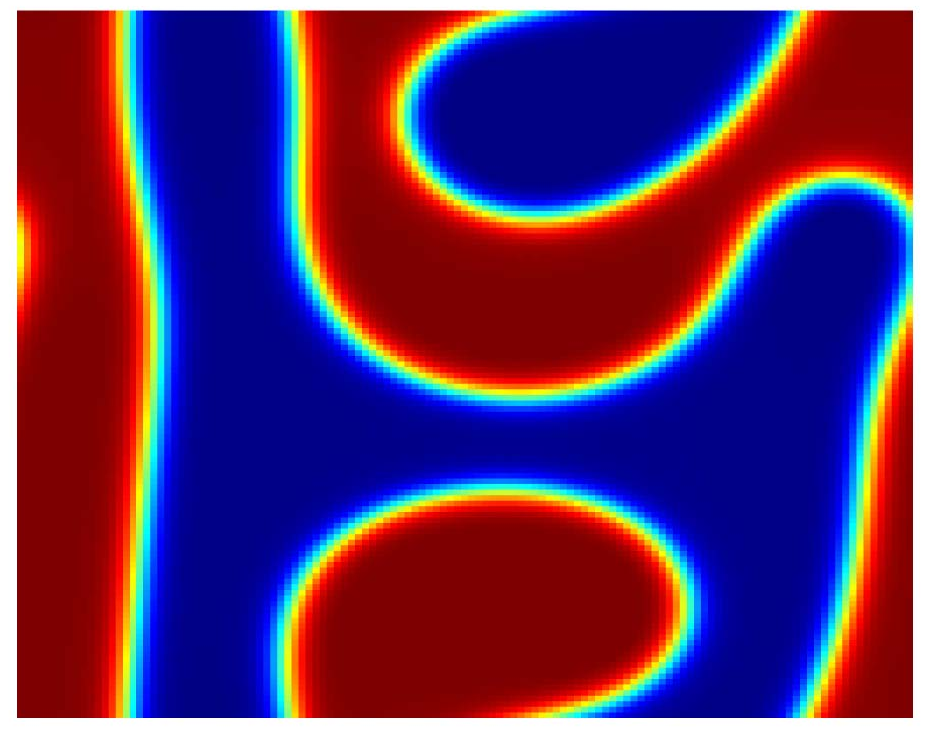}\includegraphics[scale=0.16]{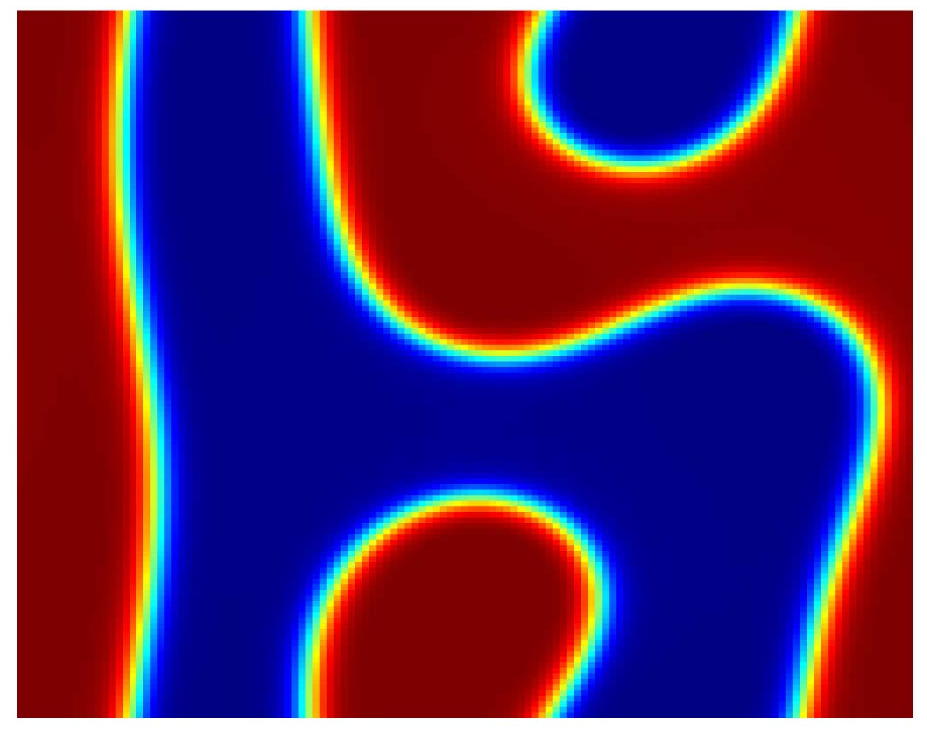}\includegraphics[scale=0.16]{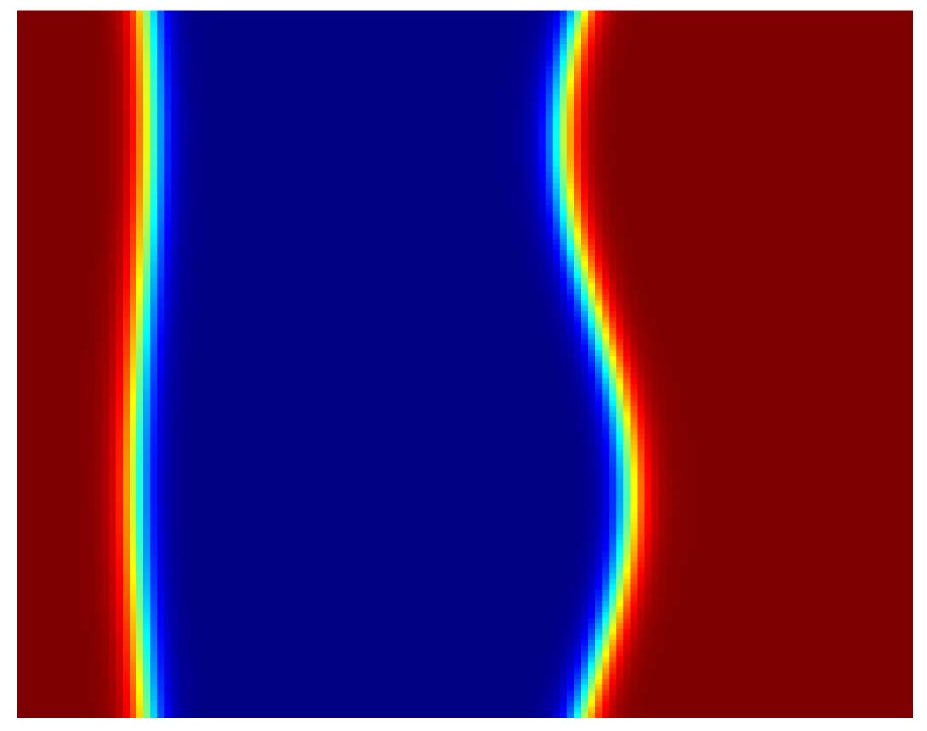}\includegraphics[scale=0.16]{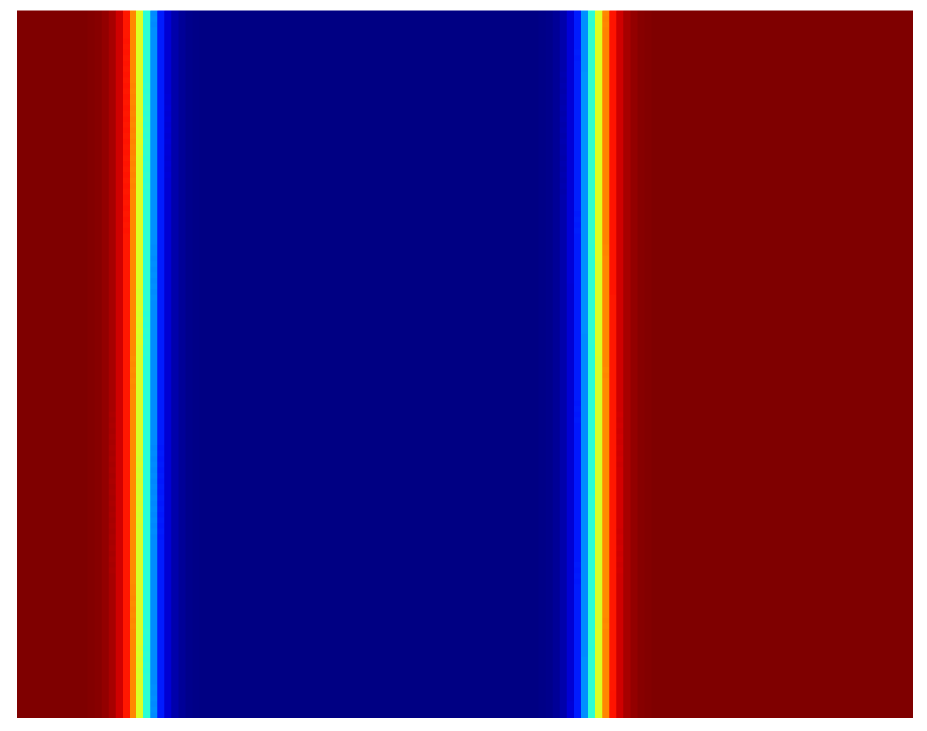}}
\centerline{\small (a) The fixed time step size $\tau=0.1$}
\vskip 1mm
\centerline{ \includegraphics[scale=0.16]{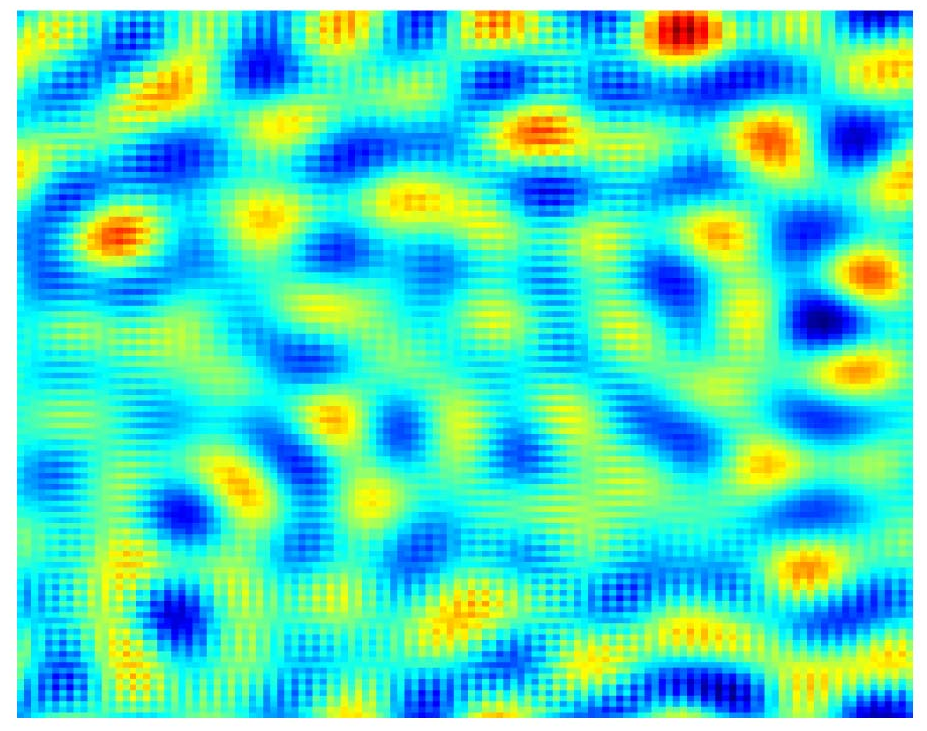}\includegraphics[scale=0.16]{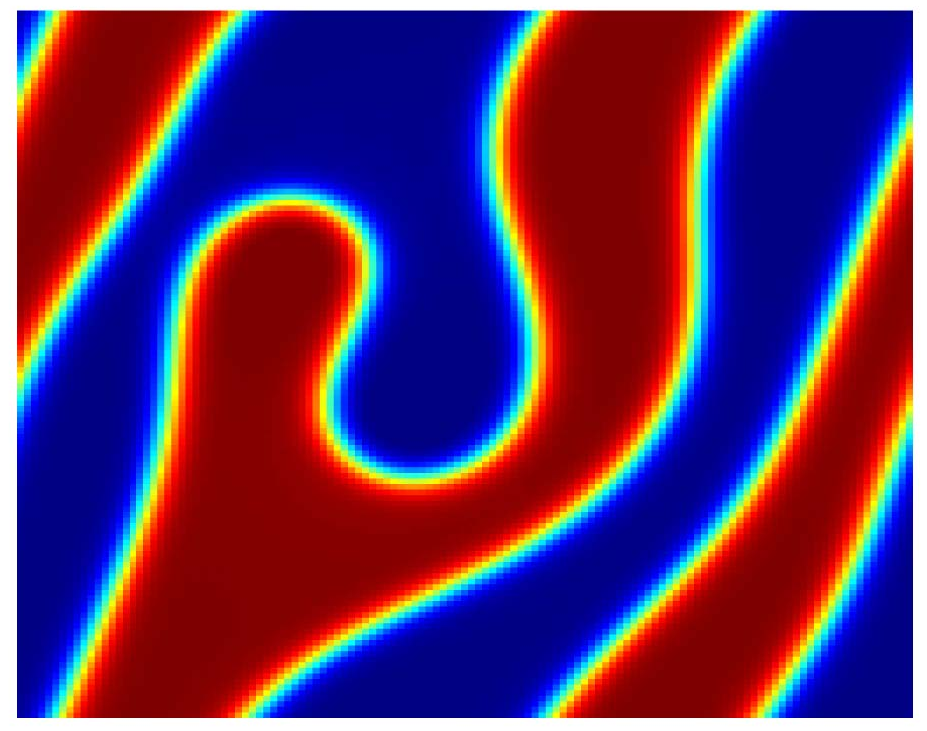}\includegraphics[scale=0.16]{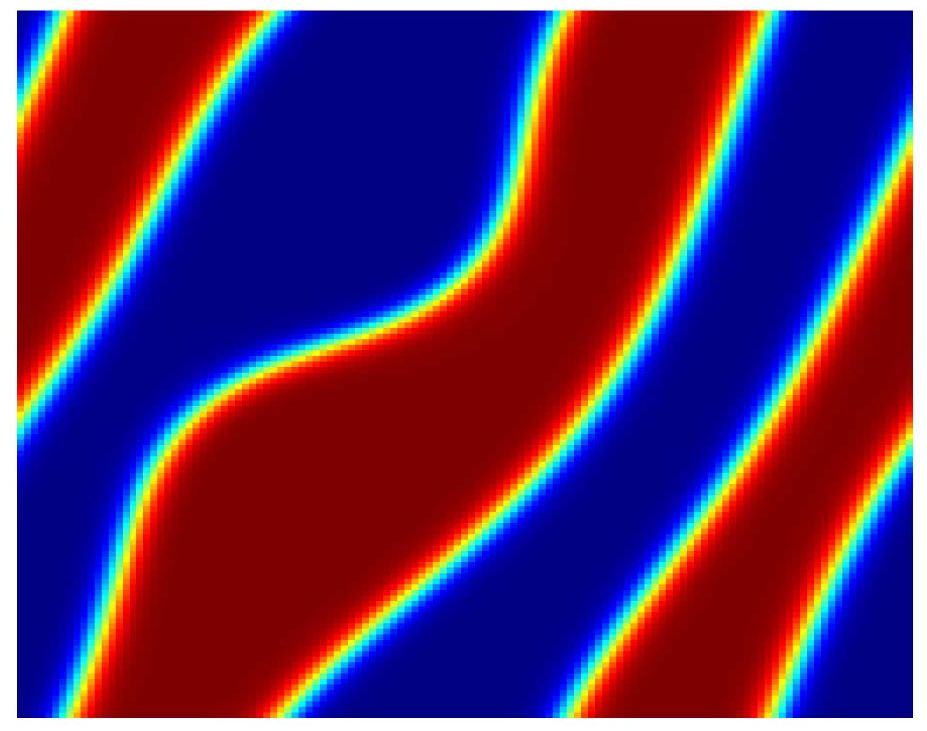}\includegraphics[scale=0.16]{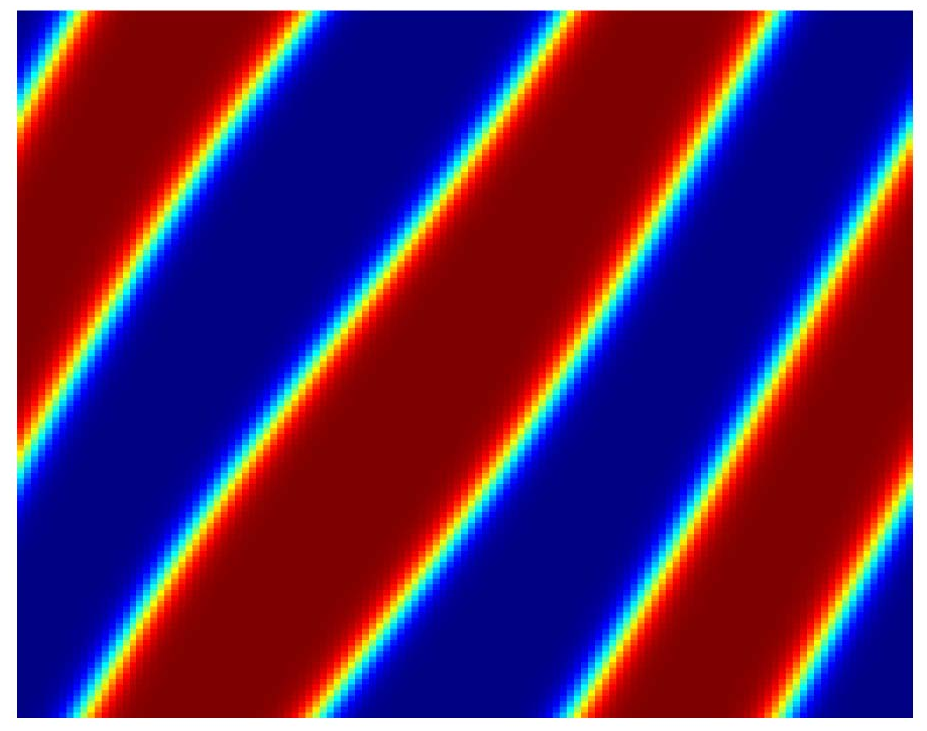}\includegraphics[scale=0.16]{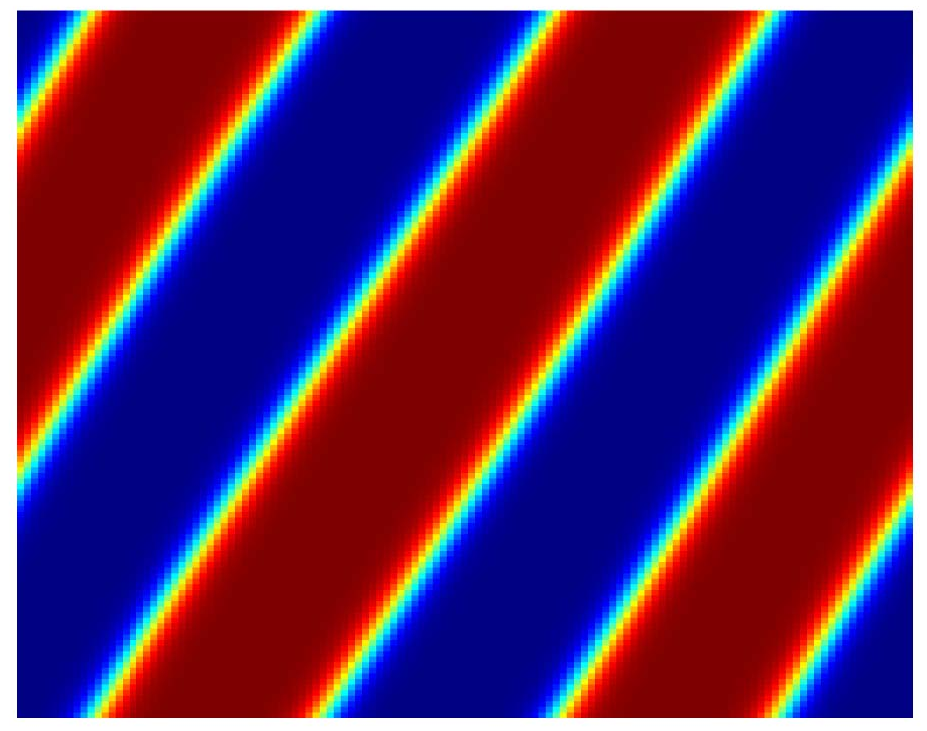}}
\centerline{\small (b)  The time-adaptive strategy  \eqref{adp}   with $\tau_{min}=10^{-4}, \tau_{max}=0.1$, and $\gamma=10^{5}$}
\vskip 1mm
\centerline{ \includegraphics[scale=0.16]{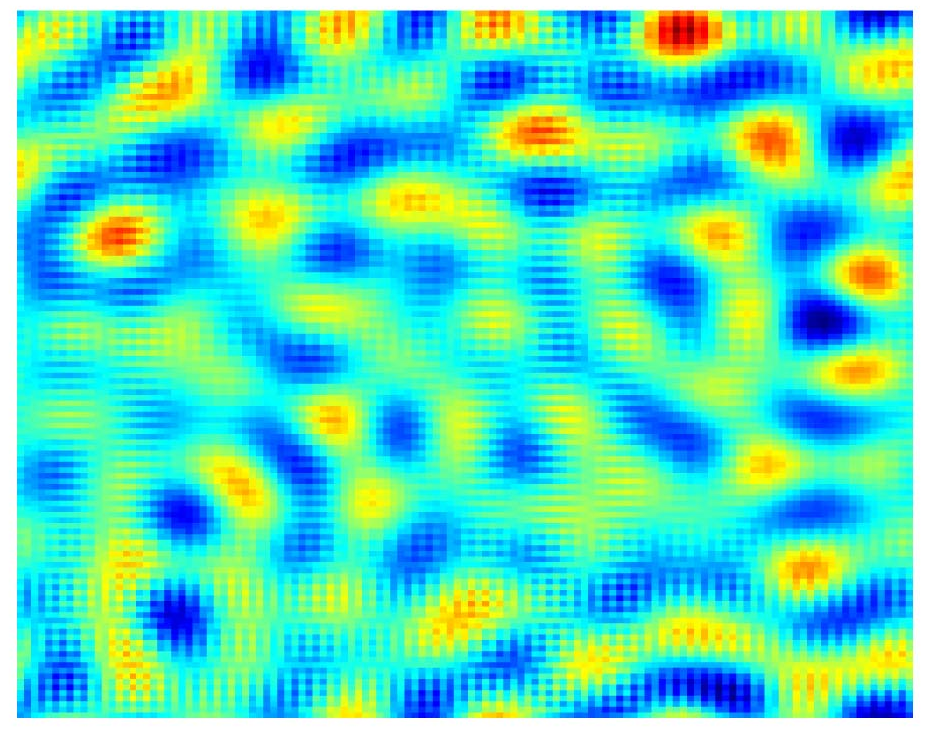}\includegraphics[scale=0.16]{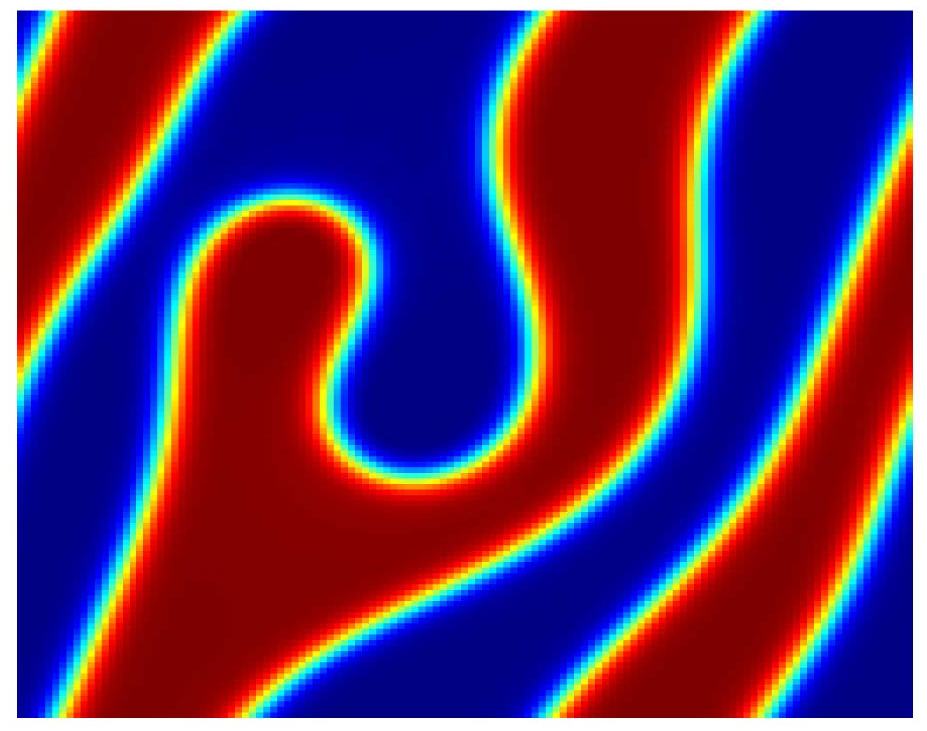}\includegraphics[scale=0.16]{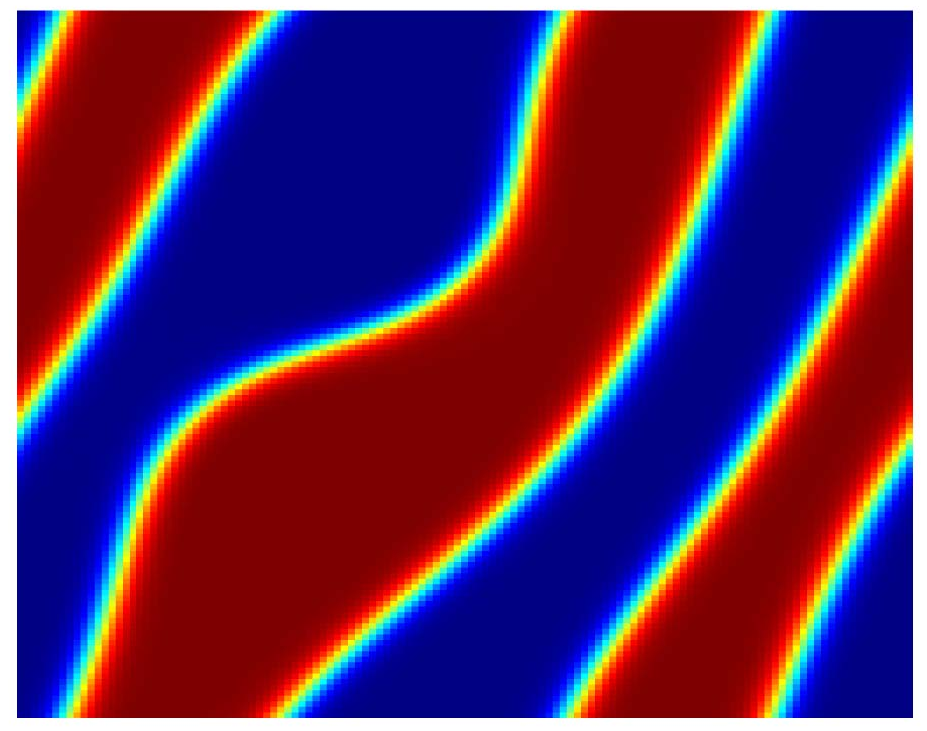}\includegraphics[scale=0.16]{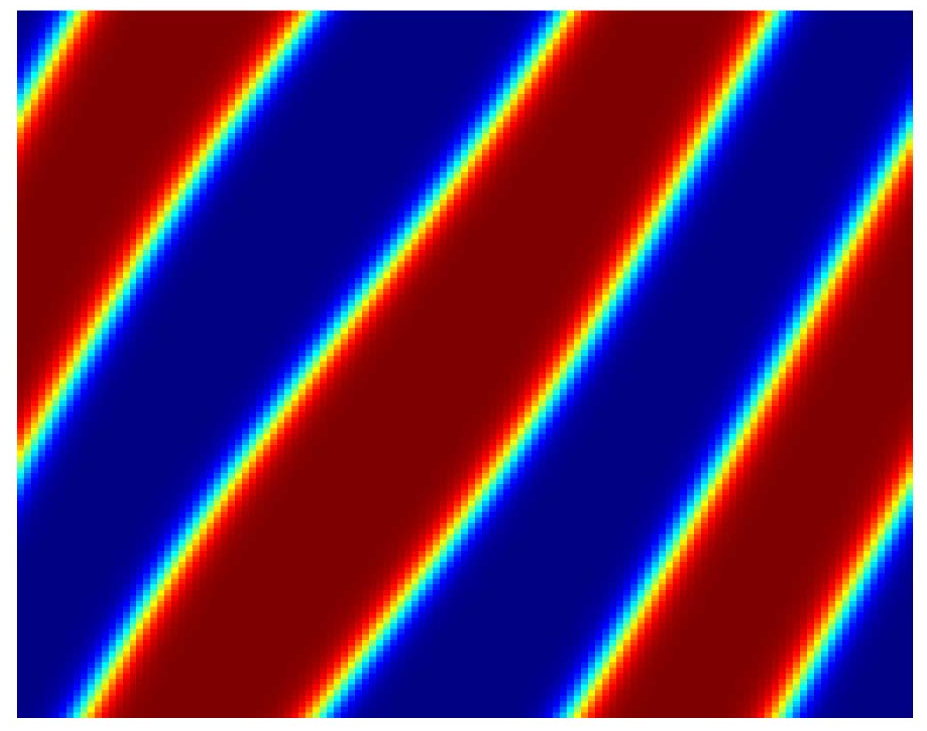}\includegraphics[scale=0.16]{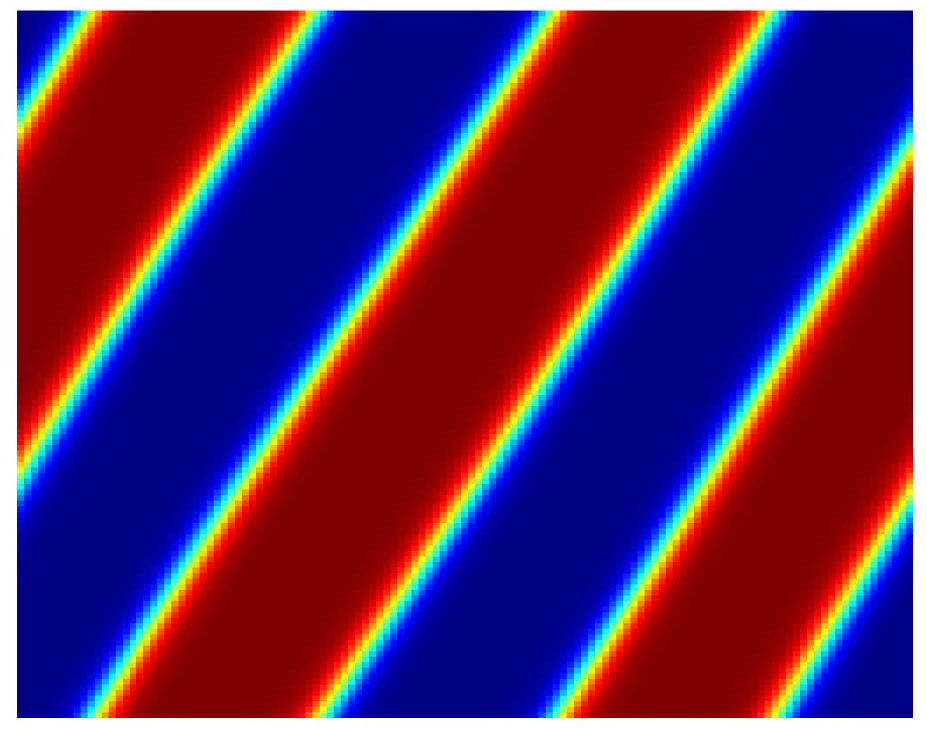}}
\centerline{\small (c) The fixed time step size $\tau=10^{-4}$}
\caption{Snapshots of the simulated coarsening process produced by the TDSR-ETD3 scheme \eqref{full-Ne}  for  the $H^{-1}$ gradient flow at $t=0.1, 10, 20, 100, 200$, respectively.
}\label{fig5}
\end{figure*}

\subsubsection{Coarsening dynamic by the Cahn-Hilliard equation}
\begin{figure}[!t]
\begin{minipage}[t]{0.49\linewidth}
\centerline{\includegraphics[scale=0.3]{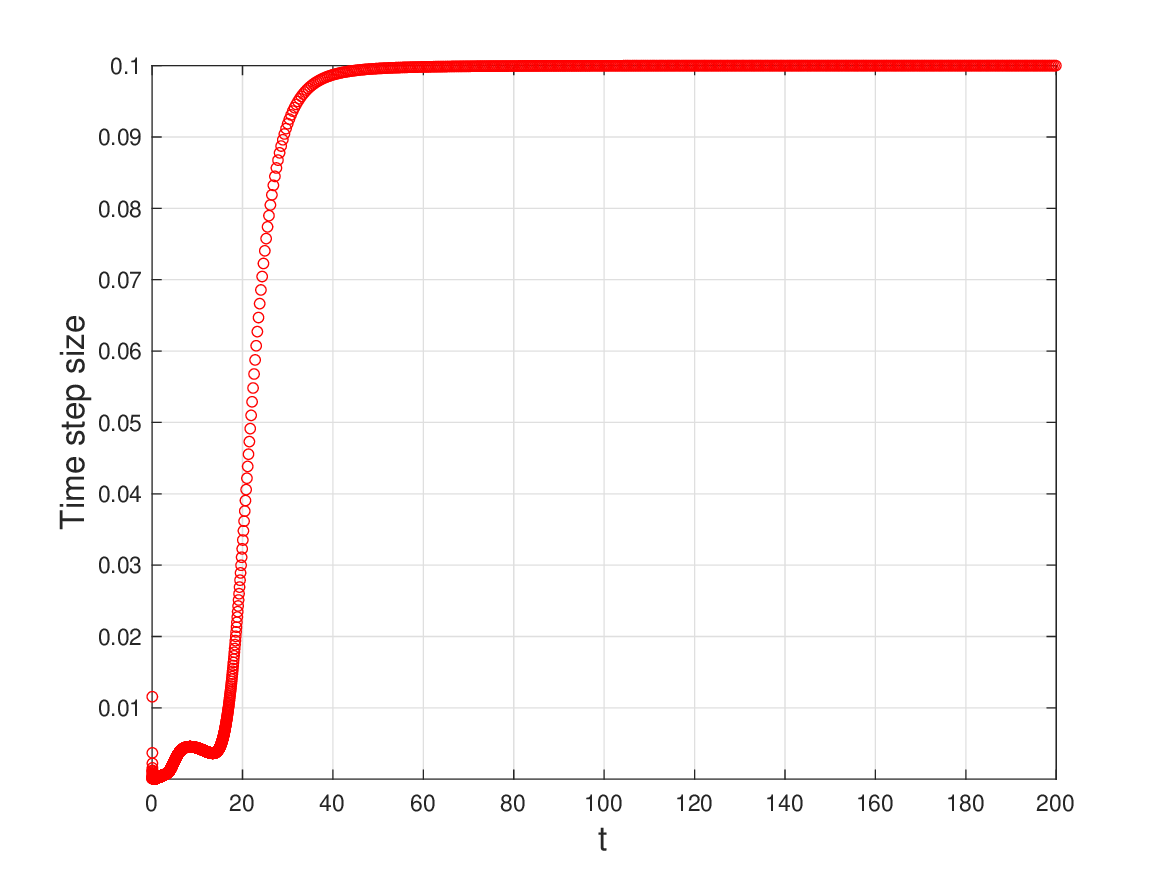}}
\centerline{\small (a) The adaptive time step size}
\end{minipage}
\begin{minipage}[t]{0.49\linewidth}
\centerline{\includegraphics[scale=0.3]{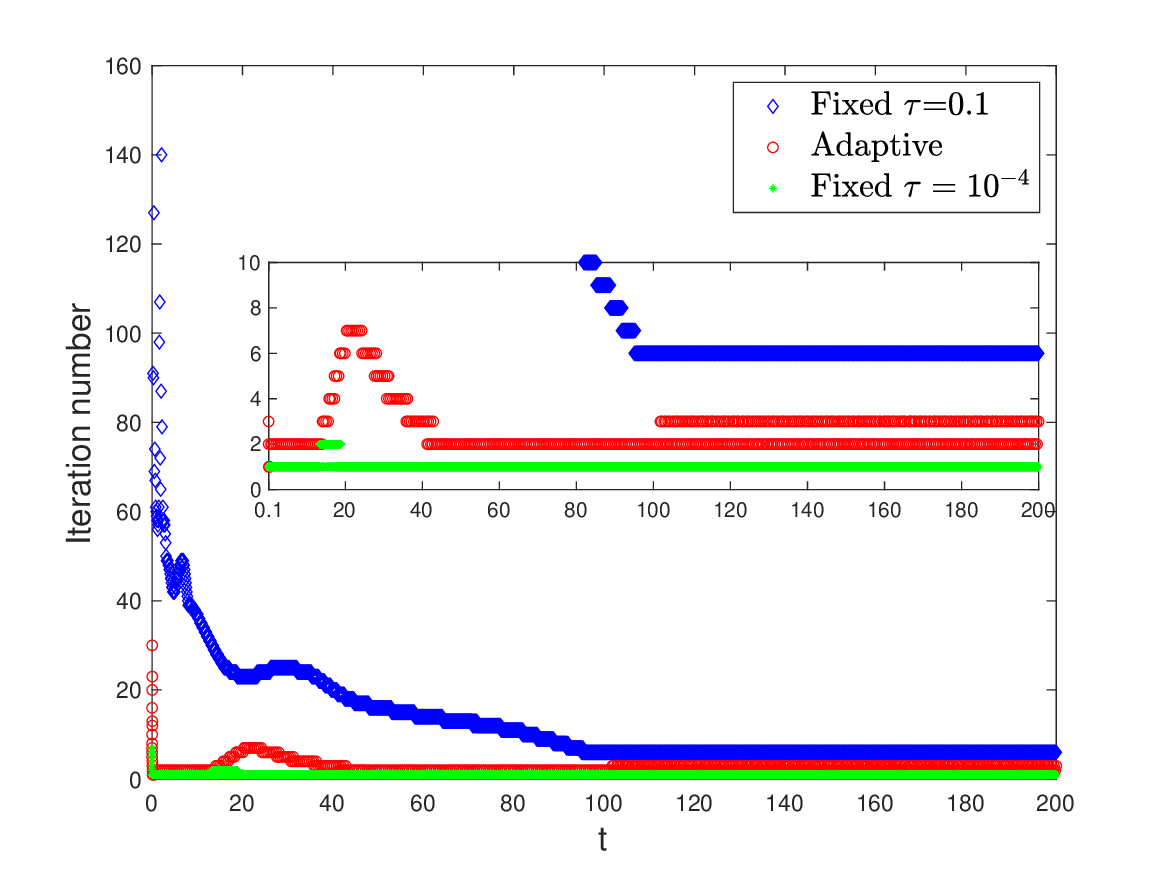}}
\centerline{\small (b) The number of Picard iterations}
\end{minipage}
\vskip -1mm
\begin{minipage}[t]{0.49\linewidth}
\centerline{\includegraphics[scale=0.3]{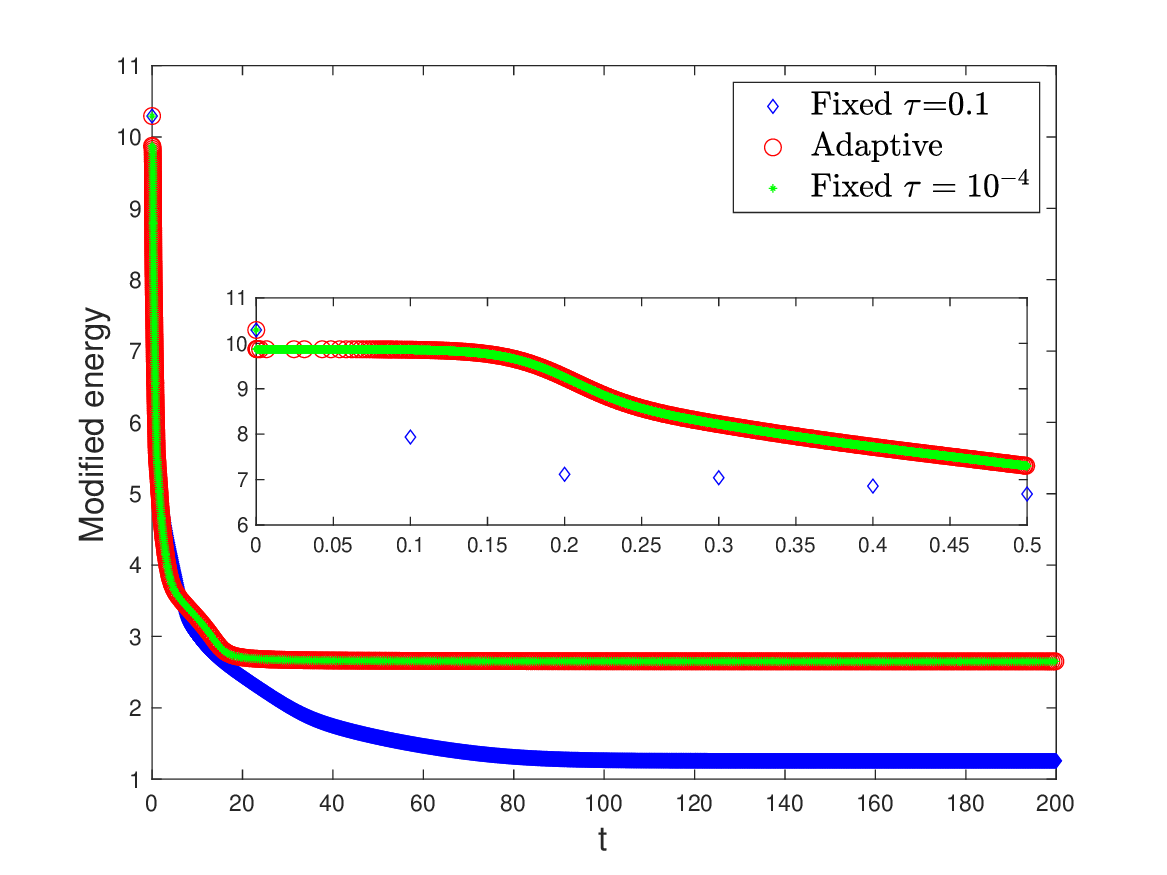}}
\centerline{\small (c) The modified energy}
\end{minipage}
\begin{minipage}[t]{0.49\linewidth}
\centerline{\includegraphics[scale=0.3]{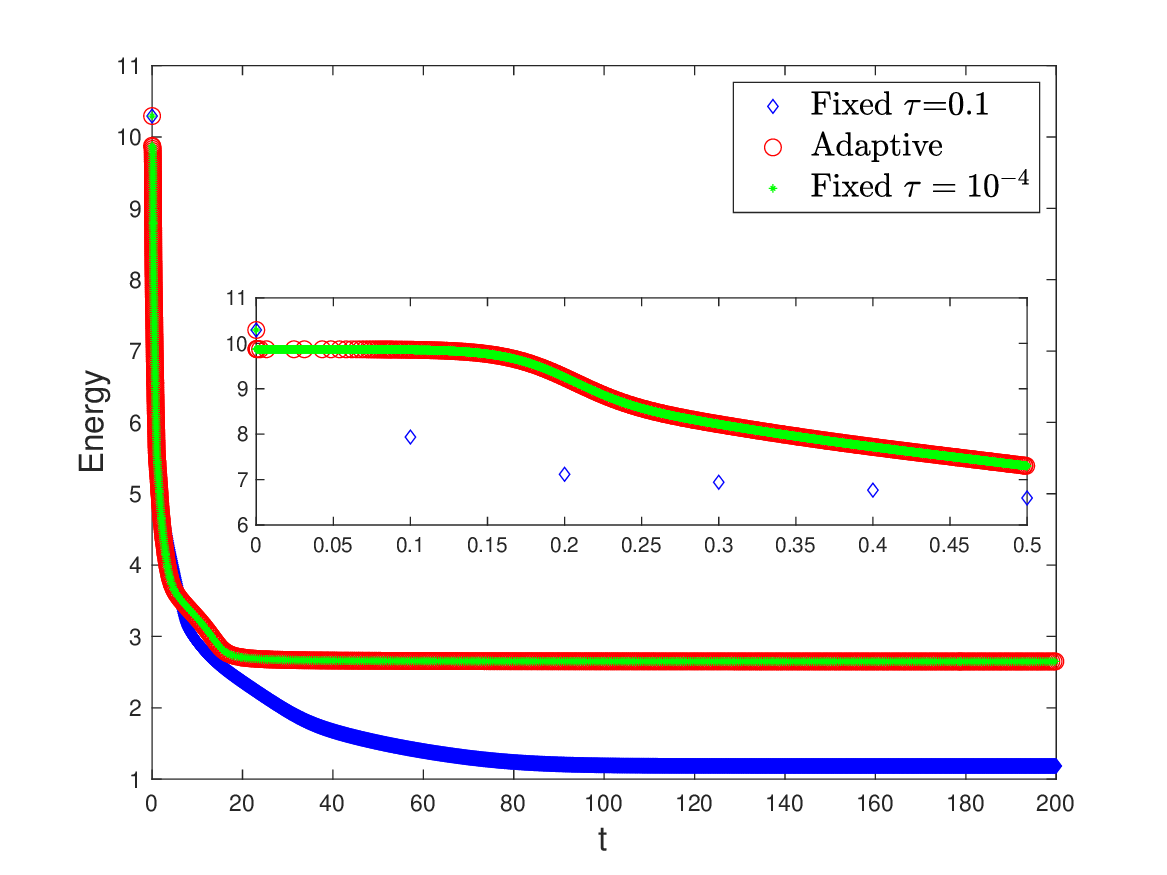}}
\centerline{\small (d) The original energy}
\end{minipage}
\caption{Evolutions of   (a) the  time step size by the time-adaptive strategy \eqref{adp}, (b) the number of Picard iterations, (c) the modified energy $E\big[\phi^{n}_{N}\big]+\theta \big[(R^{n})^2-1]$, and (d) the original energy  $E\big[\phi^{n}_{N}\big]$ produced by the TDSR-ETD3 scheme \eqref{full-Ne}  with  the time-adaptive strategy \eqref{adp}  for the $H^{-1}$ gradient flow.}\label{fig6}
\end{figure}

We first consider  a long-time coarsening dynamic process governed by  $H^{-1}$ gradient flow \eqref{prob} (i.e., the Cahn-Hilliard equation) with the double-well potential, subject to the periodic boundary condition with $\varepsilon^{2}=0.01$ and the domain $\Omega=(0,2\pi)^{2}$.
 The initial condition   is given by random data ranging from $-0.05$ to $0.05$ in the domain, i.e., $\phi(x,y,0) = 0.05\,(2\,\text{rand}(\cdot)-1)$ where
 $\text{rand}(\cdot)$ randomly generates a number between 0 and 1.
The simulation of the coarsening dynamics  is performed by the TDSR-ETD3 scheme with $128\times128$ Fourier modes.
Figure \ref{fig5} presents a comparison on the numerical solutions (up to $T=200$) produced by
the uniform time steps with $\tau=0.1$, the time-adaptive strategy \eqref{adp} with parameters $\Dt_{min}=10^{-4},\Dt_{max}=0.1,$ and $\gamma=10^{5}$, and the uniform time stepping  with $\tau=10^{-4}$.
It is observed that TDSR-ETD3 with  the  large time step size  $\tau=0.1$  leads to an inaccurate phase transition process, while
the time-adaptive strategy gives the correct coarsening process which matches very well with the results produced by  the small time step size  $\Dt=10^{-4}$.
The evolution of the time step size used by TDSR-ETD3 with the time-adaptive strategy \eqref{adp} is shown in Figure \ref{fig6}-(a).
We observe that  the large time step size $\tau=0.1$ is chosen during the time interval $[40,200]$, which clearly demonstrates the efficiency of TDSR-ETD3 with the time-adaptive strategy.
The evolutions of
the number of Picard iterations at each time step are  displayed in Figure \ref{fig6}-(b) for all the above three different types of time stepping schemes.
It shows that TDSR-ETD3 with the time-adaptive strategy also significantly  reduces the number of Picard iterations at each time step compared to the case of fixed large time step size $\Dt=0.1$ while still maintaining  the same accuracy as the case of the fixed small time step size $\Dt=10^{-4}$. Except for a very short initial time period and the time interval $[16,40]$, the number of Picard iterations is only either 2 or 3. In addition, the evolutions of the simulated modified and original free energies are plotted in Figure \ref{fig6}-(c) and (d), respectively, which show they are almost identical, and the time-adaptive  strategy and  the uniform time stepping  with $\tau=10^{-4}$ both give accurate results.

\subsubsection{Thin film epitaxial growth}
\begin{figure*}[!t]
\begin{minipage}[t]{0.99\linewidth}
\centerline{\includegraphics[scale=0.3]{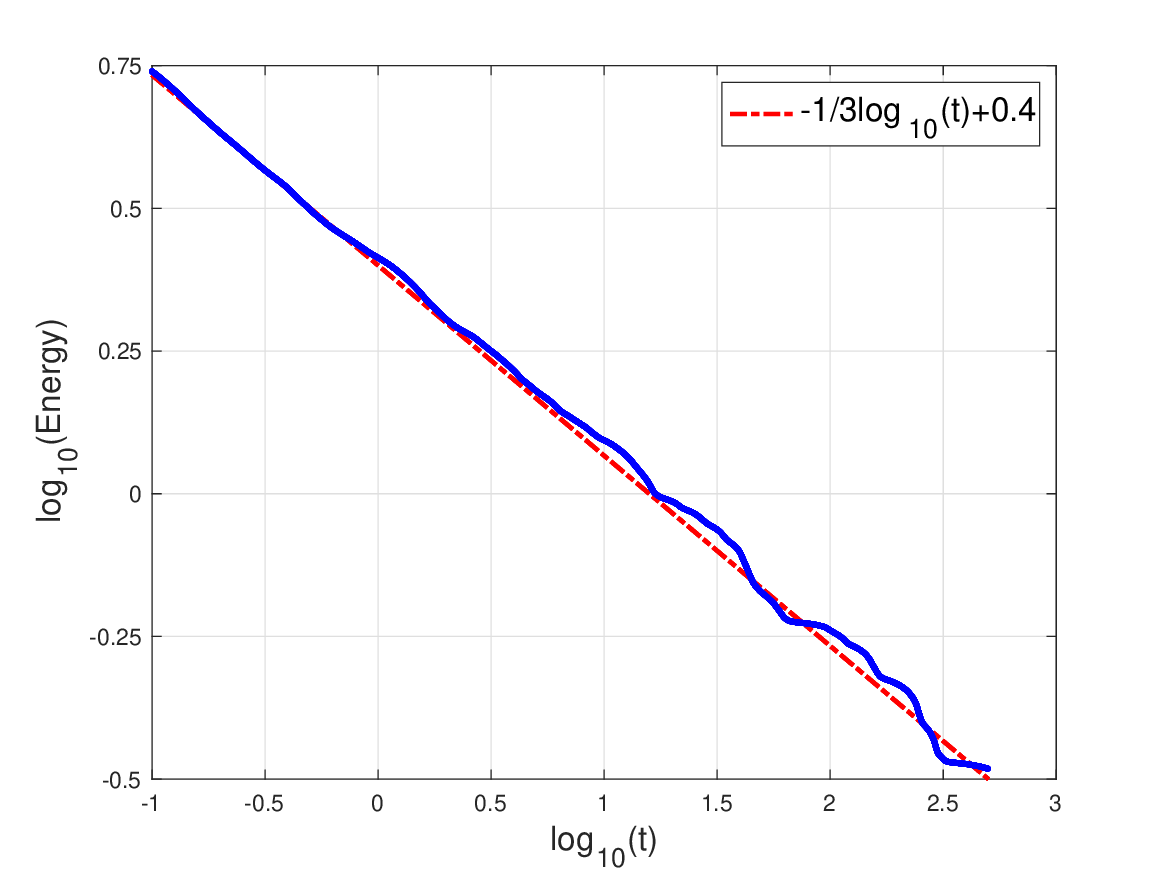}\includegraphics[scale=0.3]{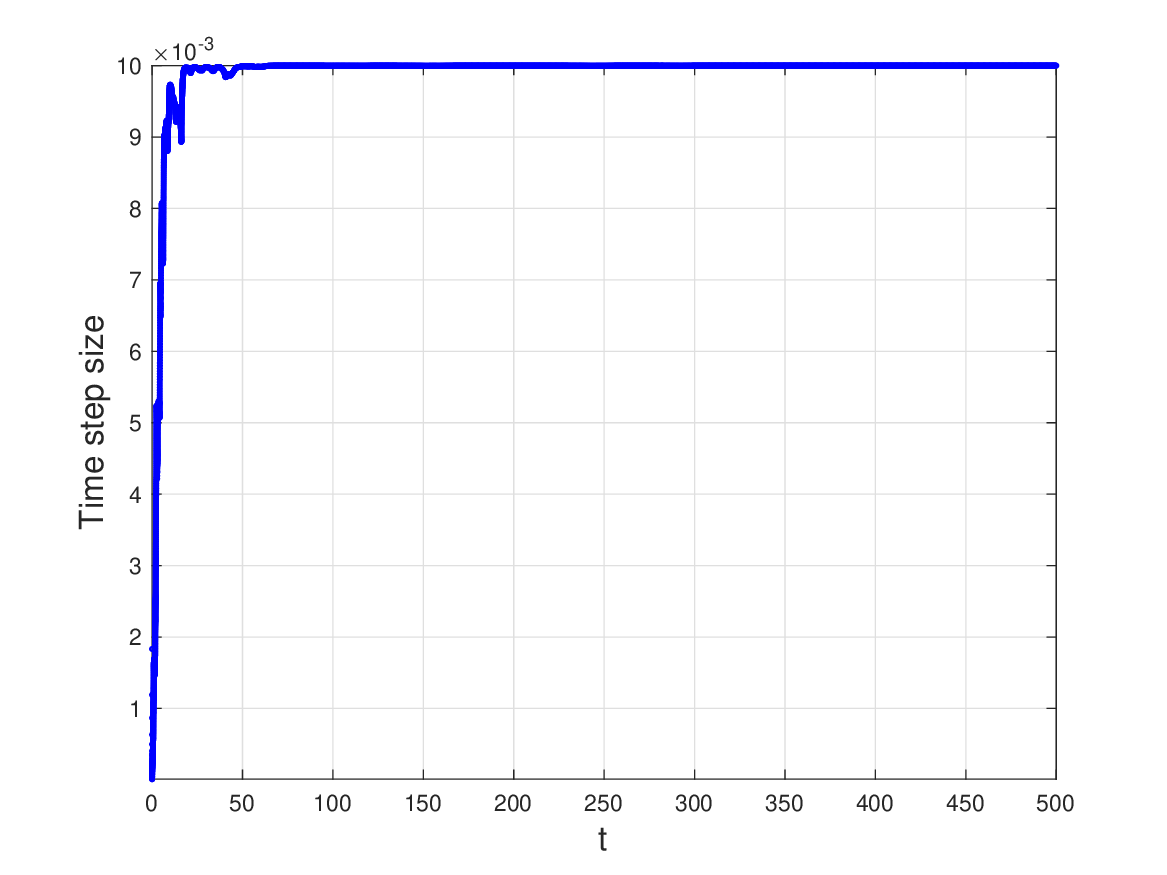}}
\centerline{\small (a) With slope selection. Left: energy, right: time step size.}
\end{minipage}\\
\begin{minipage}[t]{0.99\linewidth}
\centerline{\includegraphics[scale=0.3]{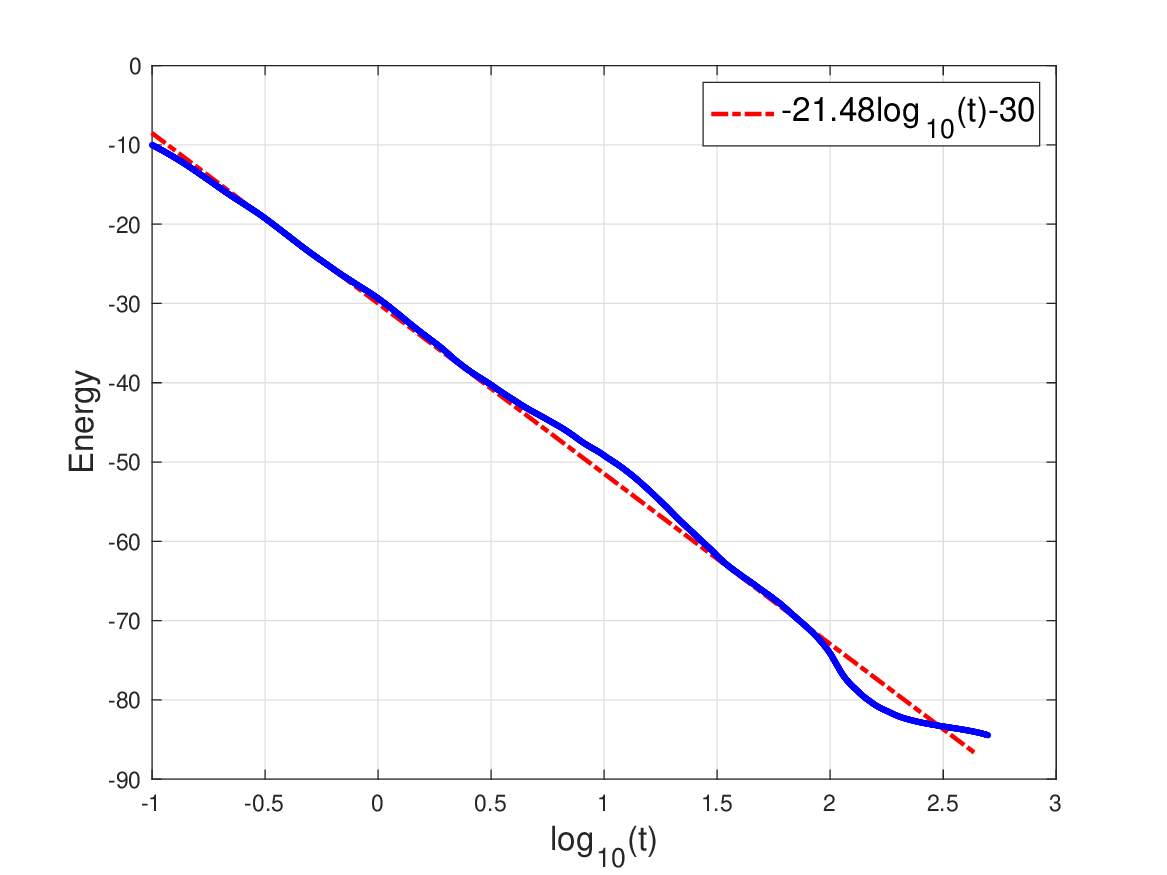}\includegraphics[scale=0.3]{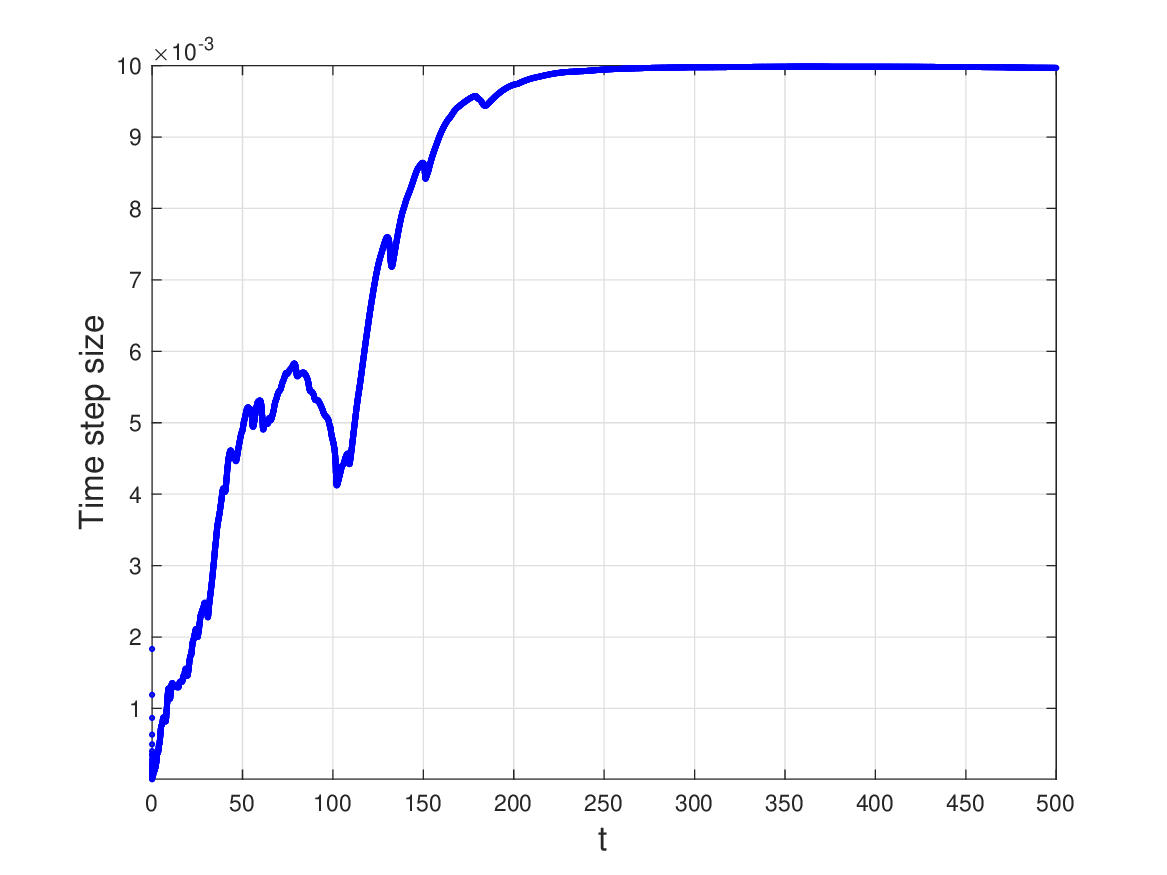} }
\centerline{\small (b) Without slope selection. Left: energy, right: time step size}
\end{minipage}
\caption{Evolutions of the simulated energy  and the time step size  by the TDSR-ETD3 scheme \eqref{full-Ne}  with  the time-adaptive strategy \eqref{adp}  for  the MBE model.
}\label{fig7}
\end{figure*}

\begin{figure*}[!t]
\centerline{\includegraphics[scale=0.22]{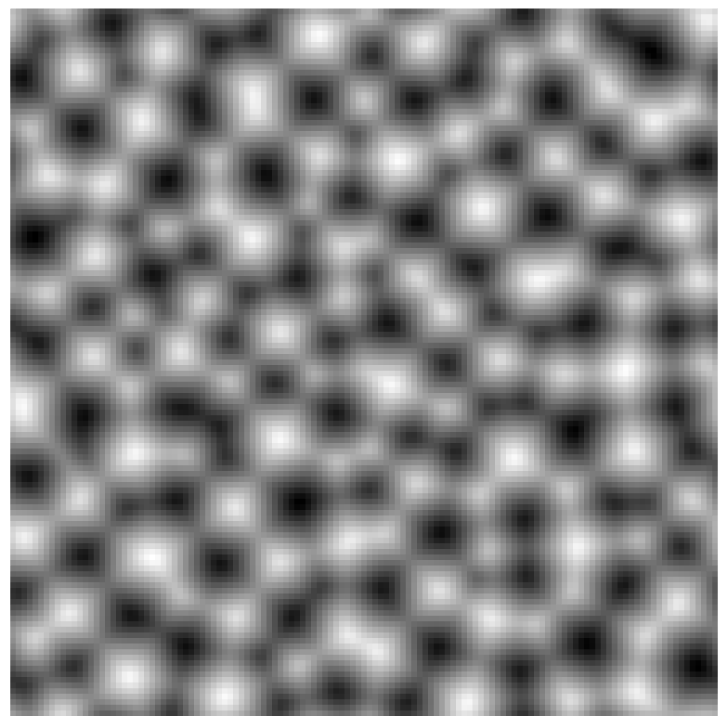}\includegraphics[scale=0.22]{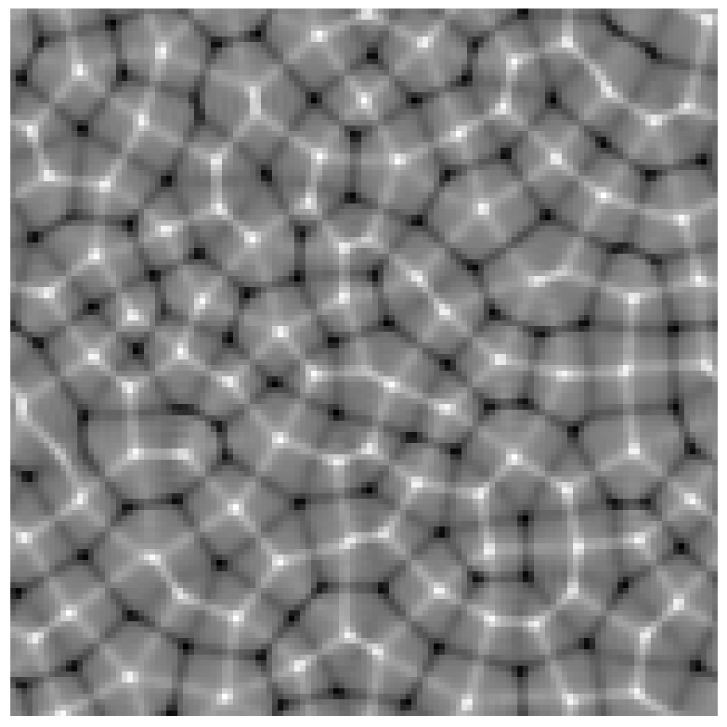}\hspace{0.5cm}\includegraphics[scale=0.22]{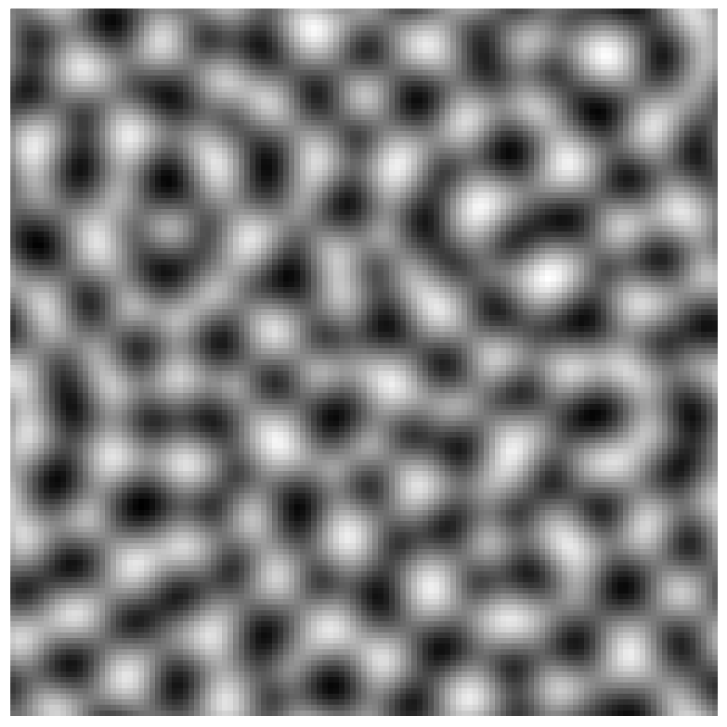}\includegraphics[scale=0.22]{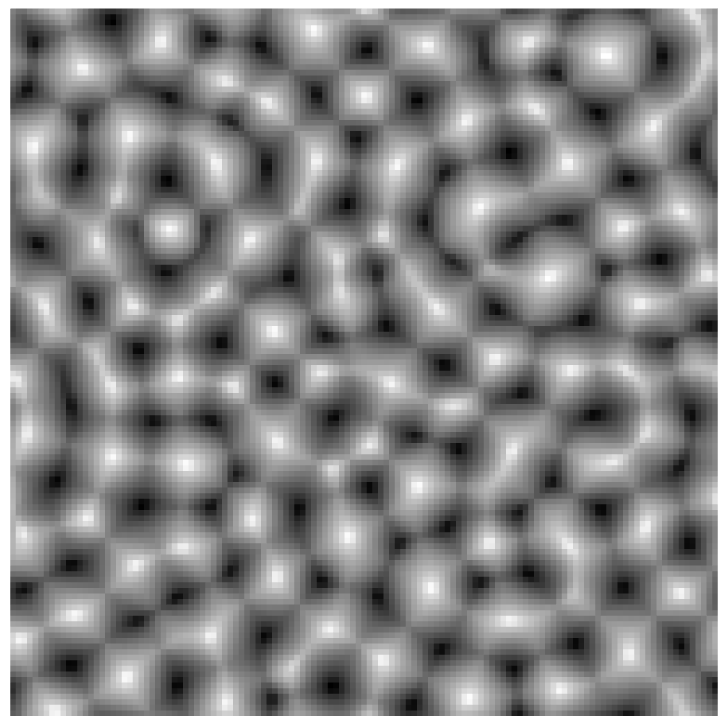}}
\centerline{\includegraphics[scale=0.22]{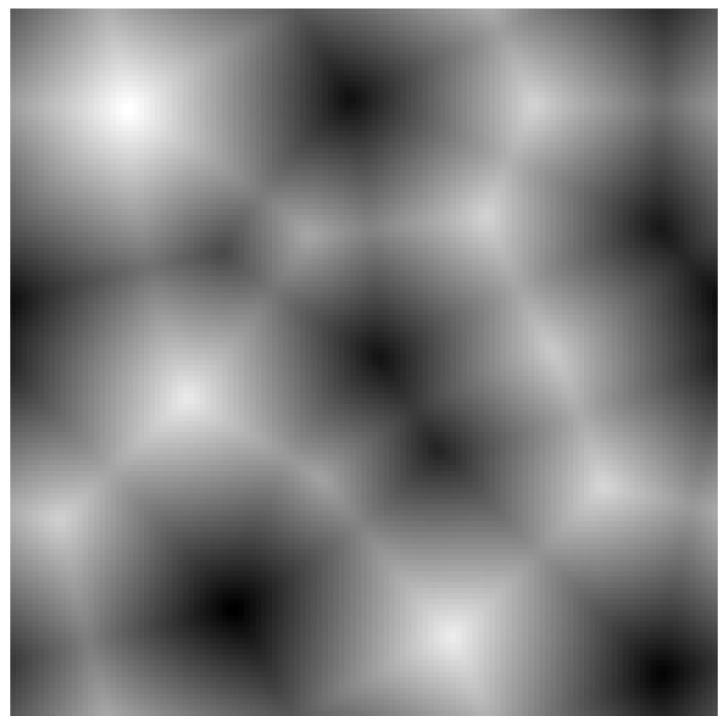}\includegraphics[scale=0.22]{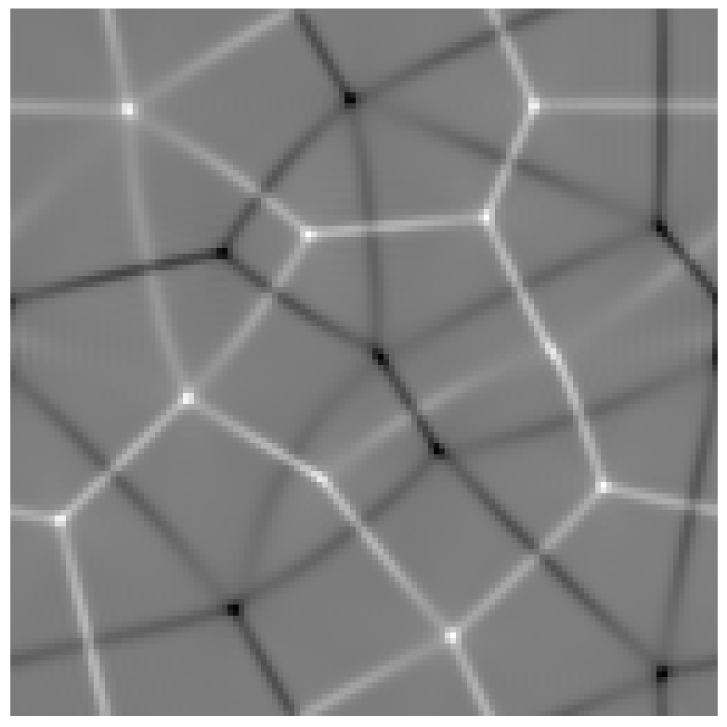}\hspace{0.5cm}\includegraphics[scale=0.22]{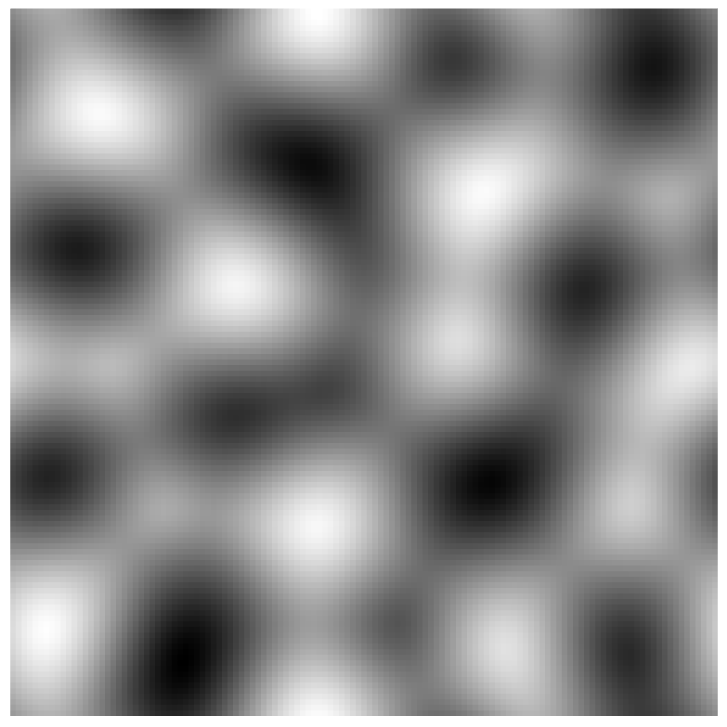}\includegraphics[scale=0.22]{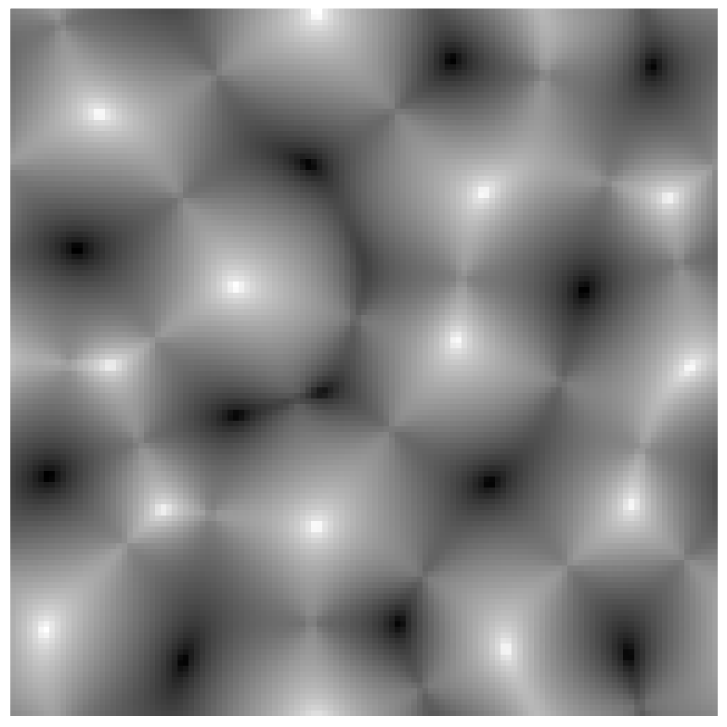}}
\centerline{\includegraphics[scale=0.22]{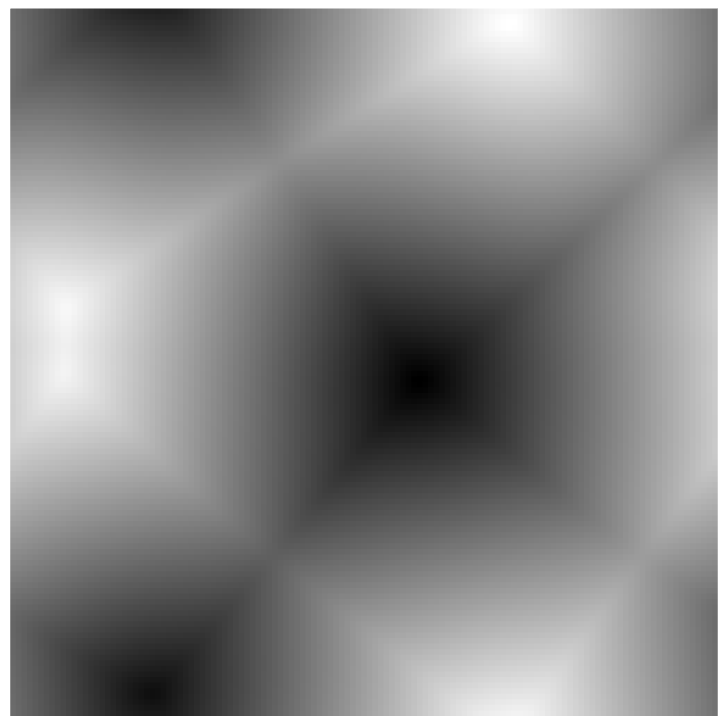}\includegraphics[scale=0.22]{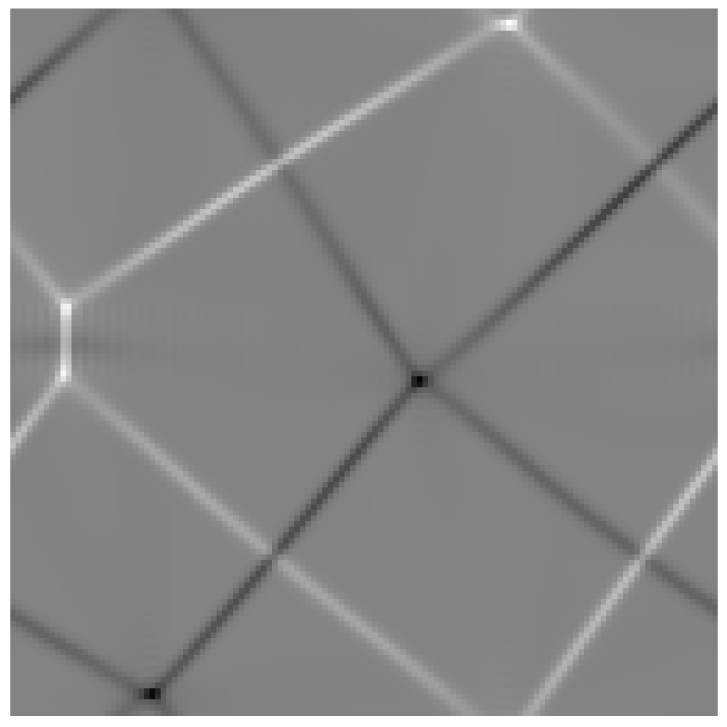}\hspace{0.5cm}\includegraphics[scale=0.22]{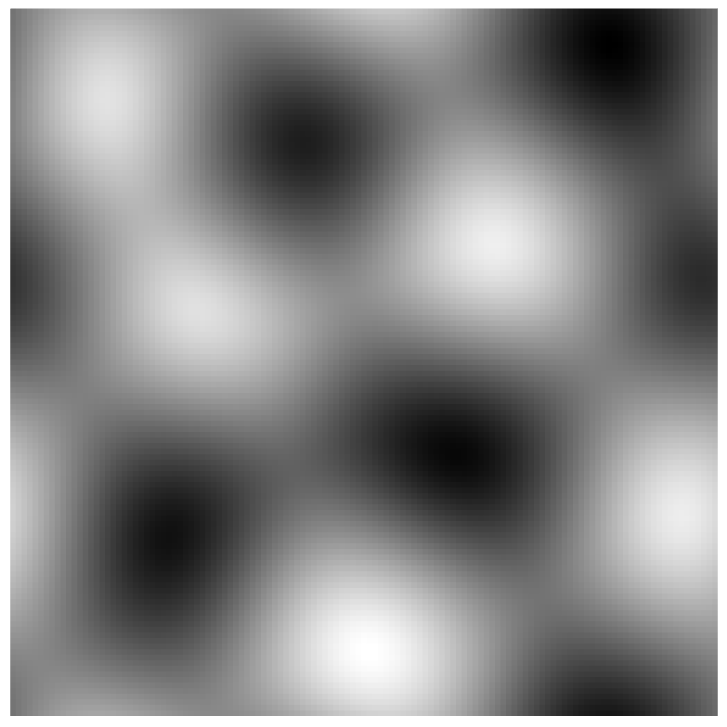}\includegraphics[scale=0.22]{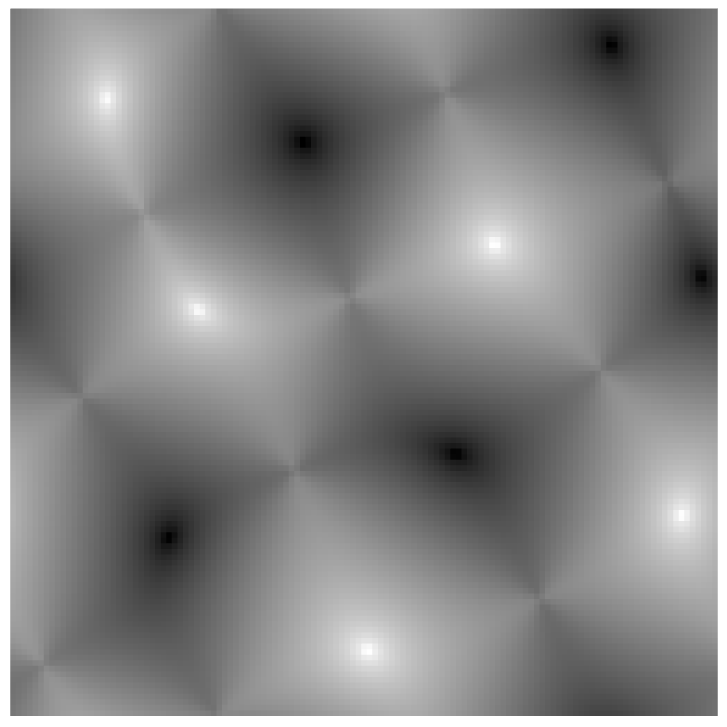}}
\centerline{\includegraphics[scale=0.22]{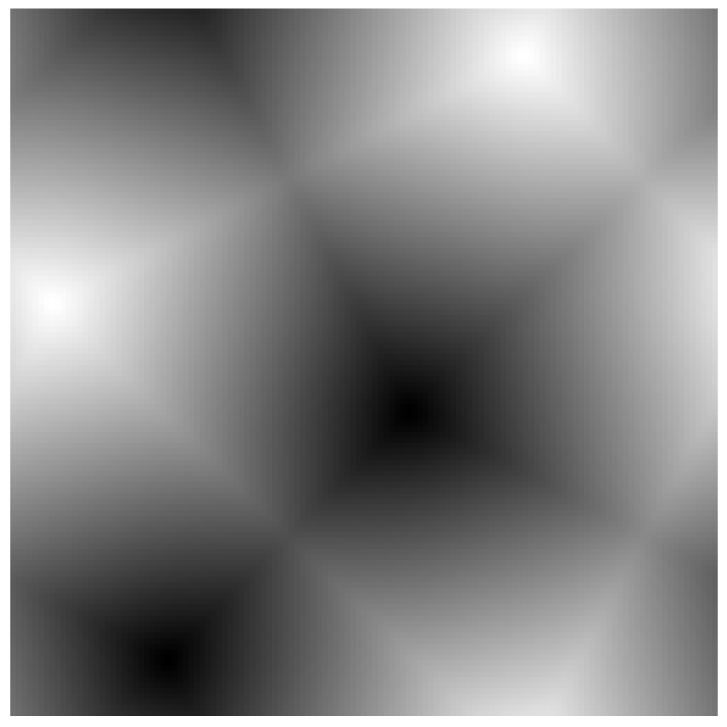}\includegraphics[scale=0.22]{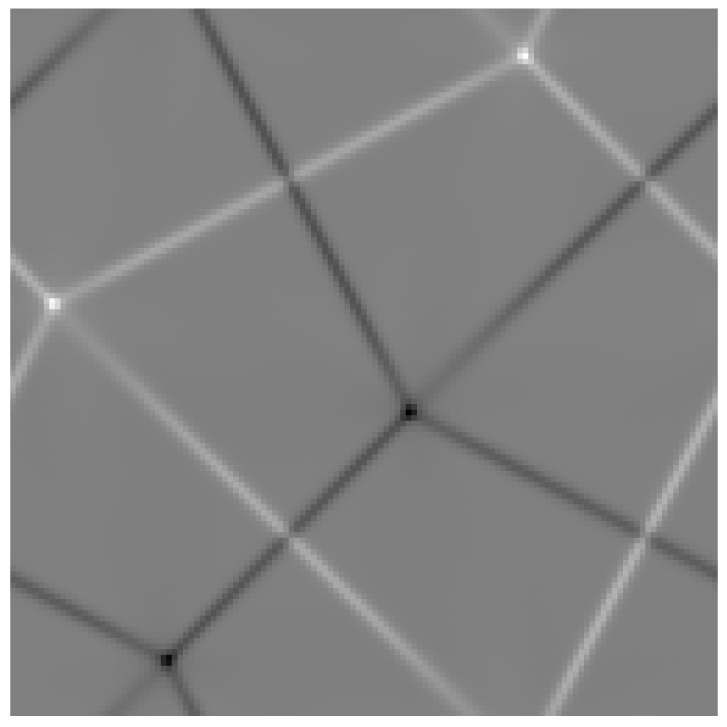}\hspace{0.5cm}\includegraphics[scale=0.22]{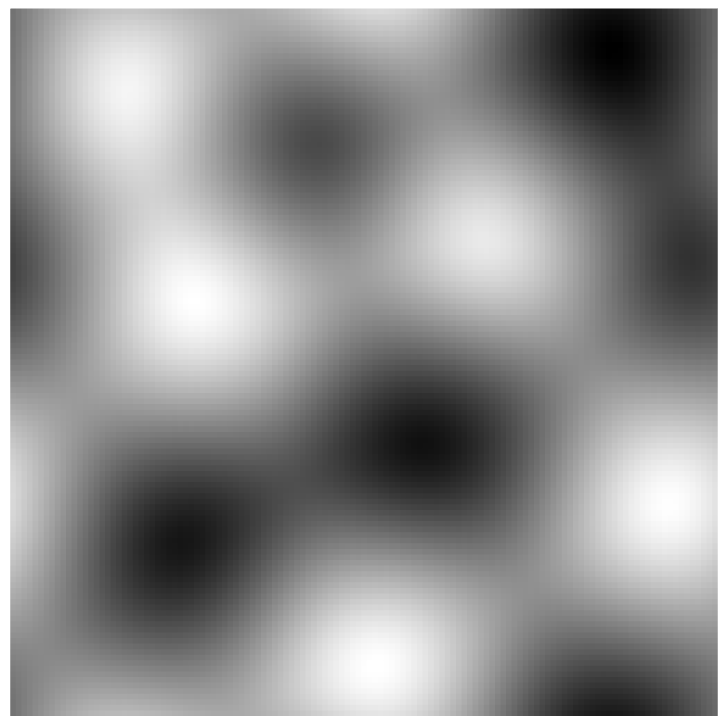}\includegraphics[scale=0.22]{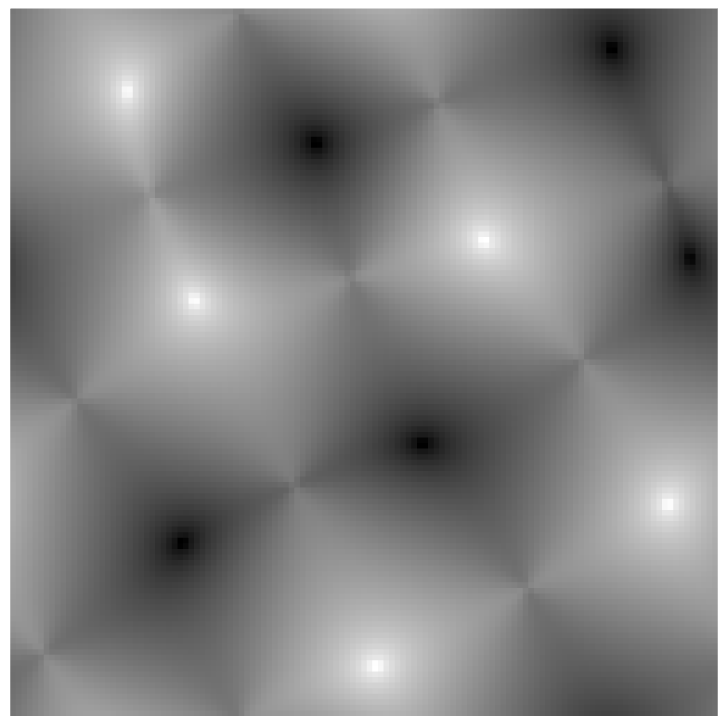}}
\caption{Isolines of the simulated height function $\phi$  and its
Laplacian $\Delta \phi$  for the MBE model with slope selection (left two columns) and without slope selection (right two columns) around the times $t=1$, $50$, $100$, and $500$ (from top to bottom).}
\label{fig8}
\end{figure*}

We next numerically investigate the power laws of the height growth function $\phi$ governed by the MBE model \eqref{MBE} with slope selection or without slope selection.
The domain  and the width parameter are set to be $\Omega=(0,2\pi)^{2}$ and $\varepsilon =0.03$, respectively. The initial condition is a random state with values ranging from $-0.001$ to $0.001$, i.e., $\phi(x,y,0) = 0.001(2\,\text{rand}(\cdot)-1)$.
The thin film epitaxial growth process usually requires a long time simulation before reaching a steady state, and the free energy \eqref{EF2} of the MBE model usually undergoes some very fast and slow rates of change during the process. Therefore, robust, high-order, and energy stable numerical schemes are crucial for such simulation, in order to guarantee the accuracy  and easily be  employed with time-adaptive strategies to speed up the computations.
We take the TDSR-ETD3 scheme \eqref{full-ds-mbe} and  the time-adaptive strategy \eqref{adp} with $\Dt_{min}=10^{-5},\Dt_{max}=10^{-2}$ and $\gamma=100$. The Fourier spectral method with $128\times128$ basis modes is used for spatial discretization.
Figure \ref{fig7} presents the evolutions of the energy and the time step size for the models with slope section and without slope selection, respectively.
It is observed that the energy of the MBE model with slope selection  dissipates closely at a level of $O(t^{-1/3})$ and the one of the case without slope selection is approximately $O(-\log_{10}(t)),$ which are consistent well with those reported in \cite{QZT11,Xu06}.
We also see that a large time step size $\Dt=10^{-2}$ is automatically chosen by the time-adaptive strategy in the later period of the simulations.
In Figure \ref{fig8}, we display some  snapshots of the height function $\phi$
and its Laplacian $\Delta\phi$ for both models obtained around the times $t=1, 50, 100,$ and $500$.

\subsubsection{Crystal growth in a supercooled liquid}
\begin{figure*}[!t]
\centerline{\includegraphics[scale=0.3]{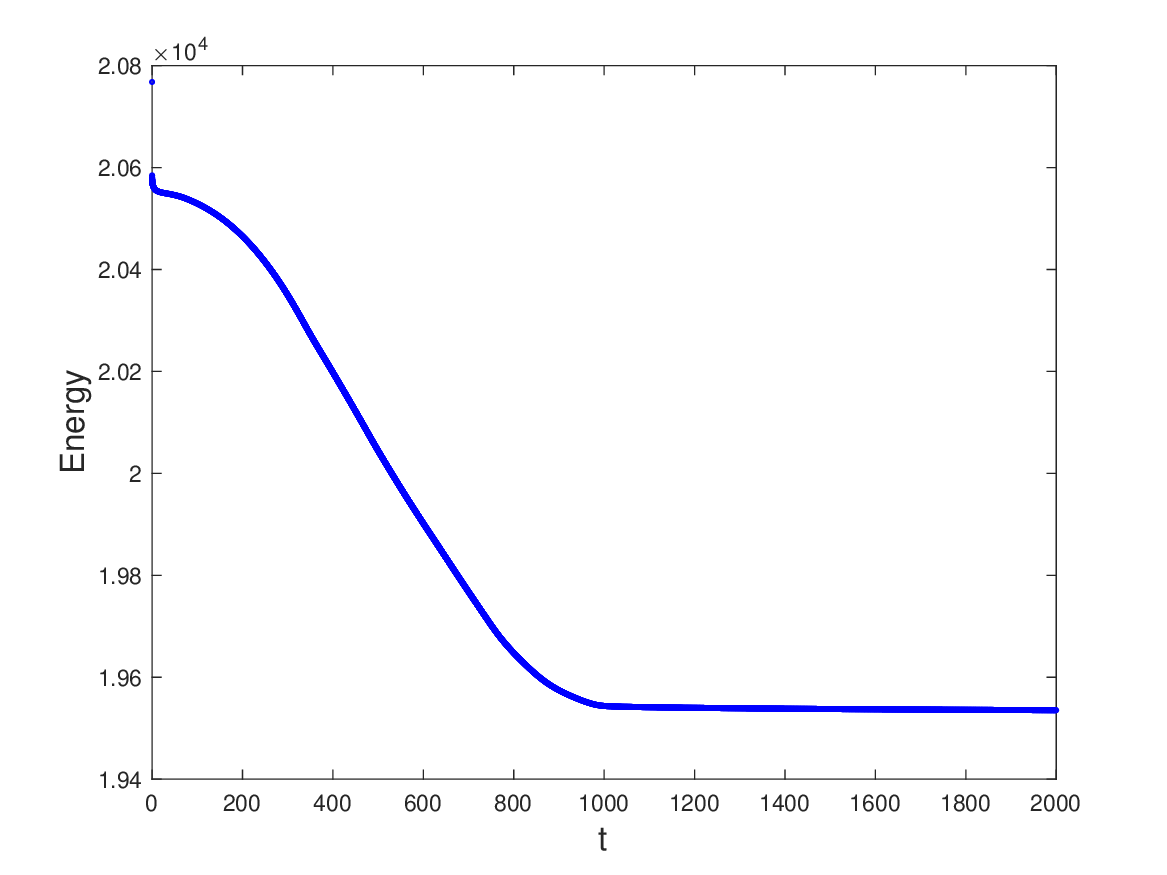}\includegraphics[scale=0.3]{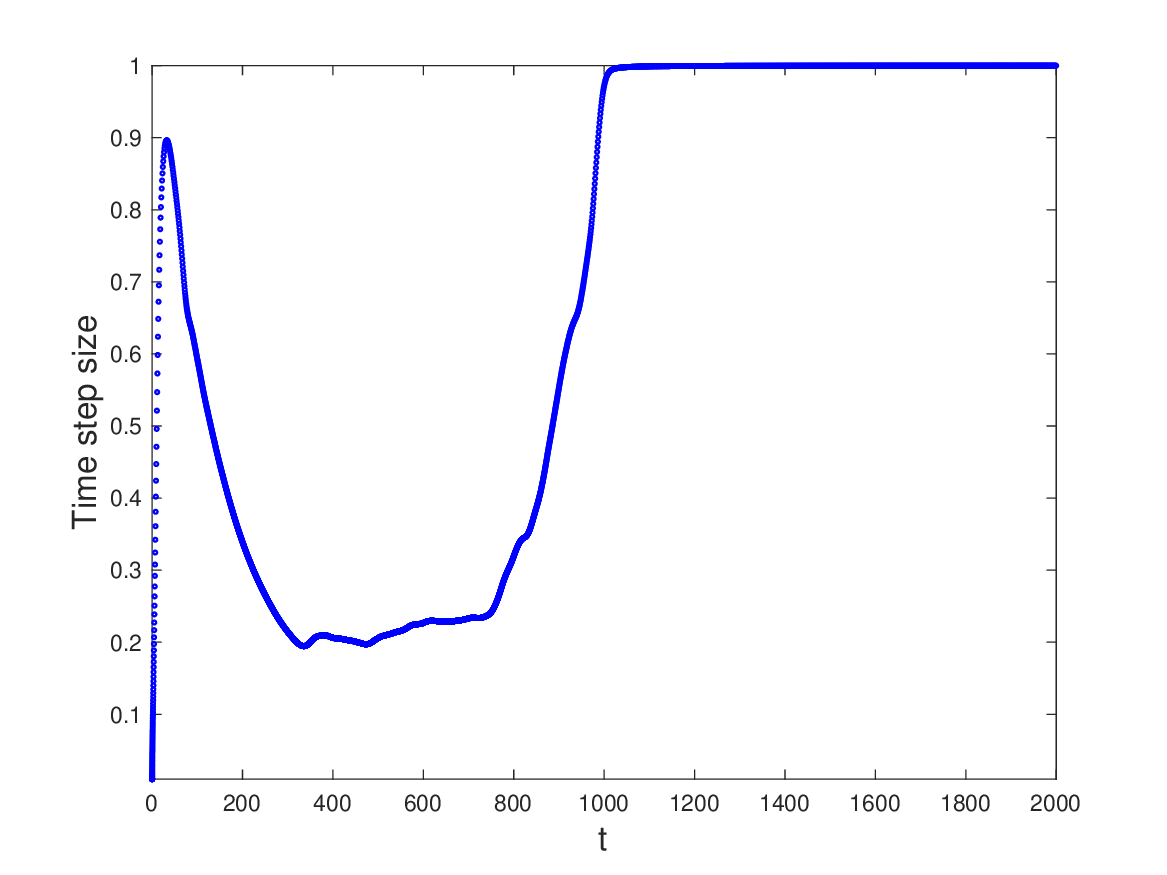}}
\vspace{-0.2cm}
\caption{Evolutions of the simulated energy (left) and the time step size (right) by the TDSR-ETD3 scheme \eqref{full-ds-pfc} with the time-adaptive strategy \eqref{adp} for the growth of three crystallites governed by the 2D PFC model.
 }\label{fig11}
\end{figure*}

\begin{figure*}[!ht]
\centerline{ \includegraphics[scale=0.22]{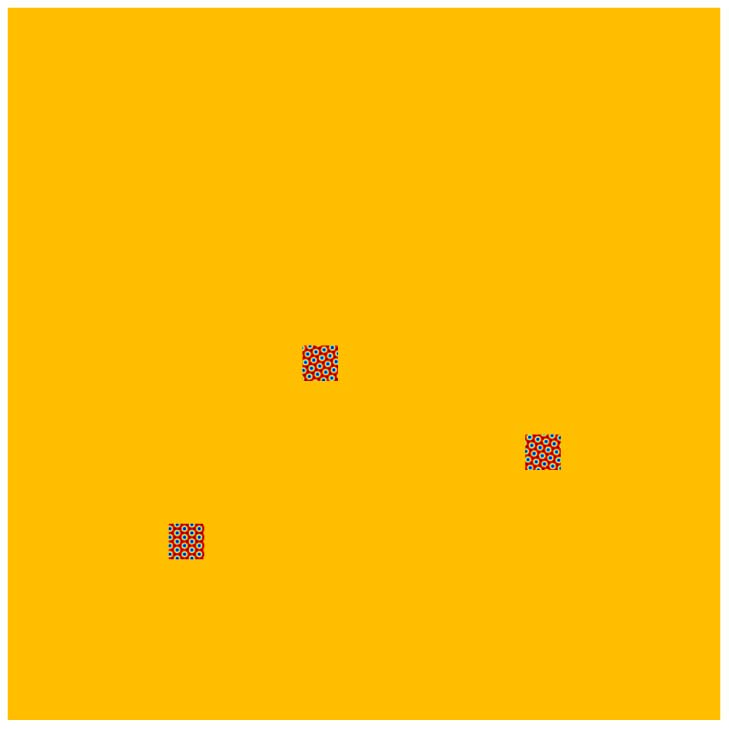} \includegraphics[scale=0.22]{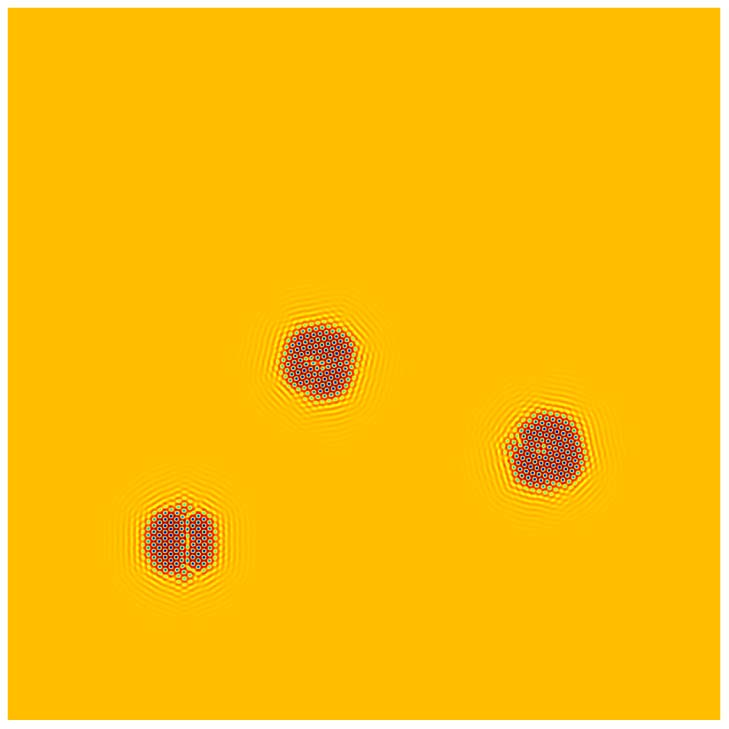}  \includegraphics[scale=0.22]{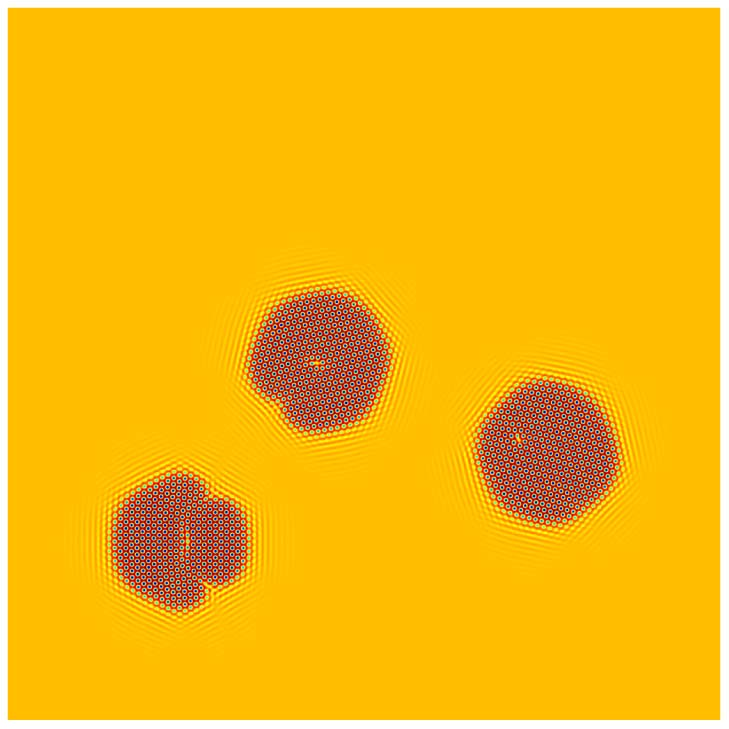}   \includegraphics[scale=0.22]{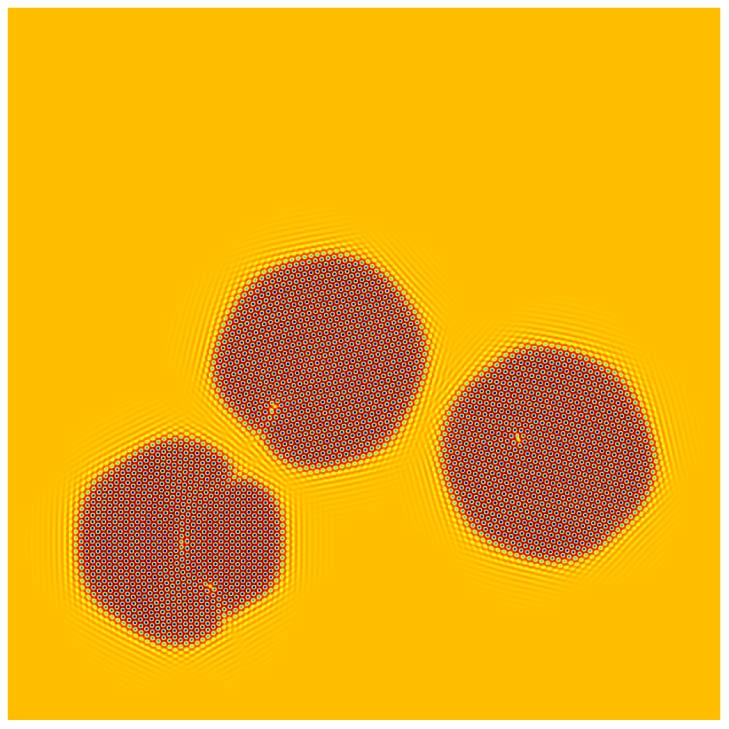}}
\centerline{ \includegraphics[scale=0.22]{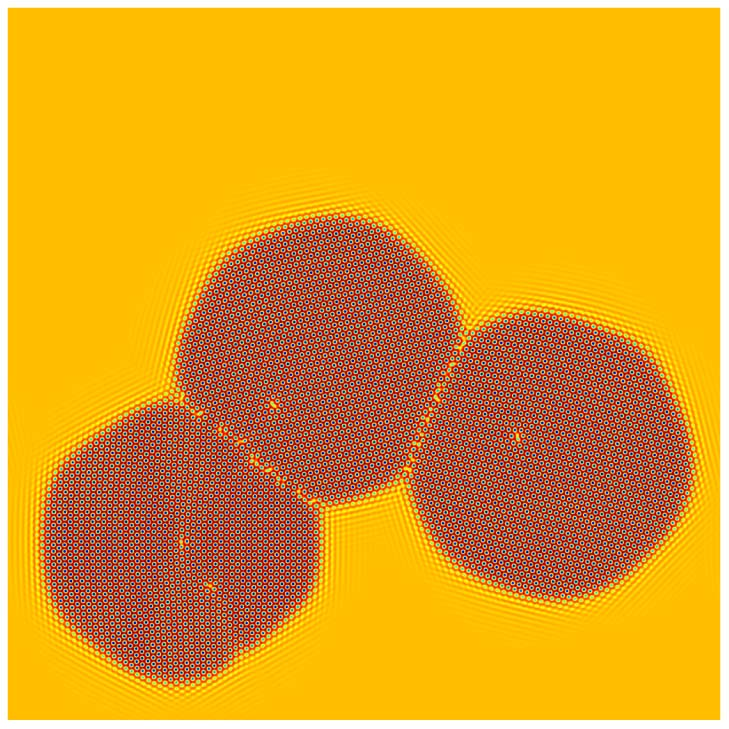} \includegraphics[scale=0.22]{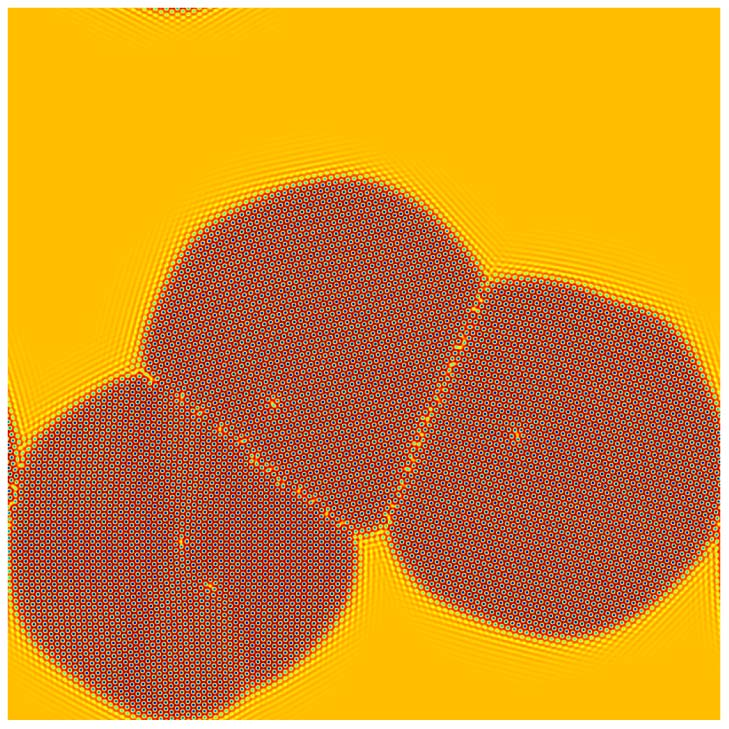}  \includegraphics[scale=0.22]{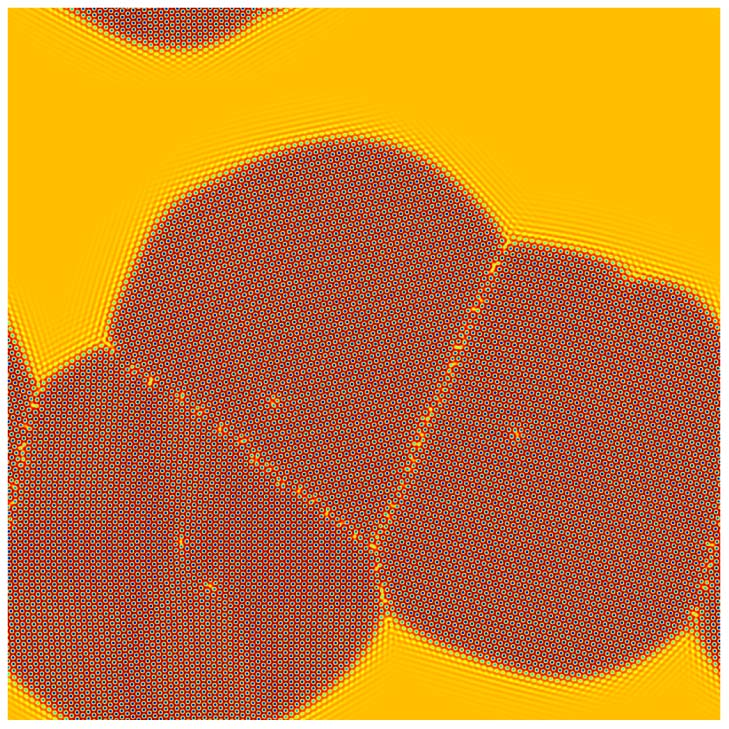}   \includegraphics[scale=0.22]{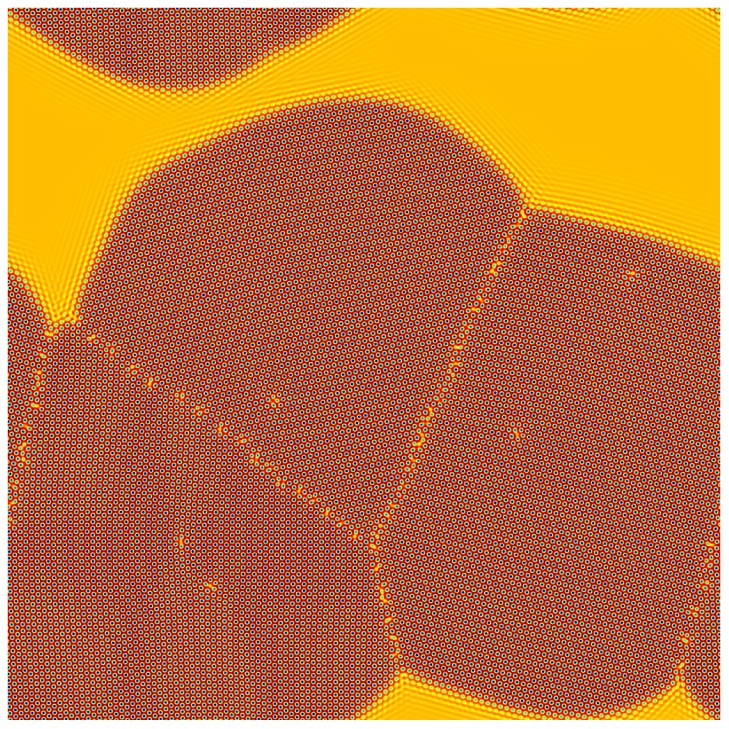}}
\centerline{ \includegraphics[scale=0.22]{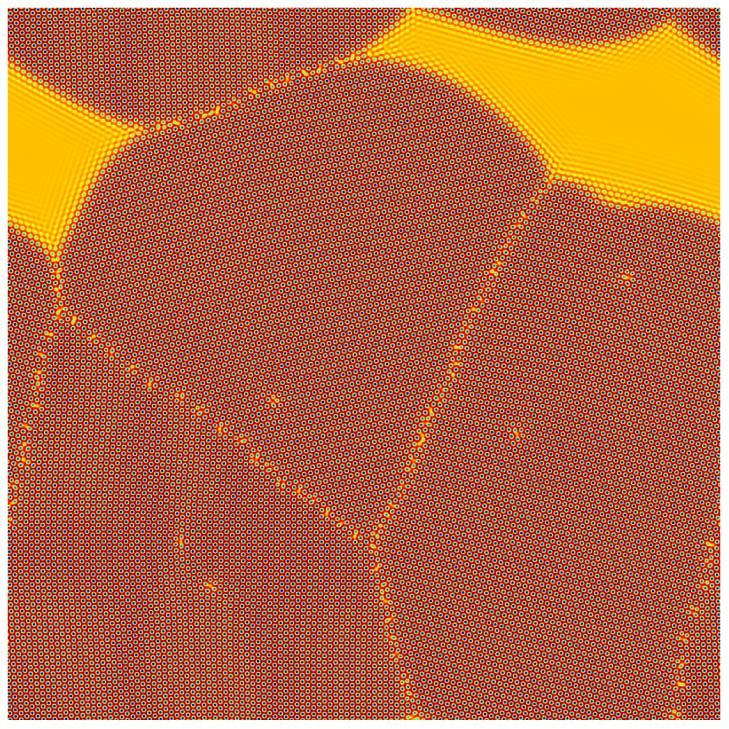} \includegraphics[scale=0.22]{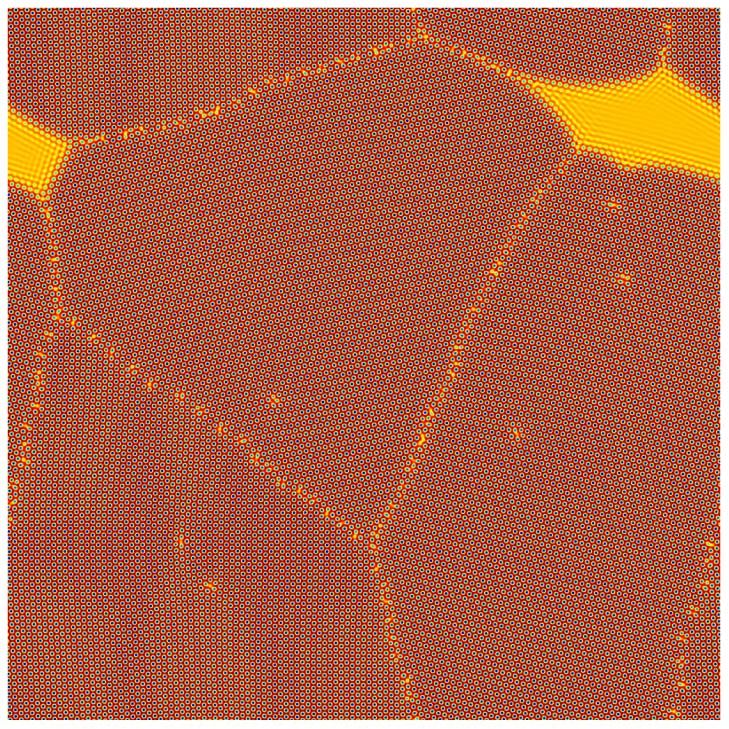}  \includegraphics[scale=0.22]{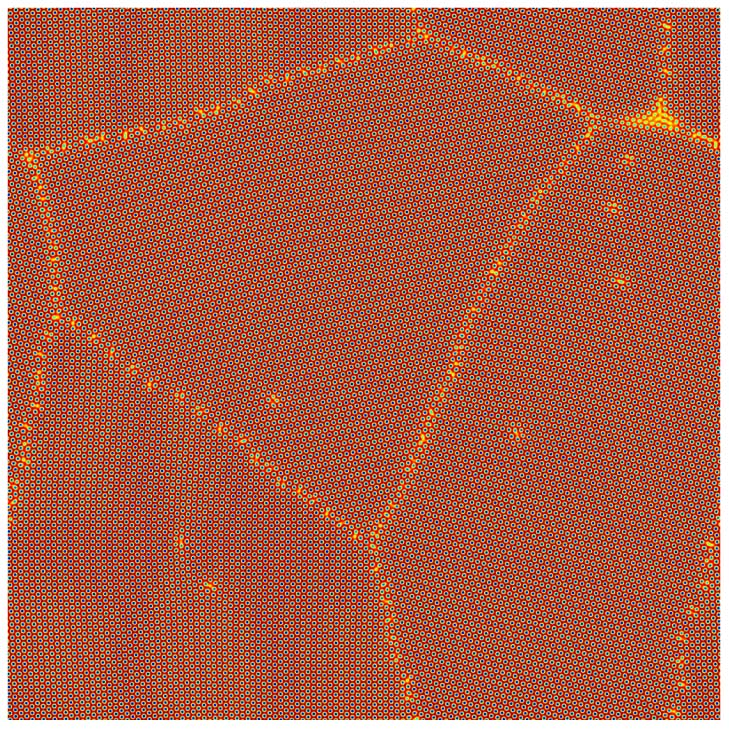}   \includegraphics[scale=0.22]{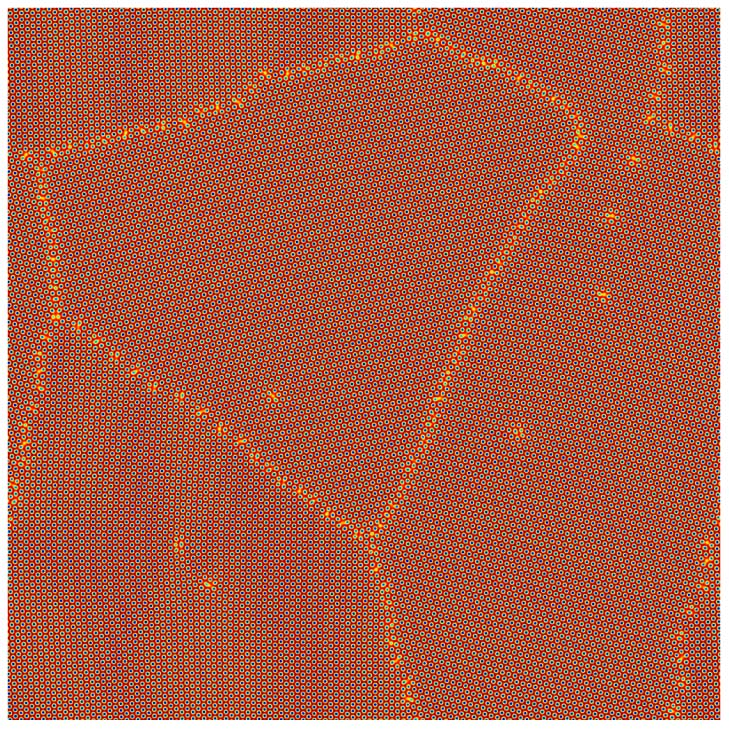}}
\vspace{-0.2cm}
\caption{Snapshots of the simulated growth of three crystallites  governed by the 2D PFC model around the times $t=0, 100, 200, 300, 400, 500, 600, 700, 800,$ $900, 1000$, and $2000$ respectively.
}\label{fig10}
\end{figure*}

Now we simulate the polycrystal growth governed by the PFC model \eqref{PFC} in a supercooled liquid.
Let us first consider the evolution of three crystallites with different orientations in two dimensions \cite{ZS22}, whose initial condition is given by
\beq
\phi(x_{l},y_{l},0)=0.285+0.446\Big(\cos\Big(\frac{0.66}{\sqrt{3}}y_{l}\Big)\cos(0.66x_{l})-0.5\cos\Big(\frac{1.32}{\sqrt{3}}y_{l}\Big)\Big), ~~l=1,2,3,
\eeq
where $x_l$ and $y_{l}$ define a local system of cartesian coordinates that is oriented with the crystallite lattice. The corresponding parameters are set to be $\sigma=1$ and $\delta=0.25$.
As shown in the first snapshot of Figure \ref{fig10}, the three initial crystallites in different orientations are located in the blocks with a length of $40$ in the domain $\Omega=(0,800)^{2},$  which are obtained via different definitions of the local coordinates $(x_{l},y_{l})$ using the following affine transformation of the global coordinates $(x,y)$:
\bryl
\begin{cases}
x_{l}=x\sin(\alpha)+y\cos(\alpha),\\
y_{l}=-x\cos(\alpha)+y\sin(\alpha),\\
\end{cases}
\eryl
where $\alpha=-\pi/4, 0, \pi/4$ for $l=1,2,3$, respectively. Some similar numerical examples  can  be found in \cite{Eld02,GN12,YH17}.
We adopt the TDSR-ETD3 scheme \eqref{full-ds-pfc} with $1024\times1024$ Fourier modes for spatial discretization
and the time-adaptive strategy \eqref{adp} with $\Dt_{min}=0.02, \Dt_{max}=1$ and $\gamma=10$. The simulation of the polycrystal growth lasts up to $T=2000$ in a supercooled liquid.
In Figure \ref{fig11}, we plot the evolutions of the simulated energy and the  time step sizes, which again  validate the stability and the efficiency of the proposed scheme \eqref{full-ds-pfc} with the time-adaptive strategy \eqref{adp}.
Figure \ref{fig10} displays  snapshots  of the crystal growth  around the times $t=0, 100, 200,300,400,500,600,700,800,900,1000$, and $2000$.
The numerical  results match very well with those reported in  \cite{ZS22},  which was produced by the relaxed SAV-BDF2 scheme with the fixed small time step size $\tau=0.02$.
The growth of the crystalline phase and the movement of the well-defined crystal-liquid interface are clearly observed, and
furthermore, it shows that different arrangements of crystallites may lead to defects and dislocations. Such crystal growth phenomenon has also been observed in \cite{Eld02,GN12,YH17}.

\begin{figure*}[!t]
\centerline{\includegraphics[scale=0.27]{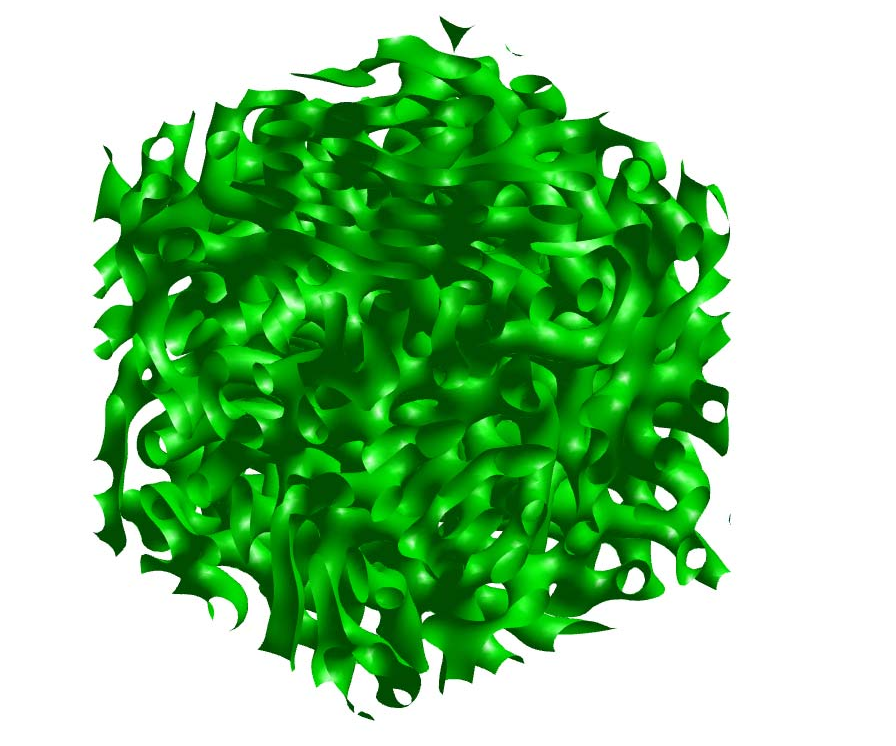} \includegraphics[scale=0.27]{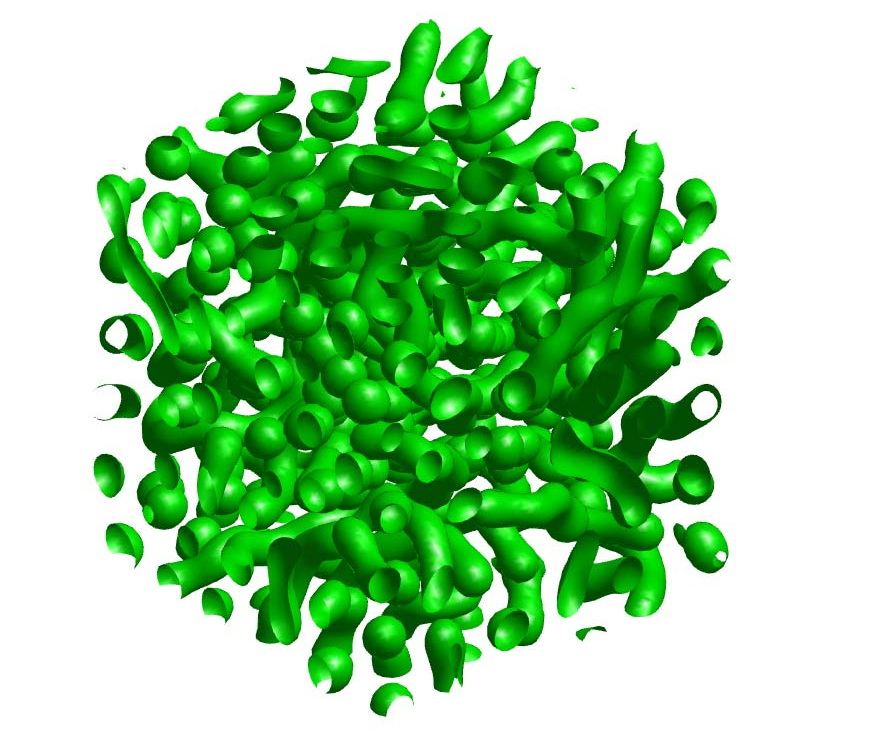} \includegraphics[scale=0.27]{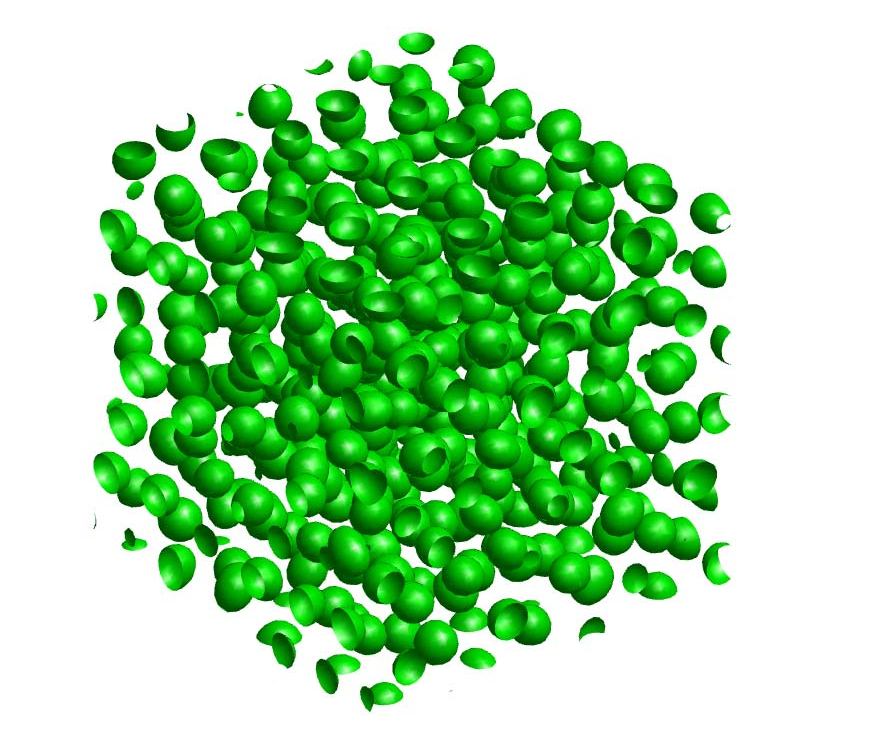}}
\centerline{\includegraphics[scale=0.25]{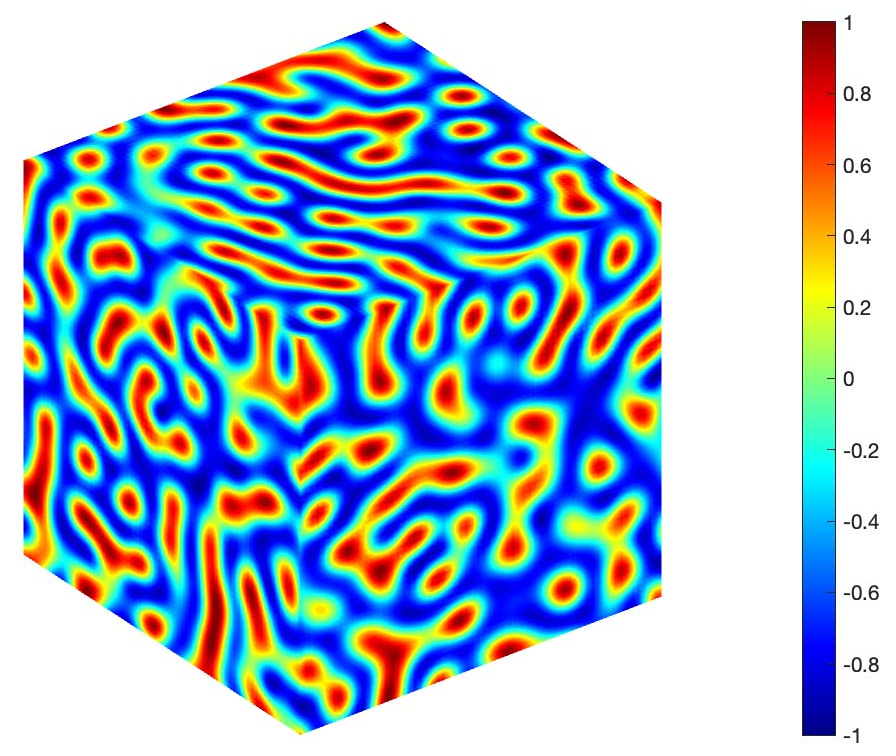}\includegraphics[scale=0.25]{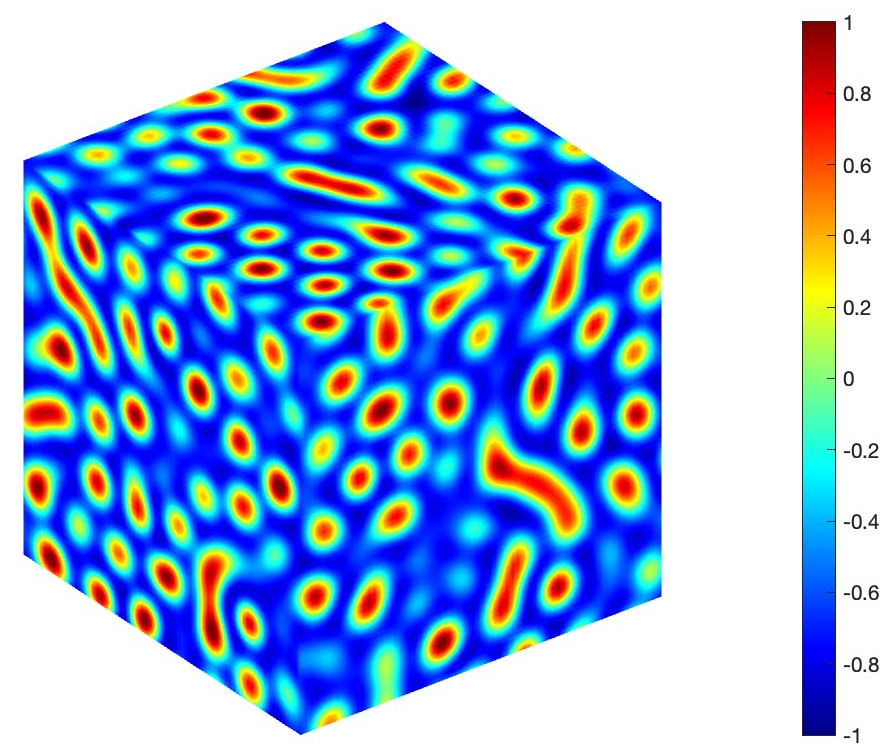} \includegraphics[scale=0.25]{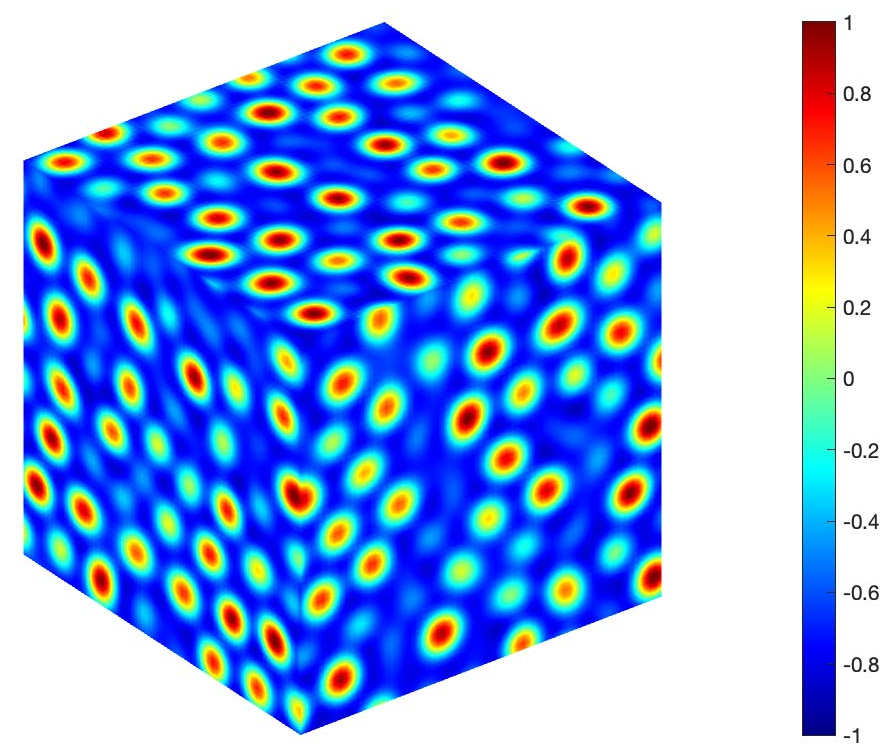}}
\vspace{-0.2cm}
\caption{Snapshots of the simulated  iso-surfaces of $\phi=0$  (top row) and   density field $\phi$ (bottom row) at $T=3000$ for the 3D PFC model with $\overline{\phi}=-0.2, -0.35$, and $-0.43$ (from left to right), respectively.}
\label{fig12}
\end{figure*}

Next we simulate the crystal growth in three dimensions with $\Omega=(0,50)^{3}$ and initial data $\phi(x,y,z,0)=\overline{\phi}+0.01 \text{rand}\,(\cdot)$ where $\overline{\phi}$ is a constant.
The parameters  are set to be $\sigma=1$ and $\delta=0.56$ and $64\times64\times64$ Fourier modes are applied to the spatial discretization in the  TDSR-ETD3 scheme \eqref{full-ds-pfc}. In Figure \ref{fig12},  snapshots of the iso-surfaces of the $\phi=0$ and the density field $\phi$ at the time $t=3000$ are displayed for the PFC model with different values of
$\overline{\phi}=-0.2,-0.35$, and $-0.43$. These results are consistent with those reported in \cite{GN12,PDACSGE07,ZS22}.

\section {Conclusions}
In this paper we have developed and analyzed a family of  high-order energy dissipative TDSR-ETD schemes for solving general gradient flows, utilizing a combination of the time-dependent spectral renormalization approach and the exponential integrator method. By incorporating the TDSR factor into the system, the discrete energy dissipation laws for the proposed schemes can be easily derived. In addition, the introduction of a new  enforcing term for updating  the TDSR factor significantly improves the time-step restriction to the numerical solution of  the resulting nonlinear systems. Extensive numerical experiments in two and three dimensions have been conducted to demonstrate the accuracy and efficiency of the proposed TDSR-ETD method, especially when adopting a time-adaptive strategy based on energy variation. It is also worth noting that the complexity and nonlinearity of the proposed schemes require nontrivial error analysis, which will require further deeper study.
 Moreover, the nonlinear term was treated implicitly in our proposed TDSR-ETD schemes, leading to the need of solving a nonlinear system at each time step. Therefore, it would also be  very interesting and important to investigate whether the proposed TDSR approach can be further extended to  construct  linear high-order schemes with similar energy dissipative properties.

%

\section{Acknowledgments}
We would like to express our sincere appreciation to Prof. Ziad Musslimani of Florida State University for his inspiring work and presentations on  the spectral renormalization method in physics, which stimulated the main idea behind our research.

\bibliographystyle{siamplain}
\bibliography{ref}
\end{document}

%% file: ex_shared.tex
\usepackage{lipsum}
\usepackage{amsfonts}
\usepackage{graphicx}
\usepackage{epstopdf}
\usepackage{algorithmic}
\usepackage{amssymb}
\ifpdf
  \DeclareGraphicsExtensions{.eps,.pdf,.png,.jpg}
\else
  \DeclareGraphicsExtensions{.eps}
\fi

\newsiamremark{remark}{Remark}
\newsiamremark{hypothesis}{Hypothesis}
\crefname{hypothesis}{Hypothesis}{Hypotheses}
\newsiamthm{claim}{Claim}

\headers{Spectral renormalization exponential integrator method}{D. Hou, L. Ju, and Z. Qiao}

\title{Energy-dissipative spectral renormalization exponential integrator method for \\gradient flow problems\thanks{
Submitted to the editors DATE.
\funding{The work of D. Hou's work is partially supported by NSFC grant 12001248, NSF of Jiangsu Province grant
BK20201020, NSF of Jiangsu Province Universities grant 20KJB110013 and the Hong Kong Polytechnic University grant 1-W00D.
L. Ju's work is partially supported by US National Science Foundation grant DMS-2109633.
Z. Qiao's work is partially supported by the Hong Kong Research Grants Council RFS grant RFS2021-5S03 and GRF grant 15302122, the Hong Kong Polytechnic University grant 4-ZZPF, and CAS AMSS-PolyU Joint Laboratory of Applied Mathematics.
}}}

\author{Dianming Hou\thanks{School of Mathematics and Statistics, Jiangsu Normal University, Xuzhou, Jiangsu 221116, China.
  (\email{dmhou@stu.xmu.edu.cn}).}
  \and Lili Ju\thanks{Department of Mathematics, University of South Carolina, Columbia, SC 29208, USA.
  (\email{ju@math.sc.edu}).}
\and Zhonghua Qiao\thanks{Department of Applied Mathematics, The Hong Kong Polytechnic University, Hung Hom, Kowloon, Hong Kong.
  (\email{zqiao@polyu.edu.hk}).}
}

\usepackage{amsopn}
